\documentclass[final,3p,times]{elsarticle}
\biboptions{sort&compress}
\usepackage{txfonts}
\usepackage{fleqn}

\usepackage{amsthm}
\usepackage{bm}
\usepackage{multirow}
\usepackage{array}
\usepackage{xcolor}
\colorlet{myrem}{red}
\usepackage{algorithm}% http://ctan.org/pkg/algorithms
\usepackage{algpseudocode}% http://ctan.org/pkg/algorithmicx

\usepackage[%
pdftex,%
colorlinks=true,%
bookmarks=true,%
citecolor=blue,%
urlcolor=blue]{hyperref} %pdflatex
\usepackage[caption=false]{subfig}

\usepackage{amsmath}
\usepackage{mathrsfs}  

\newcommand{\field}[1]{\mathbb{#1}}
\newcommand{\vbs}[1]{\boldsymbol{#1}}

\newcommand{\OP}[1]{\mathscr{#1}}

\newcommand{\tp}{\intercal}% transpose operation
\DeclareMathOperator*{\supp}{supp}
\let\Re\relax
\DeclareMathOperator{\Re}{Re}
\let\Im\relax
\DeclareMathOperator{\Im}{Im}

\DeclareMathOperator{\Span}{Span}

\graphicspath{{./figs/}}
\DeclareGraphicsExtensions{.pdf,.jpeg,.png}

\newcommand{\fs}[1]{\mathsf{#1}}

\DeclareMathOperator{\diag}{diag}

\newcommand{\ovl}[1]{\overline{#1}}
\newcommand{\melem}[1]{\mathfrak{#1}}

\let\Re\relax
\DeclareMathOperator{\Re}{Re}
\let\Im\relax
\DeclareMathOperator{\Im}{Im}
\newcommand{\vv}[1]{\boldsymbol{#1}}
\newcommand{\vs}[1]{\boldsymbol{#1}}

\DeclareMathOperator{\tvec}{vec}

\newcommand{\what}[1]{\widehat{#1}}

\newtheorem{rem}{Remark}
\usepackage{enumitem}

\journal{Communications in Nonlinear Science and Numerical Simulation}

\begin{document}

\begin{frontmatter}

%% Title, authors and addresses

%% use the tnoteref command within \title for footnotes;
%% use the tnotetext command for theassociated footnote;
%% use the fnref command within \author or \address for footnotes;
%% use the fntext command for theassociated footnote;
%% use the corref command within \author for corresponding author footnotes;
%% use the cortext command for theassociated footnote;
%% use the ead command for the email address,
%% and the form \ead[url] for the home page:
%% \title{Title\tnoteref{label1}}
%% \tnotetext[label1]{}
%% \author{Name\corref{cor1}\fnref{label2}}
%% \ead{email address}
%% \ead[url]{home page}
%% \fntext[label2]{}
%% \cortext[cor1]{}
%% \address{Address\fnref{label3}}
%% \fntext[label3]{}
\title{Transparent boundary condition and its high frequency approximation 
for the Schr\"{o}dinger equation on a rectangular computational domain}

%% use optional labels to link authors explicitly to addresses:
%% \author[label1,label2]{}
%% \address[label1]{}
%% \address[label2]{}

\author[phy]{Samardhi Yadav}
\ead{Samardhi@physics.iitd.ac.in}
\author[phy,opc]{Vishal Vaibhav}
\ead{vishal.vaibhav@gmail.com}

\address[phy]{Department of Physics, Indian Institute of Technology Delhi, Hauz
Khas, New Delhi–110016, India}%
\address[opc]{Optics and Photonics Center, Indian Institute of Technology Delhi,
Hauz Khas, New Delhi–110016, India}%
%=============================================================================%
%=============================================================================%
\begin{abstract}
This paper addresses the numerical implementation of the transparent boundary 
condition (TBC) and its various approximations for the free Schr\"{o}dinger equation 
on a rectangular computational domain. In particular, we consider the exact TBC 
and its spatially local approximation under high frequency assumption along with 
an appropriate corner condition. For the spatial discretization, we use a 
Legendre-Galerkin spectral method where Lobatto polynomials serve as the basis. 
Within variational formalism, we first arrive at the time-continuous dynamical 
system using spatially discrete form of the initial boundary-value problem 
incorporating the boundary conditions. This dynamical system is then discretized 
using various time-stepping methods, namely, the backward-differentiation formula 
of order 1 and 2 (i.e., BDF1 and BDF2, respectively) and the trapezoidal rule (TR) to obtain a 
fully discrete system. Next, we extend this approach to the novel 
Pad\'e based implementation of the TBC presented by Yadav and Vaibhav 
[\href{http://arxiv.org/abs/2403.07787} {\path{arXiv:2403.07787 (2024)}}]. 
Finally, several numerical tests are presented to demonstrate the effectiveness 
of the boundary maps (incorporating the corner conditions) where we study the 
stability and convergence behavior empirically. 

\end{abstract}
%=============================================================================%
%=============================================================================%


\end{frontmatter}
\tableofcontents

%=============================================================================%
%=============================================================================%
%=============================================================================%
%=============================================================================%

\section*{Notations}
\label{sec:notations}
The set of non-zero positive real numbers ($\field{R}$) is denoted by
$\field{R}_+$. For any complex number $\zeta$, $\Re(\zeta)$ and $\Im(\zeta)$ refer to 
the real and the imaginary parts of $\zeta$, respectively.
The open interval $(-1,1)$ is denoted by $\field{I}$.
 
%=============================================================================%
%=============================================================================%
%=============================================================================%
%=============================================================================%
\section{Introduction}
In this article, we address the numerical solution of the initial value 
problem (IVP) corresponding to the free Schr\"{o}dinger equation formulated on 
$\field{R}^2$ given by
\begin{equation}\label{eq:2D-SE}
i\partial_tu+\triangle u=0,\quad(\vv{x},t)\in\field{R}^2\times\field{R}_+,
\end{equation}
with initial condition $u(\vv{x},0)=u_0(\vv{x})$ which is assumed to be compactly 
supported. A Sommerfeld-like radiation condition at infinity is imposed in order 
to avoid any incoming waves which in turn ensures the uniqueness of the 
solution~\cite{CiCP2008}. We are interested in solving the IVP on a rectangular 
computational domain denoted by $\Omega_i$ such that 
$\supp u_0(\vv{x})\subset\Omega_i$. The Schr\"odinger 
equation appears in many areas of physical and technological interests, 
for instance, in quantum mechanics, underwater acoustics~\cite{LM1998}, 
and electromagnetic wave propagation~\cite{Levy2000}. This simple model can be 
augmented with a space and time-varying potential and nonlinear terms to model 
various systems of physical interest such as optical fibers~\cite{AK2003}. For the 
rectangular domain, development of an efficient numerical scheme 
for~\eqref{eq:2D-SE} which is capable of handling the boundary and corner 
conditions is key to addressing more general systems in this hierarchy.

%=============================================================================%
The exact TBC for the free Schr\"{o}dinger equation for convex computational 
domains with smooth boundaries have been developed by 
several authors~\cite{S2002,HH2004,S2005,V2019}. Among the existing approaches 
for a more general version of the Schr\"odinger equation, the pseudo-differential 
approach provides an approximation for the exact nonreflecting boundary 
operators~\cite{AB2001,ABM2004,S2005,ABK2012,ABK2013} for convex computational 
domains with smooth boundaries. Let us emphasize that the aforementioned 
methods are incapable of treating the corners in the boundary of the computational 
domain which rules out the rectangular domain. To the best of our knowledge, 
within the class of methods which employ nonreflecting boundary maps, the problem 
of corners has only been addressed in~\cite{Menza1997,FP2011,V2019,YV2024}. 
The transparent boundary operator for the Schr\"odinger equation on a 
rectangular domain has the form $\sqrt{\partial_t-i\triangle_{\Gamma}}$ and is 
nonlocal, both in time as well as space. The earliest approach for the numerical 
implementation of this operator was given by Menza~\cite{Menza1997} which is 
based on Pad\'e approximant for the square root function, however, this approach 
became problematic at the corners of the domain.
This problem was resolved by Feshchenko and Popov~\cite{FP2011} who developed 
a representation in terms of fractional operators. The numerical approach presented 
by Feschenko and Popov becomes computationally expensive on account of the 
increasing number of linear systems to be solved on the boundary segments 
with increasing time-steps. In the recent work of Yadav and Vaibhav~\cite{YV2024},
a novel-Pad\'e approach was developed which happens to be computationally efficient and
capable of handling the corners of the domain adequately. The generalization
of the novel-Pad\'e approach for more general systems (i.e., including exterior
potentials and nonlinear terms) is yet to be explored. For such systems, the 
use of pseudo-differential approach leads to approximate boundary 
conditions (BCs) involving the operators of the form 
$(\partial_t-i\triangle_{\Gamma})^{\alpha},\;\alpha=1/2,-1/2,-1,\ldots,$ 
where $\triangle_{\Gamma}$ is the Laplace-Beltrami operator. Implementation of 
such BCs is not yet addressed on a rectangular domain because of non-availability 
of the corner conditions which is an open problem. However, spatially local 
form of these BCs under high-frequency approximation readily offer the possibility 
of constructing the~\emph{corner conditions} for the case of free as well as 
general Schr\"odinger equation~\cite{V2019}. Let us emphasize the fact 
that the use of the high-frequency BCs may 
impact accuracy but they are easier to implement and affords comparatively 
lower computational complexity which makes them attractive\footnote{Note that 
use of certain type of time-stepping methods are only possible under periodic 
or Dirichlet BCs which require sufficiently large computational 
domain. In this work, we would like to explore if high-frequency BCs can offer 
a better trade-off between complexity and accuracy.}. 
%=============================================================================%

Having made the aforementioned observations, we now turn to the main objective 
of this work which is twofold: (i) We address the numerical solution of the IVP
in~\eqref{eq:2D-SE} with the BCs obtained as a result high-frequency approximation 
of the transparent boundary operator of this problem on a rectangular domain, and, 
(ii) we discuss the novel-Pad\'e based implementation of the TBCs which serves as 
golden standard for comparison\footnote{The method of perfectly matched layers 
will not be considered in this work for the purpose of comparison 
because it belongs to a different class of 
methods than presented in this article. Besides, a fair comparison can only be made if 
we consider an IVP that is fairly general so that no specific method has any 
distinct advantage. Therefore, this comparison is outside the scope of this article.}.
For space discretization, we use a Legendre-Galerkin method where Lobatto 
polynomials serve as basis functions. Within the variational formalism, we 
present a time-continuous route to arrive at a fully discrete system for the IBVP 
in~\eqref{eq:2D-SE-CT}. To this end, we start with a spatially discrete 
form of the initial boundary-value problem to obtain a time-continuous dynamical
system by incorporating the boundary conditions. This dynamical system is nonlocal 
on account of presence of time-fractional operators. Temporal discretization of this
system is addressed with convolution quadrature and also by employing
Pad\'e approximant based rational approximation for these time-fractional
operators. Note that the former remains nonlocal while the latter becomes effectively 
local. This formalism presents a straightforward way 
of arriving at a dynamical system which can be further studied using the theory 
of autonomous ODEs (discussed in~\ref{app:auto-sys}) to assess the stability of 
the overall system. Let us note that this new formalism presents the possibility of 
using any higher order time-stepping method for the temporal discretization which 
can improve the accuracy of overall system. Finally, we extend this approach to 
the novel-Pad\'e based implementation of TBCs within variational formalism.   

The paper is organized in the following manner: Sec.~\ref{sec:tbcs} presents 
high-frequency approximation of the TBCs followed by a discussion of numerically 
tractable forms of the TBCs. This section also addresses the development of corner 
conditions under high frequency approximation. In Sec.~\ref{sec:variational-prelim}, 
we collect all the facts about the Lobatto polynomial based Legendre-Galerkin spectral 
method for the numerical solution of the initial boundary-value problem (IBVP).
In Sec.~\ref{sec:vform-hf}, we discuss the numerical solution of the IBVP using
the boundary maps obtained as a result of high-frequency approximation
using Legendre-Galerkin method. In Sec.~\ref{sec:vform-np}, we discuss the numerical 
solution of the IBVP (in context of a Galerkin method) with effectively local form of 
the TBCs obtained as a result of novel Pad\'e approach. We test the efficiency of
our numerical schemes and confirm the order of convergence empirically with several 
numerical tests presented in Sec.~\ref{sec:numerical-experiments}. 
Finally, we conclude this paper in Sec.~\ref{sec:conclusion}.

%=============================================================================%
%=============================================================================%
%=============================================================================%
%=============================================================================%
\section{Transparent Boundary Conditions}\label{sec:tbcs}

\begin{figure}[!htbp]
\begin{center}
\def\myscale{1}
\includegraphics[scale=\myscale]{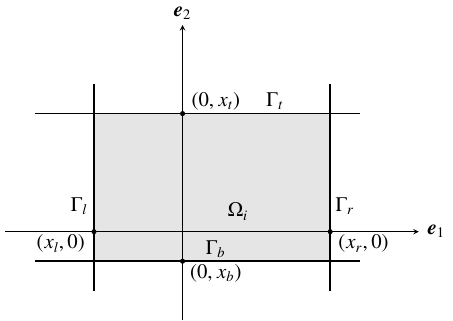}
\end{center}
\caption{\label{fig:rect-domain} The figure shows a rectangular domain
with boundary segments parallel to one of the axes.}
\end{figure}

Consider a rectangular computational domain ($\Omega_i$) with boundary 
segments parallel to one of the axes (see Fig.~\ref{fig:rect-domain}) with
$\Omega_i =(x_l,x_r)\times(x_b,x_t)$ referred to as the~\emph{interior} 
domain. Consider the decomposition of the field $u(\vv{x},t)$ such that
$u(\vv{x},t)\in\fs{L}^2(\field{R}^2)=\fs{L}^2(\Omega_i)\oplus \fs{L}^2(\Omega_e)$ 
where $\Omega_{e}=\field{R}^2\setminus\overline{\Omega}_i$ is referred to as 
the~\emph{exterior} domain.
%=============================================================================%
%=============================================================================%
An equivalent formulation of the IVP in~\eqref{eq:2D-SE} on the computational 
domain $\Omega_i$ can be stated as~\cite{V2019,YV2024}:
\begin{equation}\label{eq:2D-SE-CT}
\left\{\begin{aligned}
&i\partial_tu+\triangle u=0,\quad (\vv{x},t)\in\Omega_i\times\field{R}_+,\\
&u(\vv{x},0)=u_0(\vv{x})\in \fs{L}^2(\Omega_i),\quad\supp\,u_0\subset\Omega_i,\\
&\partial_{n}{u}+
e^{-i\pi/4}(\partial_{t}-i\partial^2_{x_2})^{1/2}u=0,
\quad\vv{x}\in\Gamma_l\cup\Gamma_r,\,t>0,\\
&\partial_{n}{u}+
e^{-i\pi/4}(\partial_{t}-i\partial^2_{x_1})^{1/2}u=0,
\quad\vv{x}\in\Gamma_b\cup\Gamma_t,\,t>0.
\end{aligned}\right.
\end{equation}
Recalling from~\cite{YV2024}, the operator $(\partial_{t}-i\partial^2_{x_1})^{1/2}$ 
is defined as
\begin{equation}
(\partial_t-i\partial_{x_2}^2)^{1/2}f
=(\partial_t-i\partial_{x_2}^2)[(\partial_t-i\partial_{x_2}^2)^{-1/2}f],
\end{equation}
where
\begin{equation}\label{eq:Op-ISP-TBC}
(\partial_t-i\partial_{x_2}^2)^{-1/2}f({x}_2,t)\\
=\frac{1}{\sqrt{\pi}}\int_0^t\int_{\field{R}}f(x'_2,\tau)
\frac{\mathcal{G}(x_2-x'_2,t-\tau)}{\sqrt{t-\tau}}dx_2'd\tau,
\end{equation}
with the kernel $\mathcal{G}$ function given by
\begin{equation}
\mathcal{G}(x_2,t)=\frac{e^{-i\pi/4}}{\sqrt{4\pi t}}
\exp\left[i\frac{x_2^2}{4t}\right].
\end{equation}
Let us note that the fractional operators present in the DtN maps 
are non-local both in space and time. For an efficient numerical implementation
of the IBVP~\eqref{eq:2D-SE-CT}, we present a spatially local form of these DtN-maps
by employing a high-frequency approximation with respect to temporal frequencies.
This approach is presented in Sec.~\ref{sec:CT-HF} which also includes corner
conditions obtained as a result of high-frequency approximation.
Before we address the numerical aspects, we revisit various 
numerically tractable representations for these operators at the continuous level 
in Sec.~\ref{sec:CT-Frac} and Sec.~\ref{sec:CT-NP}.
Finally, we present our main contribution from Sec.~\ref{sec:variational-prelim}
onward where a numerical recipe is described to arrive at a spatially discrete
time-continuous dynamical system in context of variational formulation.  
This dynamical system can be further discretized temporally with any choice 
of time-stepping method to arrive at a fully discrete system for the 
IBVP~\eqref{eq:2D-SE-CT}.
%=============================================================================%
%=============================================================================%
\subsection{High frequency approximation}\label{sec:CT-HF}
The boundary conditions in~\eqref{eq:2D-SE-CT} can be further simplified by making 
them local in terms of spatial variables under high-frequency approximation with
respect to the temporal frequencies. We reproduce the steps from~\cite{V2019} for 
the sake of completeness. Let $(\zeta_1,\zeta_2)$ be the covariables corresponding 
to $(x_1,x_2)$ to be used in the two dimensional Fourier transform. Let $z$ denote 
the complex variable in the Laplace transform. The operator 
$(\partial_{t}-i\partial^2_{x_1})^{-1/2}$ 
in the transformed $1/\sqrt{z+i\zeta_2^2}$ such that $\sqrt{\cdot}$ denotes the 
branch with $\Im(\alpha)>0$. Consider
\begin{equation}
\OP{L}^{-1}\left[\frac{1}{\sqrt{z+i\zeta_2^2}}\right]
=\frac{1}{2\pi i}\int \frac{e^{zt}}{\sqrt{z+i\zeta_2^2}}dz.
\end{equation}
Setting $\xi=zt$, we have
\begin{equation*}
\begin{split}
\OP{L}^{-1}\left[\frac{1}{\sqrt{z+i\zeta_2^2}}\right]
&=\frac{t^{-\frac{1}{2}}}{2\pi i}\int_{a+i\field{R}}\frac{e^{\xi}}{\sqrt{\xi+i\zeta^2_2t}}d\xi
=\frac{t^{-\frac{1}{2}}}{2\pi i}\int_{a+i\field{R}}
\left(1+\frac{i\zeta^2_2t}{\xi}\right)^{-\frac{1}{2}}\frac{e^{\xi}d\xi}{\sqrt{\xi}}\quad(a>0)\\
&\sim\frac{t^{-\frac{1}{2}}}{2\pi i}\int_{a+i\field{R}}
\left(\frac{1}{\xi^{\frac{1}{2}}}
-\frac{i\zeta^2_2t}{2\xi^{\frac{3}{2}}}
-\frac{3\zeta^4_2t^2}{8\xi^{\frac{5}{2}}}
+\ldots\right)e^{\xi}d\xi
\sim \frac{1}{\Gamma(\tfrac{1}{2})}t^{-\frac{1}{2}}
-\frac{i}{2\Gamma(\tfrac{3}{2})}t^{\frac{1}{2}}\zeta_2^2
-\frac{3}{8\Gamma(\tfrac{5}{2})}t^{\frac{3}{2}}\zeta_2^4
+\ldots
\end{split}
\end{equation*}
which yields the expansion
\begin{equation*}
\begin{split}
&(\partial_{t}-i\partial^2_{x_1})^{-1/2}
\sim \partial_t^{-1/2}
+i\frac{1}{2}\partial_t^{-3/2}\partial^2_{x_2}
-\frac{3}{8}\partial_t^{-5/2}\partial^4_{x_2}
+\ldots,\\
&(\partial_{t}-i\partial^2_{x_1})^{1/2}
=(\partial_{t}-i\partial^2_{x_1})^{-1/2}(\partial_{t}-i\partial^2_{x_1})
\sim \partial_t^{1/2}
-i\frac{1}{2}\partial_t^{-1/2}\partial^2_{x_2}
+\frac{1}{8}\partial_t^{-3/2}\partial^4_{x_2}
+\ldots
\end{split}
\end{equation*}
We obtain the following form of
DtN map as a result of high-frequency approximation on $\Gamma_r$:
\begin{equation}\label{eq:ABC-2D-HF}
\partial_{x_1}u+e^{-i\pi/4}\partial_t^{1/2}u
-e^{i\pi/4}\frac{1}{2}\partial_{x_2}^2\partial_t^{-1/2}u=0\mod{(\partial_t^{-3/2})}.
\end{equation}
Following a similar approach for rest of the segments of the rectangular domain, we 
obtain the following approximate BCs:
\begin{equation}\label{eq:ABC-2D-HF-rect}
\begin{split}
&\partial_{n}u+e^{-i\pi/4}\partial_t^{1/2}u
-e^{i\pi/4}\frac{1}{2}\partial_{x_2}^2\partial_t^{-1/2}u=0,
\quad\vv{x}\in\Gamma_r\cup\Gamma_l,\\
&\partial_{n}u+e^{-i\pi/4}\partial_t^{1/2}u
-e^{i\pi/4}\frac{1}{2}\partial_{x_1}^2\partial_t^{-1/2}u=0,
\quad\vv{x}\in\Gamma_b\cup\Gamma_t.
\end{split}
\end{equation}
These boundary conditions do not take into account the corners of the rectangular
domain and become problematic at corners. To understand this, we consider the weak formulation
of the original IVP as follows: Consider a test function $\psi(\vv{x})\in\fs{L}^2(\Omega_i)$, 
taking the inner product with the equation~\eqref{eq:2D-SE}, we have
\begin{equation}
\int_{\Omega_i}(i\partial_t u+\nabla^2 u)\psi d^2\vv{x}
=i\partial_t\int_{\Omega_i}u\psi d^2\vv{x}
-\int_{\Omega_i}(\nabla u)\cdot(\nabla\psi)d^2\vv{x}
+\int_{\partial\Omega_i}\psi(\nabla u)\cdot d\vbs{\sigma}.
\end{equation}
The boundary integrals on top and right segments are given by
\begin{equation*}
\begin{split}
\int_{\Gamma_r}\psi\partial_{{x_1}}u dx_2 +\int_{\Gamma_t}\psi\partial_{{x_2}}u dx_1
&=-e^{-i\pi/4}\int_{\Gamma_r\cup\Gamma_l}\psi\partial^{1/2}_tu
+\frac{1}{2}e^{i\pi/4}\left[\int_{\Gamma_r}\psi\partial^2_{{x_2}}\partial^{-1/2}_tu dx_2
+\int_{\Gamma_t}\psi\partial^2_{{x_1}}\partial^{-1/2}_tu dx_1\right]\\
&=-e^{-i\pi/4}\int_{\Gamma_r\cup\Gamma_l}\psi\partial^{1/2}_tu
+\frac{1}{2}e^{i\pi/4}\biggl[
\left.\psi\partial_{x_2}\partial^{-1/2}_tu\right|_{x_2=x_b}^{x_t}+
\left.\psi\partial_{x_1}\partial^{-1/2}_tu\right|_{x_1=x_l}^{x_r}\\
&\quad-\int_{\Gamma_r}(\partial_{x_2}\psi)(\partial_{{x_2}}\partial^{-1/2}_tu) dx_2
-\int_{\Gamma_t}(\partial_{x_1}\psi)(\partial_{{x_1}}\partial^{-1/2}_tu)dx_1\biggl].
\end{split}
\end{equation*}
Consider the terms which correspond to the top-right corner in the above
equation:
\begin{equation}
\left(\partial_{x_2}\partial^{-1/2}_tu
+\partial_{x_1}\partial^{-1/2}_tu\right)_{\Gamma_r\cap\Gamma_t}
=\partial_t^{-1/2}\left(\partial_{x_2}u+\partial_{x_1}u\right)_{\Gamma_r\cap\Gamma_t}.
\end{equation}
They are problematic on account of the fact that the BCs in the current 
form cannot be used to evaluate them. We can evaluate these terms by applying
$\partial^{-1/2}_t$ to the IVP in~\eqref{eq:2D-SE} as
\begin{equation}
i\partial_t^{1/2}u+(\partial^2_{x_1}+\partial^2_{x_2})\partial_t^{-1/2}u=0,
\,\,(x_1,x_2)\in\Gamma_r\cap\Gamma_t.
\end{equation}
Here, the fact that the field is zero at the corner at $t=0$ is explicitly used
to arrive at the fractional derivative. Using BCs in~\eqref{eq:ABC-2D-HF-rect} and
the last equation, we get the following corner condition:
\begin{equation}
\partial_{x_1}u+\partial_{x_2}u+\frac{3}{2}e^{-i\pi/4}\partial_t^{1/2}u=0,
\quad (x_1,x_2)\in\Gamma_r\cap\Gamma_t.
\end{equation}
Following a similar approach, we can construct the corner conditions for rest of the
corners of the rectangular domain:
\begin{equation}
\partial_{n}u|_{\Gamma_i}+\partial_{n}u|_{\Gamma_j}+\frac{3}{2}e^{-i\pi/4}\partial_t^{1/2}u=0,
\quad (x_1,x_2)\in \Gamma_i\cap\Gamma_j,
\end{equation}
where $i\neq j$ and $i,j\in\{r,t,l,b\}$. Let us introduce the set of corner points 
of the computational domain
as $\Gamma_C=\{\Gamma_{rt},\Gamma_{rb},\Gamma_{lt},\Gamma_{lb}\}$ where 
$\Gamma_{ij}=\Gamma_i\cap\Gamma_j$ with $i,j\in\{r,t,l,b\}$. Under high-frequency 
approximation, we have the following equivalent formulation 
for~\eqref{eq:2D-SE-CT} on $\Omega_i$: 
\begin{equation}\label{eq:2D-SE-HF}
\left\{\begin{aligned}
&i\partial_tu+\triangle u=0,\quad(\vv{x},t)\in\Omega_i\times\field{R}_+,\\
&u(\vv{x},0)=u_0(\vv{x})\in \fs{L}^2(\Omega_i),\quad\supp\,u_0\subset\Omega_i,\\
&\partial_{n}u+e^{-i\pi/4}\partial_t^{1/2}u
-e^{i\pi/4}\frac{1}{2}\triangle_{\Gamma}\partial_t^{-\frac{1}{2}}u=0,
\quad\vv{x}\in\Gamma\setminus\Gamma_C,\\
&\partial_{n_1}u+\partial_{n_2}u +\frac{3}{2}e^{-i\pi/4}\partial_t^{\frac{1}{2}}u=0,
\;\vv{x}\in\Gamma_C,\;t>0.
\end{aligned}\right.
\end{equation}

%=============================================================================%    
\begin{figure}[!hb]
\begin{center}
\def\myscale{0.75}
\includegraphics[scale=\myscale]{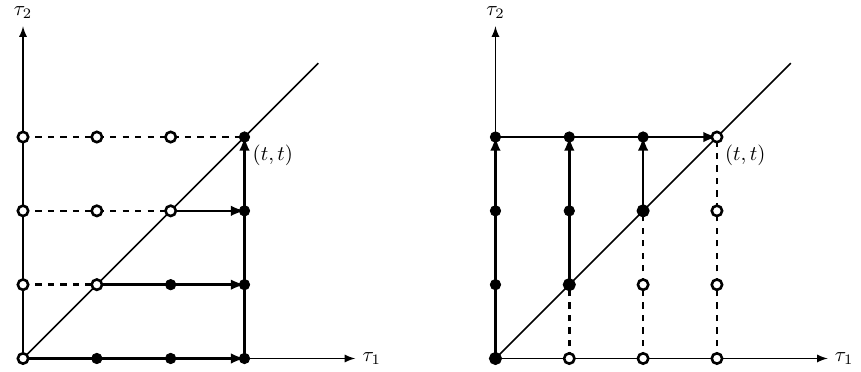}
\end{center}
\caption{\label{fig:IVP-auxi}A schematic depiction of the evolution of the 
auxiliary field $\varphi(x_1,x_2,\tau_1,\tau_2)$ in the 
$(\tau_1,\tau_2)$-plane is provided in this figure where the plot on the right 
corresponds $\vv{x}\in\Gamma_r\cup\Gamma_l$ and the plot on the left corresponds 
$\vv{x}\in\Gamma_t\cup\Gamma_b$. The filled circles depict the evolution of the 
auxiliary field $\varphi(x_1,x_2,\tau_1,\tau_2)$ either above or below the diagonal
in the $(\tau_1,\tau_2)$-plane starting from the diagonal which also serves as 
initial conditions for solving IVPs corresponding to the auxiliary function. 
The TBCs for the auxiliary field require the history of the auxiliary field 
at the corner points which makes empty circles relevant. Note that these
values at the corners can be taken from the adjacent segment of the boundary 
where it is already being computed and this is depicted by broken lines. 
Note that the vertical/horizontal lines where the arrows end corresponds to the
history of the auxiliary field needed for the TBCs on $\partial\Omega_i$ 
in the current time ($t$).}
\end{figure}
%=============================================================================%
%=============================================================================%    
%=============================================================================%    
\subsection{Fractional operator based formulation}\label{sec:CT-Frac}
In this section, we provide a summary of the fractional operator based approach 
which we presented in~\cite{YV2024}. The DtN maps in terms of a time-fractional 
derivative operator for the IVP~\eqref{eq:2D-SE-CT} can be expressed as
\begin{equation}\label{eq:maps-cq}
\begin{split}
&\partial_{x_1}{u}(\vv{x},t)=\left.-{e^{-i\pi/4}}\partial^{1/2}_{\tau_1}
\varphi(x_1,x_2,\tau_1,\tau_2)\right|_{\tau_1,\tau_2=t},\quad (x_1,x_2)\in\Gamma_r,\\
&\partial_{x_1}{u}(\vv{x},t)=\left.+{e^{-i\pi/4}}\partial^{1/2}_{\tau_1}
\varphi(x_1,x_2,\tau_1,\tau_2)\right|_{\tau_1,\tau_2=t},\quad (x_1,x_2)\in\Gamma_l,\\
&\partial_{x_2}{u}(\vv{x},t)=\left.-{e^{-i\pi/4}}\partial^{1/2}_{\tau_2}
\varphi(x_1,x_2,\tau_1,\tau_2)\right|_{\tau_1,\tau_2=t},\quad (x_1,x_2)\in\Gamma_t,\\
&\partial_{x_2}{u}(\vv{x},t)=\left.+{e^{-i\pi/4}}\partial^{1/2}_{\tau_2}
\varphi(x_1,x_2,\tau_1,\tau_2)\right|_{\tau_1,\tau_2=t},\quad (x_1,x_2)\in\Gamma_b.
\end{split}
\end{equation}
the auxiliary functions satisfy the IVPs given by
\begin{equation}\label{eq:ivps-cq}
\begin{split}
&[i\partial_{\tau_1}+\partial_{x_1}^2]\varphi(x_1,x_2,\tau_1,\tau_2)=0,
\quad (x_1,x_2)\in\Gamma_t\cup\Gamma_b,\;\tau_1\in(\tau_2,t],\\
&[i\partial_{\tau_2}+\partial_{x_2}^2]\varphi(x_1,x_2,\tau_1,\tau_2)=0,
\quad (x_1,x_2)\in\Gamma_l\cup\Gamma_r,\;\tau_2\in(\tau_1,t],
\end{split}
\end{equation}
with the initial condition $\varphi(x_1,x_2,\tau_1,\tau_1)=u(x_1,x_2,\tau_1)$ and 
$\varphi(x_1,x_2,\tau_2,\tau_2)=u(x_1,x_2,\tau_2)$. Let us introduce the set of corner 
points of the computational domain
as $\Gamma_C=\{\Gamma_{rt},\Gamma_{rb},\Gamma_{lt},\Gamma_{lb}\}$ where 
$\Gamma_{ij}=\Gamma_i\cap\Gamma_j$ with $i,j\in\{r,t,l,b\}$. The transparent boundary conditions 
at corner points for the IVPs listed in \eqref{eq:ivps-cq} are as follows
\begin{equation}\label{eq:maps-cq-auxi}
\begin{split}
&\partial_{x_1}\varphi(x_1,x_2,\tau_1,\tau_2)
+e^{-i\pi/4}\partial^{1/2}_{\tau_1}\varphi(x_1,x_2,\tau_1,\tau_2)=0,
\quad (x_1,x_2)\in\{\Gamma_{rt},\Gamma_{rb}\},\\
&\partial_{x_1}\varphi(x_1,x_2,\tau_1,\tau_2)
-e^{-i\pi/4}\partial^{1/2}_{\tau_1}\varphi(x_1,x_2,\tau_1,\tau_2)=0,
\quad (x_1,x_2)\in\{\Gamma_{lt},\Gamma_{lb}\},\\
&\partial_{x_2}\varphi(x_1,x_2,\tau_1,\tau_2)
+e^{-i\pi/4}\partial^{1/2}_{\tau_2}\varphi(x_1,x_2,\tau_1,\tau_2)=0,
\quad (x_1,x_2)\in\{\Gamma_{rt},\Gamma_{lt}\},\\
&\partial_{x_2}\varphi(x_1,x_2,\tau_1,\tau_2)
-e^{-i\pi/4}\partial^{1/2}_{\tau_2}\varphi(x_1,x_2,\tau_1,\tau_2)=0.
\quad (x_1,x_2)\in\{\Gamma_{rb},\Gamma_{lb}\}.\\
\end{split}
\end{equation}
The history required for the fractional derivatives present in the DtN maps
can be understood  with the help of the schematic shown in Fig.~\ref{fig:IVP-auxi}.
The two IVPs for the auxiliary field $\varphi(x_1,x_2,\tau_1,\tau_2)$ advances the 
field either above or below the diagonal in the $(\tau_1,\tau_2)$-plane starting 
from the diagonal which also serves as initial conditions for solving IVPs. This 
is depicted by filled circles where arrows denote the direction of evolution. 
The TBCs present in~\eqref{eq:maps-cq-auxi} to solve the IVPs for the auxiliary 
functions on the boundary segments require the history of the auxiliary function 
at corners from the start of the computations which makes the empty circles 
relevant. Note that these values at the corners can be taken from the adjacent 
segment of the boundary where it is already being computed. This is depicted by 
broken lines in Fig.~\ref{fig:IVP-auxi}. 

%=============================================================================%
%=============================================================================%
\subsection{Novel Pad\'e approach}\label{sec:CT-NP}
In this section, we provide a summary of the novel-Pad\'e approach 
which we presented in~\cite{YV2024}. This approach relies on the Pad\'e approximant 
based representation for the $1/2$-order temporal derivative in~\eqref{eq:maps-cq}. 
Let the $K$-th order diagonal Pad\'{e} approximant based rational approximation for the 
function $\sqrt{z}$ be denoted by $R_K^{(1/2)}(z)$ is given by
\begin{equation}\label{eq:sqrt-pade}
R^{(1/2)}_K(z)
=b_0-\sum_{k=1}^K\frac{b_k}{z+\eta^2_k},\quad\text{where}\quad 
\left\{\begin{aligned}
& b_0=2K+1,\;b_k = \frac{2\eta^2_k(1+\eta^2_k)}{2K+1},\\
&\eta_k = \tan\theta_k,\;\theta_k =\frac{k\pi}{2K+1},\;k=1,2,\ldots,K.
\end{aligned}\right.
\end{equation}
Following~\cite{YV2024}, we introduce the following notations for the auxiliary 
function:
\begin{equation}
\begin{split}
&\varphi_{a_1}(x_2,\tau_1,\tau_2)=\varphi(x_{a_1},x_2,\tau_1,\tau_2),\quad a_1\in\{l,r\},\\
&\varphi_{a_2}(x_1,\tau_1,\tau_2)=\varphi(x_1,x_{a_2},\tau_1,\tau_2),\quad a_2\in\{b,t\}.
\end{split}
\end{equation}
Corresponding to the each of the partial fractions in~\eqref{eq:sqrt-pade}, we
introduce the auxiliary functions $\varphi_{k,a_1}$ and 
$\varphi_{k,a_2}$ such that they satisfy the following ODEs: 
\begin{equation}\label{eq:ode-auxi-npade1}
\begin{split}
&(\partial_{\tau_1}+\eta^2_k)\varphi_{k,a_1}(x_2,\tau_1,\tau_2)=
\varphi_{a_1}(x_2,\tau_1,\tau_2),\\
&(\partial_{\tau_2}+\eta^2_k)\varphi_{k,a_2}(x_1,\tau_1,\tau_2)=
\varphi_{a_2}(x_1,\tau_1,\tau_2),
\end{split}
\end{equation}
with the initial conditions assumed to be $\varphi_{k,a_1}(x_2,0,\tau_2)=0$ and
$\varphi_{k,a_2}(x_1,\tau_1,0)=0$, respectively. The solution can be stated as 
\begin{equation}\label{eq:phi_k}
\begin{split}
& \varphi_{k,a_1}(x_2,\tau_1,\tau_2)
=\int_0^{\tau_1}e^{-\eta^2_k(\tau_1-s_1)}\varphi_{a_1}(x_2,s_1,\tau_2)ds_1,
\quad a_1\in\{l,r\},\\
& \varphi_{k,a_2}(x_1,\tau_1,\tau_2)
=\int_0^{\tau_2}e^{-\eta^2_k(\tau_2-s_2)}\varphi_{a_2}(x_1,\tau_1,s_2)ds_2,
\quad a_2\in\{t,b\}.
\end{split}
\end{equation}
The DtN maps for the interior problem present in~\eqref{eq:maps-cq} now reads as 
\begin{equation}\label{eq:maps-pade}
\begin{split}
&\partial_{x_1}{u}(\vv{x},t)\pm e^{-i\pi/4}\left[ b_0 u(\vv{x},t)
-\sum_{k=1}^K b_k\varphi_{k,a_1}(x_2,t,t)\right]\approx 0,\\
&\partial_{x_2}{u}(\vv{x},t)\pm e^{-i\pi/4}\left[ b_0 u(\vv{x},t)
-\sum_{k=1}^K b_k\varphi_{k,a_2}(x_1,t,t)\right]\approx 0.
\end{split}
\end{equation}
It is straightforward to conclude that auxiliary fields satisfy the same IVPs which are 
satisfied by auxiliary function
\begin{equation}\label{eq:ivp-pade-auxi}
\begin{split}
&[i\partial_{\tau_2}+\partial_{x_2}^2]\varphi_{k,a_1}(x_2,\tau_1,\tau_2)=0,
\quad (x_1,x_2)\in\Gamma_l\cup\Gamma_r,\;\tau_2\in(\tau_1,t],\\
&[i\partial_{\tau_1}+\partial_{x_1}^2]\varphi_{k,a_2}(x_1,\tau_1,\tau_2)=0,
\quad (x_1,x_2)\in\Gamma_t\cup\Gamma_b,\;\tau_1\in(\tau_2,t].
\end{split}
\end{equation}
The transparent boundary conditions for these IVPs can be written as:
\begin{equation}\label{eq:maps-pade-auxi}
\begin{split}
&\partial_{x_2}\varphi_{k,a_1}(x_b,\tau_1,\tau_2)
-e^{-i\pi/4}\partial^{1/2}_{\tau_2}\varphi_{k,a_1}(x_b,\tau_1,\tau_2)=0,\\
&\partial_{x_2}\varphi_{k,a_1}(x_t,\tau_1,\tau_2)
+e^{-i\pi/4}\partial^{1/2}_{\tau_2}\varphi_{k,a_1}(x_t,\tau_1,\tau_2)=0,\\
&\partial_{x_1}\varphi_{k,a_2}(x_l,\tau_1,\tau_2)
-e^{-i\pi/4}\partial^{1/2}_{\tau_1}\varphi_{k,a_2}(x_l,\tau_1,\tau_2)=0,\\
&\partial_{x_1}\varphi_{k,a_2}(x_r,\tau_1,\tau_2)
+e^{-i\pi/4}\partial^{1/2}_{\tau_1}\varphi_{k,a_2}(x_r,\tau_1,\tau_2)=0.
\end{split}
\end{equation}
Once again, we use the Pad\'e approximants based representation for the 
$1/2$-order temporal derivative operator present in the TBCs for 
the auxiliary fields described in~\eqref{eq:maps-pade-auxi}.
Introducing the auxiliary fields $\psi_{k,k',a_1,a_2}(\tau_1,\tau_2)$ and 
$\psi_{k,k',a_2,a_1}(\tau_1,\tau_2)$ at the end points of the boundary segments 
$\Gamma_{a_1}$ and $\Gamma_{a_2}$, respectively, such that 
In the physical space, every $\psi_{k,k',a_1,a_2}$ and $\psi_{k,k',a_2,a_1}$ satisfy
the following ODEs:
\begin{equation}\label{eq:ode-auxi-npade2}
\begin{split}
&(\partial_{\tau_2}+\eta^2_{k'})\psi_{k,k',a_1,a_2}(\tau_1,\tau_2)
=\varphi_{k,a_1}(x_{a_2},\tau_1,\tau_2),\\
& (\partial_{\tau_1}+\eta^2_{k'})\psi_{k,k',a_2,a_1}(\tau_1,\tau_2)
=\varphi_{k,a_2}(x_{a_1},\tau_1,\tau_2),
\end{split}
\end{equation}
with the initial conditions assumed to be $\psi_{k,k',a_1,a_2}(\tau_1,0)=0$ and
$\psi_{k,k',a_2,a_1}(0,\tau_2)=0$, respectively. The solution to the 
ODEs~\eqref{eq:ode-auxi-npade2} reads as
\begin{equation}
\begin{split}
\psi_{k,k',a_1,a_2}(\tau_1,\tau_2)
&=\int_0^{\tau_2}e^{-\eta^2_{k'}(\tau_2-s_2)}\varphi_{k,a_1}(x_{a_2},\tau_1,s_2)ds_2,\\
\psi_{k,k',a_2,a_1}(\tau_1,\tau_2)
&=\int_0^{\tau_1}e^{-\eta^2_{k'}(\tau_1-s_1)}\varphi_{k,a_2}(x_{a_1},s_1,\tau_2)ds_1.\\
\end{split}
\end{equation}
The DtN maps for the auxiliary fields~\eqref{eq:maps-pade-auxi} now reads as 
\begin{equation}\label{eq:dtn-using-psi}
\begin{split}
& \partial_{x_2}\varphi_{k,a_1}(x_{a_2},\tau_1,\tau_2)\pm e^{-i\pi/4}\left[
  b_0\varphi_{k,a_1}(x_{a_2},\tau_1,\tau_2) -
  \sum_{k'=1}^K b_{k'}\psi_{k,k',a_1,a_2}(\tau_1,\tau_2)\right]=0,\\
& \partial_{x_1}\varphi_{k,a_2}(x_{a_1},\tau_1,\tau_2)\pm e^{-i\pi/4}\left[ 
  b_0\varphi_{k,a_2}(x_{a_1},\tau_1,\tau_2) -
  \sum_{k'=1}^K b_{k'}\psi_{k,k',a_2,a_1}(\tau_1,\tau_2)\right]=0.
\end{split}
\end{equation}
It is interesting to note that fields $\psi_{k,k',a_1,a_2}(\tau_1,\tau_2)$ and 
$\psi_{k,k',a_2,a_1}(\tau_1,\tau_2)$ are transpose of each other~\cite{YV2024}. 
This observation allows us to develop a rather efficient numerical scheme from a 
storage point of view. These ODEs for the auxiliary fields can now be summarized as
\begin{equation}\label{eq:ode-auxi-npade3}
\begin{split}
&(\partial_{\tau_1}+\eta^2_{k})\psi_{k,k',a_1,a_2}(\tau_1,\tau_2)
=\varphi_{k',a_2}(x_{a_1},\tau_1,\tau_2),\\
&(\partial_{\tau_2}+\eta^2_{k})\psi_{k,k',a_2,a_1}(\tau_1,\tau_2)
=\varphi_{k',a_1}(x_{a_2},\tau_1,\tau_2).
\end{split}
\end{equation}
%=============================================================================%
\begin{figure}[!htbp]
\begin{center}
\includegraphics[scale=0.75]{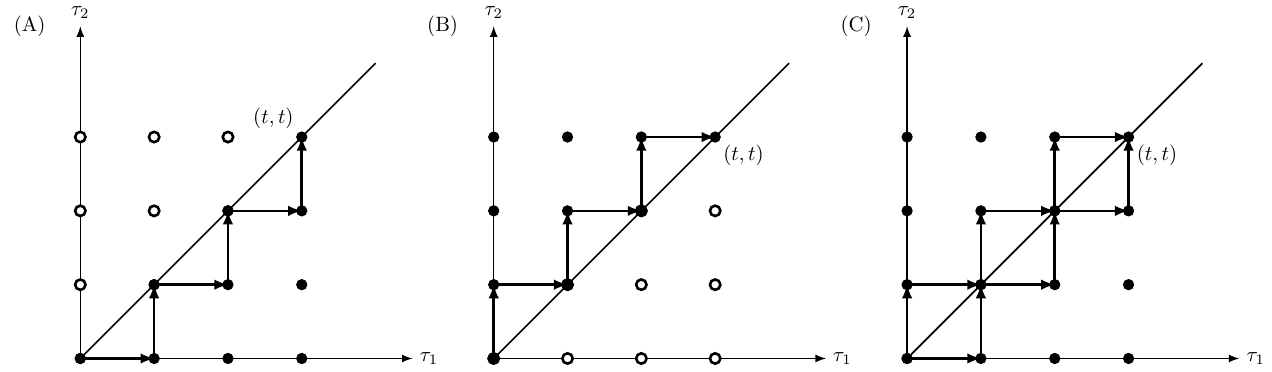}
\end{center}
\caption{\label{fig:IVP-auxi-pade}A schematic depiction of the evolution of the 
auxiliary fields $\varphi_{k,a_1}(x_2,\tau_1,\tau_2),\;
\varphi_{k,a_2}(x_1,\tau_1,\tau_2)$ and $\psi_{k,k',a_1,a_2}(\tau_1,\tau_2)$ in
the $(\tau_1,\tau_2)$-plane is provided in this figure. The plots (A)
and (B) depict the evolution of the fields $\varphi_{k,a_2}(x_1,\tau_1,\tau_2)$ 
and $\varphi_{k,a_1}(x_2,\tau_1,\tau_2)$ on the boundary segments $\Gamma_{a_2}$
and $\Gamma_{a_1}$, respectively. The plot (C) depicts the evolution of the field 
$\psi_{k,k',a_1,a_2}(\tau_1,\tau_2)$ which can be achieved by moving either below 
or above the diagonal.} 
\end{figure}
%=============================================================================%    
The evolution of the auxiliary fields $\varphi_{k,a_1}(x_2,\tau_1,\tau_2),\;
\varphi_{k,a_2}(x_1,\tau_1,\tau_2)$ and $\psi_{k,k',a_1,a_2}(\tau_1,\tau_2)$
in the $(\tau_1,\tau_2)$-plane can be understood with the help of the schematic 
shown in Fig.~\ref{fig:IVP-auxi-pade}. In this schematic, Fig.~\ref{fig:IVP-auxi-pade}(A)
and Fig.~\ref{fig:IVP-auxi-pade}(B) depict the evolution of the fields 
$\varphi_{k,a_2}(x_1,\tau_1,\tau_2)$ and $\varphi_{k,a_1}(x_2,\tau_1,\tau_2)$ 
on the boundary segments $\Gamma_{a_2}$
and $\Gamma_{a_1}$, respectively. The diagonal to diagonal computation 
of the fields $\varphi_{k,a_1}(x_2,\tau_1,\tau_2)$ and
$\varphi_{k,a_2}(x_1,\tau_1,\tau_2)$ consists of first advancing the fields using the IVPs
established in~\eqref{eq:ivp-pade-auxi} and then using the ODEs in~\eqref{eq:ode-auxi-npade1} 
for the second movement. Our novel Pad\'e approach makes the numerical scheme efficient 
from a storage point of view which is obvious from the schematic presented. Similarly,
Fig.~\ref{fig:IVP-auxi-pade}(C) depicts the evolution of the field 
$\psi_{k,k',a_1,a_2}(\tau_1,\tau_2)$ which can be achieved by moving either below or 
above the diagonal. 

%=============================================================================%    
%=============================================================================%    
%=============================================================================%    
\section{Variational Formulation: Preliminaries}\label{sec:variational-prelim}
In this section, we recall certain preliminary facts about the Lobatto
polynomials based Legendre-Galerkin spectral method for the numerical solution 
of the initial boundary-value problem (IBVP) stated in~\eqref{eq:2D-SE-CT}. For 
the computational domain $\Omega_i$, we introduce a reference domain 
$\Omega_i^{\text{ref.}}=\field{I}\times\field{I}$ where $\field{I}=(-1,1)$. In 
order to describe the associated linear maps between the reference domain and the 
actual computational domain, we introduce the variables 
$y_1,y_2\in\Omega_i^{\text{ref.}}$ such that
\begin{equation}
\left\{\begin{aligned}
&x_1 = J_1 y_1+\bar{x}_1,\quad J_1 = \frac{1}{2}(x_r-x_l),\quad
\bar{x}_1=\frac{1}{2}(x_l+x_r), & \beta_1=J_1^{-2},\\
&x_2 = J_2 y_2+\bar{x}_2,\quad J_2 = \frac{1}{2}(x_t-x_b),\quad
\bar{x}_2=\frac{1}{2}(x_b+x_t), & \beta_2=J_2^{-2}.
\end{aligned}\right.
\end{equation}
Let $L_n(y)$ denotes the Legendre polynomial of degree $n$, then  
the Lobatto polynomials can be defined as
\begin{equation}
\begin{split}
&\phi_0(y)=\frac{1}{2}[L_0(y)-L_1(y)],\quad \phi_1(y)=\frac{1}{2}[L_0(y)+L_1(y)],\\
&\phi_k(y)
=\frac{1}{\sqrt{2(2k-1)}}\left[L_k(y)-L_{k-2}(y)\right]
=c_k\left[L_k(y)-L_{k-2}(y)\right],\quad k\geq 2.
\end{split}
\end{equation}
Observe that $\phi_k(\pm1)=0,\,k\geq2$. Introduce the polynomial space 
\begin{equation}
\fs{P}_N = \Span\left\{\phi_p(y)|\;p=0,1,\ldots,N,\;y\in\ovl{\field{I}}\right\}.
\end{equation}
Let the index set $\{0,1,\ldots,N_p\}$ be denoted by $\field{J}_p$ where $p=1,2$. 
The mass matrix and the stiffness matrix for the Lobatto basis are
denoted by $M_p = (\melem{m}_{kj})_{k,j\in\field{J}_p}$ and 
$S_p = (\melem{s}_{kj})_{k,j\in\field{J}_p}$ for $p=1,2$, respectively.
For $j,k\geq 2$, the matrix entries are given by
\begin{equation}\label{eq:sys-mat-1d}
\melem{s}_{jk}=(\phi'_j,\phi'_k)_{\field{I}}
=\begin{cases}
1,&j=k,\\
0,&\mbox{otherwise}
\end{cases},\quad
\melem{m}_{jk}=(\phi_j,\phi_k)_{\field{I}}
=\begin{cases}
-\frac{1}{\sqrt{(2k-1)(2k-5)}},&j=k-2,\\
\frac{2}{(2k+1)(2k-3)},&j=k,\\
-\frac{1}{\sqrt{(2k-1)(2k+3)}},&j=k+2,\\
0&\mbox{otherwise}.
\end{cases}
\end{equation}
The remaining matrix entries are
\begin{equation}
\begin{split} 
&\melem{s}_{00}=\melem{s}_{11}=1/2,\quad\melem{s}_{01}=\melem{s}_{10}=-1/2,\\
&\melem{m}_{00}=\melem{m}_{11}=2/3,\quad\melem{m}_{01}=\melem{m}_{10}=1/3.
\end{split}
\end{equation}
For the 2D problem, we consider the tensor product space 
$\fs{P}_{N_1}\otimes\fs{P}_{N_2}$ of the Lobatto polynomial spaces
$\fs{P}_{N_j},\,j=1,2$. The basis for the tensor product space 
$\fs{P}_{N_1}\otimes\fs{P}_{N_2}$ can be written as
\begin{equation}
\theta_{p_1,p_2}(y_1,y_2)=\phi_{p_1}(y_1)\phi_{p_2}(y_2),\quad
p_1\in\field{J}_1,\;p_2\in\field{J}_2.
\end{equation}
Within the variational formulation, our goal is to find the approximate solution 
$u_{N_1,N_2}\in\fs{P}_{N_1}\otimes\fs{P}_{N_2}$ for the interior problem in~\eqref{eq:2D-SE-CT}:
\begin{equation}
\left(i\partial_tu_{N_1,N_2},\theta_{p_1,p_2}\right)_{\Omega_i^{\text{ref}}} 
+\left(J_1^{-2}\partial^2_{y_1}u_{N_1,N_2},\theta_{p_1,p_2}\right)_{\Omega_i^{\text{ref}}}
+\left(J_2^{-2}\partial^2_{y_2}u_{N_1,N_2},\theta_{p_1,p_2}\right)_{\Omega_i^{\text{ref}}}=0,
\end{equation}
where $(u,\theta)_{\Omega_i^{\text{ref}}}=\int_{\Omega_i^{\text{ref}}}u\theta d^2\vv{y}$ is the 
scalar product in $\fs{L}^2(\Omega_i^{\text{ref}})$.
The variational formulation defined above can be restated in terms of
field $u$ (dropping the subscripts `$N_1$' and '$N_2$' for the sake of brevity) as
\begin{multline}\label{eq:vform-2d}
\left(i\partial_t u,\theta_{p_1,p_2}\right)_{\Omega_i^{\text{ref}}} 
+\left(J_1^{-2}\partial^2_{y_1}u,\theta_{p_1,p_2}\right)_{\Omega_i^{\text{ref}}}
+\left(J_2^{-2}\partial^2_{y_2}u,\theta_{p_1,p_2}\right)_{\Omega_i^{\text{ref}}}\\
=i\partial_t\left(u,\theta_{p_1,p_2}\right)_{\Omega_i^{\text{ref}}}
-J_1^{-2}\left(\partial_{y_1}u,\phi'_{p_1}\phi_{p_2}\right)_{\Omega_i^{\text{ref}}}
-J_2^{-2}\left(\partial_{y_2}u,\phi_{p_1}\phi'_{p_2}\right)_{\Omega_i^{\text{ref}}}
+\mathcal{I}_1+\mathcal{I}_2,
\end{multline}
where the integrals labelled $\mathcal{I}_1$ and $\mathcal{I}_2$ are
\begin{equation}\label{eq:boundary-integrals}
\begin{split}
\mathcal{I}_1&=J_1^{-2}\int_{\field{I}}\left(
[\theta_{p_1,p_2}\partial_{y_1}u]_{+1}-[\theta_{p_1,p_2}\partial_{y_1}u]_{-1}\right)dy_2,\\
\mathcal{I}_2&=J_2^{-2}\int_{\field{I}}\left(
[\theta_{p_1,p_2}\partial_{y_2}u]_{+1}-[\theta_{p_1,p_2}\partial_{y_2}u]_{-1}\right)dy_1.
\end{split}
\end{equation}
These terms correspond to the boundary segments so that the boundary conditions
can be used to arrive at the complete variational form for the IBVP. Let
$\what{u}_{p_1,p_2}$ denote the expansion coefficients in the Lobatto basis so
that
\begin{equation}
u_{N_1,N_2}(y_1,y_2,t)=\sum_{p_1\in\field{J}_1}\sum_{p_2\in\field{J}_2}\what{u}_{p_1,p_2}
\phi_{p_1}(y_1)\phi_{p_2}(y_2).
\end{equation}
Let us introduce the matrix 
$\what{U}=\left(\what{u}_{p_1,p_2}\right)_{p_1\in\field{J}_1,\;p_2\in\field{J}_2}$
to represent the field in the Lobatto basis on $\Omega_{\text{ref.}}$. Further, let 
$\Lambda_p=\diag(1,1,0,\ldots,0)\in\field{R}^{(N_p+1)\times(N_p+1)},\;p=1,2$, then the
representation of the restriction of the field on the boundary segments is given by
\begin{equation}\label{eq:restrict-seg}
\begin{split}
\left.\begin{aligned}
&u_l=u|_{y_1=-1}=\sum_{p_2\in\field{J}_2}\hat{u}_{0,p_2}(t)\phi_{p_2}(y_2),\\
&u_r=u|_{y_1=+1}=\sum_{p_2\in\field{J}_2}\hat{u}_{1,p_2}(t)\phi_{p_2}(y_2)
\end{aligned}\right\}\rightarrow\Lambda_1\what{U}(t);\quad
\left.\begin{aligned}
&u_b=u|_{y_2=-1}=\sum_{p_1\in\field{J}_1}\hat{u}_{p_1,0}(t)\phi_{p_1}(y_1),\\ 
&u_t=u|_{y_2=+1}=\sum_{p_1\in\field{J}_1}\hat{u}_{p_1,1}(t)\phi_{p_1}(y_1)
\end{aligned}\right\}\rightarrow\what{U}(t)\Lambda_2.
\end{split}
\end{equation}
For convenience, we introduce the vectors:
\begin{equation}
\begin{split}
& \hat{\vv{u}}_l = (\hat{u}_{0,1},\hat{u}_{0,2},\ldots,\hat{u}_{0,N_2})^{\tp},\quad
  \hat{\vv{u}}_b = (\hat{u}_{1,0},\hat{u}_{2,0},\ldots,\hat{u}_{N_1,0})^{\tp},\\
& \hat{\vv{u}}_r = (\hat{u}_{1,1},\hat{u}_{1,2},\ldots,\hat{u}_{0,N_2})^{\tp},\quad
  \hat{\vv{u}}_t = (\hat{u}_{1,1},\hat{u}_{2,1},\ldots,\hat{u}_{N_1,1})^{\tp}.
\end{split}
\end{equation}
This is made possible on account of the observation: $\phi_0(+1)=0$, $\phi_0(-1)=1$,
$\phi_1(-1)=0$ and $\phi_1(+1)=1$. The field at the corners of
$\Omega_{\text{ref.}}$ can be represented by
\begin{equation}
\left.\begin{aligned}
&u|_{y_1=-1, y_2=-1}=\hat{u}_{0,0}(t),\quad u|_{y_1=+1, y_2=-1}=\hat{u}_{1,0}(t),\\
&u|_{y_1=-1, y_2=+1}=\hat{u}_{0,1}(t),\quad u|_{y_1=+1, y_2=+1}=\hat{u}_{1,1}(t),
\end{aligned}\right\}\rightarrow\Lambda_1\what{U}(t)\Lambda_2.
\end{equation}

%=============================================================================%    
%=============================================================================%    
\section{Variational Formulation: HF}\label{sec:vform-hf}
In this section, we consider the numerical solution of the IBVP 
in~\eqref{eq:2D-SE-HF} using a Legendre-Galerkin method as laid out 
in Sec.~\ref{sec:variational-prelim}. The variational formulation begins 
with~\eqref{eq:vform-2d} where the boundary terms, namely, 
$\mathcal{I}_1$ and $\mathcal{I}_2$, need to incorporate the boundary 
conditions in~\eqref{eq:2D-SE-HF}.  
Starting the with the segments denoted 
by $\Gamma_{a_1},\;a_1\in\{l,r\}$, we have
\begin{equation}\label{eq:hf-i-one}
\begin{split}
\mathcal{I}_1&=J_1^{-2}\int_{\field{I}}\left(
[\theta_{p_1,p_2}\partial_{y_1}u]_{+1}-[\theta_{p_1,p_2}\partial_{y_1}u]_{-1}\right)dy_2\\
&=J_1^{-1}\int_{\field{I}}\phi_{p_1}(+1)\phi_{p_2}(y_2)\left(-e^{-i\pi/4}\partial_t^{1/2}u
+e^{i\pi/4}\frac{1}{2}J^{-2}_2\partial_{y_2}^2\partial_t^{-1/2}u\right)_{y_1=+1}dy_2\\
&\quad+J_1^{-1}\int_{\field{I}}\phi_{p_1}(-1)\phi_{p_2}(y_2)\left(-e^{-i\pi/4}\partial_t^{1/2}u
+e^{i\pi/4}\frac{1}{2}J^{-2}_2\partial_{y_2}^2\partial_t^{-1/2}u\right)_{y_1=-1}dy_2\\
&=-e^{-i\pi/4}J_1^{-1}\partial_t^{1/2}\left[\phi_{p_1}(+1)\int_{\field{I}}\phi_{p_2}(y_2)
\left(u\right)_{y_1=+1}dy_2
+\phi_{p_1}(-1)\int_{\field{I}}\phi_{p_2}(y_2)\left(u\right)_{y_1=-1}dy_2\right]\\
&\quad+e^{i\pi/4}\frac{1}{2}(J_1)^{-1}(J_2)^{-2}\partial_t^{-1/2}
\biggl[\underbrace{\phi_{p_1}(+1)\int_{\field{I}}\phi_{p_2}(y_2)
\left(\partial_{y_2}^2u\right)_{y_1=+1}dy_2}_{\mathcal{I}_{1,r}}
+\underbrace{\phi_{p_1}(-1)\int_{\field{I}}\phi_{p_2}(y_2)
\left(\partial_{y_2}^2u\right)_{y_1=-1}dy_2}_{\mathcal{I}_{1,l}}\biggl].
\end{split}
\end{equation}
Similarly, for the segments denoted by $\Gamma_{a_2},\;a_2\in\{b,t\}$, we have
\begin{equation}\label{eq:hf-i-two}
\begin{split}
\mathcal{I}_2&=J_2^{-2}\int_{\field{I}}\left(
[\theta_{p_1,p_2}\partial_{y_2}u]_{+1}-[\theta_{p_1,p_2}\partial_{y_2}u]_{-1}\right)dy_1\\
&=J_2^{-1}\int_{\field{I}}\phi_{p_1}(y_1)\phi_{p_2}(+1)\left(-e^{-i\pi/4}\partial_t^{1/2}u
+e^{i\pi/4}\frac{1}{2}J^{-2}_1\partial_{y_1}^2\partial_t^{-1/2}u\right)_{y_2=+1}dy_1\\
&\quad+J_2^{-1}\int_{\field{I}}\phi_{p_1}(y_1)\phi_{p_2}(-1)\left(-e^{-i\pi/4}\partial_t^{1/2}u
+e^{i\pi/4}\frac{1}{2}J^{-2}_1\partial_{y_1}^2\partial_t^{-1/2}u\right)_{y_2=-1}dy_1\\
&=-e^{-i\pi/4}J_2^{-1}\partial_t^{1/2}\left[\phi_{p_2}(+1)\int_{\field{I}}\phi_{p_1}(y_1)
\left(u\right)_{y_2=+1}dy_1
+\phi_{p_2}(-1)\int_{\field{I}}\phi_{p_1}(y_1)\left(u\right)_{y_2=-1}dy_1\right]\\
&\quad+e^{i\pi/4}\frac{1}{2}(J_1)^{-2}(J_2)^{-1}\partial_t^{-1/2}
\biggl[\underbrace{\phi_{p_2}(+1)\int_{\field{I}}\phi_{p_1}(y_1)
\left(\partial_{y_1}^2u\right)_{y_2=+1}dy_1}_{\mathcal{I}_{2,t}}
+\underbrace{\phi_{p_2}(-1)\int_{\field{I}}\phi_{p_1}(y_1)
\left(\partial_{y_1}^2u\right)_{y_2=-1}dy_1}_{\mathcal{I}_{2,b}}\biggl].
\end{split}
\end{equation}
The terms involving second order derivatives acting on the 
field $u$ in~\eqref{eq:hf-i-one} and~\eqref{eq:hf-i-two} need to be further simplified 
using integration by parts and invoking the corner conditions collectively. Consider the 
terms $\mathcal{I}_{1,a_1}$ and $\mathcal{I}_{2,a_2}$, we have
\begin{equation}
\begin{split}
\mathcal{I}_{1,r}
&=\phi_{p_1}(+1)\left[\phi_{p_2}(+1)\left(\partial_{y_2}u\right)_{y_1=+1,y_2=+1}
-\phi_{p_2}(-1)\left(\partial_{y_2}u\right)_{y_1=+1,y_2=-1}
-\int_{\field{I}}\partial_{y_2}\phi_{p_2}(y_2)\left(\partial_{y_2}u\right)_{y_1=+1}dy_2\right],\\
\mathcal{I}_{1,l}&=\phi_{p_1}(-1)\left[\phi_{p_2}(+1)\left(\partial_{y_2}u\right)_{y_1=-1,y_2=+1}
-\phi_{p_2}(-1)\left(\partial_{y_2}u\right)_{y_1=-1,y_2=-1}
-\int_{\field{I}}\partial_{y_2}\phi_{p_2}(y_2)
\left(\partial_{y_2}u\right)_{y_1=-1}dy_2\right],\\
\mathcal{I}_{2,t}
&=\phi_{p_2}(+1)\left[\phi_{p_1}(+1)\left(\partial_{y_1}u\right)_{y_1=+1,y_2=+1}
-\phi_{p_1}(-1)\left(\partial_{y_1}u\right)_{y_1=-1,y_2=+1}
-\int_{\field{I}}\partial_{y_1}\phi_{p_1}(y_1)\left(\partial_{y_1}u\right)_{y_2=+1}dy_1\right],\\
\mathcal{I}_{2,b}&=\phi_{p_2}(-1)\left[\phi_{p_1}(+1)\left(\partial_{y_1}u\right)_{y_1=+1,y_2=-1}
-\phi_{p_1}(-1)\left(\partial_{y_1}u\right)_{y_1=-1,y_2=-1}
-\int_{\field{I}}\partial_{y_1}\phi_{p_1}(y_1)
\left(\partial_{y_1}u\right)_{y_2=-1}dy_1\right].
\end{split}
\end{equation}
Collecting the corner terms from each of the expressions $\mathcal{I}_{1,a_1}$ and $\mathcal{I}_{2,a_2}$ 
and invoking the corner conditions from~\eqref{eq:2D-SE-HF}, we obtain
\begin{equation}
\begin{split}
&e^{i\pi/4}\frac{1}{2}\partial_t^{-1/2}\left[
\phi_{p_1}(+1)\phi_{p_2}(+1)\left(+J_1^{-1}\partial_{y_1}u+J_2^{-1}\partial_{y_2}u\right)_{y_1=+1,y_2=+1}
\right]=-\frac{3}{4}\left[
\phi_{p_1}(+1)\phi_{p_2}(+1)\left(u\right)_{y_1=+1,y_2=+1}\right],\\
&e^{i\pi/4}\frac{1}{2}\partial_t^{-1/2}\left[
\phi_{p_1}(-1)\phi_{p_2}(+1)\left(-J_1^{-1}\partial_{y_1}u+J_2^{-1}\partial_{y_2}u\right)_{y_1=-1,y_2=+1}\right]
=-\frac{3}{4}\left[
\phi_{p_1}(-1)\phi_{p_2}(+1)\left(u\right)_{y_1=-1,y_2=+1}\right],\\
&e^{i\pi/4}\frac{1}{2}\partial_t^{-1/2}\left[
\phi_{p_1}(+1)\phi_{p_2}(-1)\left(+J_1^{-1}\partial_{y_1}u-J_2^{-1}\partial_{y_2}u\right)_{y_1=+1,y_2=-1}\right]
=-\frac{3}{4}\left[
\phi_{p_1}(+1)\phi_{p_2}(-1)\left(u\right)_{y_1=+1,y_2=-1}\right],\\
&e^{i\pi/4}\frac{1}{2}\partial_t^{-1/2}\left[
\phi_{p_1}(-1)\phi_{p_2}(-1)\left(-J_1^{-1}\partial_{y_1}u-J_2^{-1}\partial_{y_2}u\right)_{y_1=-1,y_2=-1}\right]
=-\frac{3}{4}\left[
\phi_{p_1}(-1)\phi_{p_2}(-1)\left(u\right)_{y_1=-1,y_2=-1}\right].
\end{split}
\end{equation}
Identifying the terms corresponding to the mass and stiffness matrices 
in~\eqref{eq:vform-2d} together with 
$\mathcal{I}_{1,a_1}$ and $\mathcal{I}_{2,a_2}$, the dynamical system 
works out to be
\begin{equation}\label{eq:linear-sys-hf}
\begin{split}
i M_1\partial_t\what{U}M_2
&=J_1^{-2}S_1\what{U}M_2+J_2^{-2}M_1\what{U}S_2
+e^{-i\pi/4}\partial_t^{1/2}\left[(J_1)^{-1}\Lambda_1\what{U}M_2+(J_2)^{-1}M_1\what{U}\Lambda_2\right]\\
&+e^{i\pi/4}\frac{1}{2}\partial_t^{-1/2}\left[
(J_1)^{-1}(J_2)^{-2}\Lambda_1\what{U}S_2+(J_1)^{-2}(J_2)^{-1}S_1\what{U}\Lambda_2\right]
+\frac{3}{4}(J_1J_2)^{-1}\Lambda_1\what{U}\Lambda_2,
\end{split}
\end{equation}
where the last term contains all the contributions from the corners.

%=============================================================================%
%=============================================================================%
%=============================================================================%
%=============================================================================%
\subsection{Nonlocal temporal discretization: Convolution quadrature}\label{sec:NL-TD-HF}
In this section, we present various time-stepping methods to discretize the 
master equation~\eqref{eq:linear-sys-hf}. The standard recipe here is to use a 
linear multistep method (LMM) together with a compatible 
convolution quadrature (CQ) scheme for the fractional operators. Such a 
dicretization is evidently nonlocal exhibiting the nonlocal nature of the 
dynamical system. To this end, let us identify the nonlocal terms 
in~\eqref{eq:linear-sys-hf} as follows:
\begin{equation}\label{eq:linear-sys-hf-split}
\begin{split}
&i M_1\partial_t\what{U}M_2
=\beta_1S_1\what{U}M_2+\beta_2M_1\what{U}S_2
+\frac{3}{4}\sqrt{\beta_1\beta_2}\Lambda_1\what{U}\Lambda_2
+\underbrace{\left[\partial_t^{1/2}\what{F}_++
\partial_t^{-1/2}\what{F}_-\right]}_{\text{non-local terms}}\;,\\
&\text{where}\quad
\left\{\begin{aligned}
&\what{F}_+ = e^{-i\pi/4}\left[\sqrt{\beta_1}\Lambda_1\what{U}M_2
+\sqrt{\beta_2}M_1\what{U}\Lambda_2\right],\\
&\what{F}_- = +e^{i\pi/4}\frac{1}{2}\left[
\sqrt{\beta_1}\beta_2\Lambda_1\what{U}S_2
+\beta_1\sqrt{\beta_2}S_1\what{U}\Lambda_2\right].
\end{aligned}\right.
\end{split}
\end{equation}
Let $\Delta t$ denote the time-step so that the temporal grid becomes 
$t_j=j\Delta t,\;j\in\field{N}_0$. The fractional operators 
$\partial_t^{\pm1/2}$ can be implemented using the 
LMM based CQ method. For a given LMM, let the CQ weights be denoted by 
$\{\omega^{(\pm1/2)}_{k}|\;k\in\field{N}_0\}$, then
\begin{equation}
[\partial_t^{\pm1/2}\what{F}_{\pm}]^{j+1}
=\rho^{\pm 1/2}F_{\pm}^{j+1}
+\rho^{\pm 1/2}\sum_{k=1}^{j}\omega^{(\pm1/2)}_{j+1-k}F_{\pm}^k,
\end{equation}
where $\rho$ is a function of $\Delta t$ and it is chosen such that 
$\omega^{(\pm1/2)}_{0}=1$. In the following we restrict ourselves to the 
trapezoidal rule (TR) and backward differentiation formulae (BDF) of order 
$1$ and $2$. The recipe for computing the CQ weights for each 
of the time-stepping methods is enumerated below:
%=============================================================================%
\begin{itemize}
\item CQ--BDF1 : The discretization scheme for the time-fractional operators is said to be 
`CQ-BDF1' if the underlying time-stepping method used in designing the quadrature 
is BDF1. Let $\nu\in\{ +1/2, -1/2\} $ and $\rho =1/\Delta t$, then the quadrature weights 
can be computed as follows 
\begin{equation}
\omega^{(\nu)}_k=\left[{\left(k-1-\nu \right)}/{k}\right]\omega^{(\nu)}_{k-1},\quad k\geq 1,
\quad\text{with}\quad\omega^{(\nu)}_0=1.
\end{equation}
\item CQ--BDF2 : Let $\rho =3/(2\Delta t)$, then the quadrature weights for BDF2 based CQ
can be computed as 
\begin{equation}
\omega^{(\nu)}_k=\left[{4\left(k-1-\nu \right)}/{3k}\right]\omega^{(\nu)}_{k-1}
-\left[{\left(k-2-2\nu \right)}/{3k}\right]\omega^{(\nu)}_{k-2} ,\quad k\geq 2,
\quad\text{with}\quad\omega^{(\nu)}_0=1,\;\omega^{(\nu)}_1 = -4\nu/3.
\end{equation}
\item CQ--TR : Let $\rho =2/\Delta t$, then the quadrature weights for the TR based
CQ scheme can be computed as
\begin{equation}
\omega^{(1/2)}_k=(-1)^j\omega^{(-1/2)}_k, \quad
\omega^{(-1/2)}_k=
\begin{cases}
C_{k/2},      & k\;\text{even},\\
C_{(k-1)/2},  & k\;\text{odd},\\
\end{cases}
\quad\text{where}\quad C_n=\frac{1\cdot 3\cdots(2n-1)}{n!2^n}.
\end{equation}
The quadrature weights can also be
generated using the following recurrence relation~\cite{SV2023}:
\begin{equation}
(k+1)\omega^{(\nu)}_{k+1}=(k-1)\omega^{(\nu)}_{k-1}-2\nu\omega^{(\nu)}_k,\quad k\geq1,
\quad\text{with}\quad \omega^{(\nu)}_0=1,\;\omega^{(\nu)}_1=-2\nu.
\end{equation}
 
\end{itemize}
%=============================================================================%
\subsubsection{CQ--BDF1}
The discrete numerical scheme for the dynamical system~\eqref{eq:linear-sys-hf} 
is labelled according to the CQ scheme used for the nonlocal terms. Therefore, 
the scheme `CQ-BDF1' will use BDF1 for the time-stepping of the local terms 
in~\eqref{eq:linear-sys-hf} which is discussed below.

Recalling $\rho =1/\Delta t$ and $\alpha_k=\sqrt{\rho/\beta_k}e^{-i\pi/4}$ 
for $k=1,2$ for the BDF1 scheme, the complete discretization 
of~\eqref{eq:linear-sys-hf} reads as
\begin{equation}
\begin{split}
& M_1\what{U}^{j+1}M_2+\alpha^{-2}_1S_1\what{U}^{j+1}M_2  +\alpha^{-2}_2M_1\what{U}^{j+1}S_2
  +\left[\alpha_1^{-1}\Lambda_1\what{U}^{j+1}M_2+\alpha_2^{-1}M_1\what{U}^{j+1}\Lambda_2\right]\\
& +\frac{1}{2}\left[\alpha_1^{-1}\alpha_2^{-2}\Lambda_1\what{U}^{j+1}S_2
  +\alpha_1^{-2}\alpha_2^{-1}S_1\what{U}^{j+1}\Lambda_2\right]
  +\frac{3}{4}(\alpha_1\alpha_2)^{-1}\Lambda_1\what{U}^{j+1}\Lambda_2
  =M_1\what{U}^{j}M_2 +\mathcal{B}^{j+1}_+ +\mathcal{B}^{j+1}_-,
\end{split}
\end{equation}
where the history terms are given given by
\begin{equation}
\left\{\begin{aligned}
\mathcal{B}^{j+1}_+&=\sum_{k=1}^{j}\omega^{(1/2)}_{j+1-k}\left[
\alpha_1^{-1}\Lambda_1\what{U}^{k}M_2 + \alpha_2^{-1}M_1\what{U}^{k}\Lambda_2\right],\\
\mathcal{B}^{j+1}_-&=\frac{1}{2}\sum_{k=1}^{j}\omega^{(-1/2)}_{j+1-k}
\left[\alpha_1^{-1}\alpha_2^{-2}\Lambda_1\what{U}^kS_2 
+\alpha_1^{-2}\alpha_2^{-1}S_1\what{U}^k\Lambda_2\right].
\end{aligned}\right.
\end{equation}
We may refer to these functions as the \emph{history functions} on account of
the fact that, for the latest time-step, they can be assumed to be determined 
by the quantities that are already determined in the previous time-steps.
%=============================================================================%    
\begin{rem}\label{rem:cq-bdf1-hist}
As noted in Sec.~\ref{sec:variational-prelim}, $\Lambda_1\what{U}$ and 
$\what{U}\Lambda_2$ denote the restrictions of the field $\what{U}$ 
to the boundary segments. Therefore, for the numerical implementation 
we can rewrite the history terms defined above for 
each of the boundary segments separately as follows:
\begin{equation*}
\begin{aligned}
&\what{\vv{\mathcal{B}}}^{j+1}_{a_1,+} 
= \sum_{k=1}^j\omega^{(1/2)}_{j+1-k}\;\hat{\vv{u}}^k_{a_1},
&&\what{\vv{\mathcal{B}}}^{j+1}_{a_2,+} 
= \sum_{k=1}^j\omega^{(1/2)}_{j+1-k}\;\hat{\vv{u}}^k_{a_2}, \quad a_1\in\{l,r\},
\\
&\what{\vv{\mathcal{B}}}^{j+1}_{a_1,-} 
= \sum_{k=1}^j\omega^{(-1/2)}_{j+1-k}\;\hat{\vv{u}}^k_{a_1},
&&\what{\vv{\mathcal{B}}}^{j+1}_{a_2,-} 
=\sum_{k=1}^j\omega^{(-1/2)}_{j+1-k}\;\hat{\vv{u}}^k_{a_2}, \quad a_2\in\{b,t\}.
\end{aligned}
\end{equation*}
It is straightforward to verify that the history terms, $\mathcal{B}^{j+1}_+$ 
and $\mathcal{B}^{j+1}_-$, can now be efficiently computed as
\begin{equation}
\begin{split}
&\mathcal{B}^{j+1}_+ = \frac{1}{\alpha_1}\vv{e}^{(1)}_0\otimes
\left[\what{\vs{\mathcal{B}}}^{j+1}_{l,+}\right]^{\tp}M_2
+\frac{1}{\alpha_1}\vv{e}^{(1)}_1\otimes
\left[\what{\vs{\mathcal{B}}}^{j+1}_{r,+}\right]^{\tp}M_2
+\frac{1}{\alpha_2}
\left[M_1\what{\vs{\mathcal{B}}}^{j+1}_{b,+}\right]\otimes\left(\vv{e}_0^{(2)}\right)^{\tp}
+\frac{1}{\alpha_2}
\left[M_1\what{\vs{\mathcal{B}}}^{j+1}_{t,+}\right]\otimes\left(\vv{e}^{(2)}_1\right)^{\tp},\\
&\mathcal{B}^{j+1}_-  = \frac{1}{2\alpha_1\alpha^2_2}\left(\vv{e}^{(1)}_0\otimes
\left[\what{\vs{\mathcal{B}}}^{j+1}_{l,+}\right]^{\tp}S_2
+\vv{e}^{(1)}_1\otimes
\left[\what{\vs{\mathcal{B}}}^{j+1}_{r,+}\right]^{\tp}S_2\right)
+\frac{1}{2\alpha^2_1\alpha_2}\left(
\left[S_1\what{\vs{\mathcal{B}}}^{j+1}_{b,+}\right]\otimes\left(\vv{e}_0^{(2)}\right)^{\tp}
+\left[S_1\what{\vs{\mathcal{B}}}^{j+1}_{t,+}\right]\otimes\left(\vv{e}^{(2)}_1\right)^{\tp}
\right).
\end{split}
\end{equation}
where $\vv{e}^{(p)}_0=(\delta_{k,0})_{k\in\field{J}_p}$ and 
$\vv{e}^{(p)}_1=(\delta_{k,1})_{k\in\field{J}_p}$ for $p=1,2$ are the orthogonal 
unit vectors in $\field{R}^{N_p+1}$.
\end{rem}
%=============================================================================%
%=============================================================================%
\subsubsection{CQ--BDF2}
The discrete numerical scheme for the dynamical system~\eqref{eq:linear-sys-hf}
is labelled as `CQ-BDF2' if the underlying time-stepping method is BDF2 and 
rest of the labels follow the convention defined above. 
Recalling $\rho =3/2\Delta t$, the complete discretization of the dynamical 
system in~\eqref{eq:linear-sys-hf} reads as
\begin{equation*}
\begin{split}
&M_1\what{U}^{j+2}M_2+\alpha^{-2}_1S_1\what{U}^{j+2}M_2+\alpha^{-2}_2M_1\what{U}^{j+2}S_2
+\frac{1}{2}\left[\alpha_1^{-1}\alpha_2^{-2}\Lambda_1\what{U}^{j+2}S_2
+\alpha_1^{-2}\alpha_2^{-1}S_1\what{U}^{j+2}\Lambda_2\right]\\
&+\left[\alpha_1^{-1}\Lambda_1\what{U}^{j+2}M_2+\alpha_2^{-1}M_1\what{U}^{j+2}\Lambda_2\right]
+\frac{3}{4}(\alpha_1\alpha_2)^{-1}\Lambda_1\what{U}^{j+2}\Lambda_2
=\left[\frac{4}{3} M_1\what{U}^{j+1}M_2-\frac{1}{3}M_1\what{U}^{j}M_2\right]
+\mathcal{B}^{j+2}_+ +\mathcal{B}^{j+2}_-,\\
\end{split}
\end{equation*}
where the history terms are computed in similar manner as evaluated in BDF1 case.
%=============================================================================%
%=============================================================================%
\subsubsection{CQ--TR}
The discrete numerical scheme for the dynamical system~\eqref{eq:linear-sys-hf}
is labelled as `CQ-TR' if the underlying time-stepping method is TR and rest of 
the labels follow the convention defined above. Recalling $\rho =2/\Delta t$ and 
setting $\what{V}^{j+1}=(\what{U}^{j+1}+\what{U}^{j})/2 $, the complete 
discretization of the dynamical system in~\eqref{eq:linear-sys-hf} reads as
\begin{equation*}
\begin{split}
&M_1\what{V}^{j+1}M_2+\alpha^{-2}_1S_1\what{V}^{j+1}M_2+\alpha^{-2}_2M_1\what{V}^{j+1}S_2
+\left[\alpha_1^{-1}\Lambda_1\what{V}^{j+1}M_2+\alpha_2^{-1}M_1\what{V}^{j+1}\Lambda_2\right]\\
&+\frac{1}{2}\left[\alpha_1^{-1}\alpha_2^{-2}\Lambda_1\what{V}^{j+1}S_2
+\alpha_1^{-2}\alpha_2^{-1}S_1\what{V}^{j+1}\Lambda_2\right]
+\frac{3}{4}(\alpha_1\alpha_2)^{-1}\Lambda_1\what{V}^{j+1}\Lambda_2
=M_1\what{U}^{j}M_2+\mathcal{B}^{j+1/2}_+ +\mathcal{B}^{j+1/2}_-, \\
\end{split}
\end{equation*}
where the history terms are computed as 
$\mathcal{B}^{j+1/2}_{\pm} = ( \mathcal{B}^{j+1}_{\pm}+\mathcal{B}^{j}_{\pm})/2$.
%=============================================================================%
%=============================================================================%
\subsection{Local temporal discretization: Conventional Pad\'e}\label{sec:L-TD-HF}
In this section, we present an effectively local form of the (nonlocal) dynamical 
system stated in~\eqref{eq:linear-sys-hf-split}. This is accomplished by 
employing a Pad\'e approximant based representation for the time-fractional 
operators in~\eqref{eq:linear-sys-hf-split} to obtain a system of ODEs which 
approximate the original system.  

Let the $K$-th order diagonal Pad\'{e} approximant based rational approximation 
for the function $\sqrt{z}$ is denoted by $R_K^{(1/2)}(z)$ defined 
in~\eqref{eq:sqrt-pade}. Here, we choose the rational approximation for the 
function $1/\sqrt{z}$, denoted by $R_K^{(-1/2)}(z)$, to ensure that
\begin{equation}
\left[\frac{R_K^{(1/2)}(z)}{R_K^{(-1/2)}(z)}\right] = z
\implies R_K^{(-1/2)}(z)=\frac{R_K^{(1/2)}(z)}{z},
\end{equation}
yielding
\begin{equation}\label{eq:sqrt-pade1}
R_K^{(-1/2)}(z)
=\frac{b_0}{z}-\sum_{k=1}^K\frac{b_k}{z(z+\eta^2_k)}
=\frac{b_0}{z}-\sum_{k=1}^K \frac{b_k}{\eta_k^2}\left[\frac{1}{z}-\frac{1}{z+\eta^2_k}\right]
=\frac{d_0}{z}-\sum_{k=1}^K\frac{d_k}{z+\eta^2_k},
\end{equation}
where $d_0 = b_0-\sum_{k=1}^K b_k/\eta_k^2$ and $d_k = -b_k/\eta_k^2$.

%=============================================================================%
Let $I_k$ be the identity matrix of dimension $(N_k+1)\times (N_k+1)$ for $k=1,2$. 
In order to write the rational approximations for the 
time-fractional operators present in~\eqref{eq:linear-sys-hf-split}, we introduce 
the auxiliary fields $\what{\mathcal{A}}_{k,1}$ and $\what{\mathcal{A}}_{k,2}$ 
where $ k=0,1,\ldots,K$ such that  
\begin{equation}\label{eq:pade-apprxm}
\begin{split}
&\partial_t^{1/2}\left[\Lambda_1\what{U}\right]
\approx b_0\Lambda_1\what{U}
-\sum_{k=1}^Kb_{k}\what{\mathcal{A}}_{k,1}
=b_0\Lambda_1\what{U}-(\vv{b}^{\tp}\otimes I_1)\what{\mathcal{A}}_1,\\
&\partial_t^{1/2}\left[\what{U}\Lambda_2\right]
\approx b_0\what{U}\Lambda_2
-\sum_{k=1}^Kb_{k}\what{\mathcal{A}}_{k,2}
=b_0\what{U}\Lambda_2-\what{\mathcal{A}}_2(\vv{b}\otimes I_2),\\
&\partial_t^{-1/2}\left[\Lambda_1\what{U}\right]
\approx d_0\what{\mathcal{A}}_{0,1}
-\sum_{k=1}^Kd_{k}\what{\mathcal{A}}_{k,1}
=d_0\what{A}_{0,1}-(\vv{d}^{\tp}\otimes I_1)\what{\mathcal{A}}_1,\\
&\partial_t^{-1/2}\left[\what{U}\Lambda_2\right]
\approx d_0\what{\mathcal{A}}_{0,2}
-\sum_{k=1}^Kd_{k}\what{\mathcal{A}}_{k,2}
=d_0\what{A}_{0,2}-\what{\mathcal{A}}_2(\vv{d}\otimes I_2),
\end{split}
\end{equation}
where $\vv{b}=(b_1,b_2,\ldots,b_K)^{\tp}$, 
$\vv{d}=(d_1,d_2,\ldots,d_K)^{\tp}$ and `$\otimes$' denotes the Kronecker product.
Note that the index $k$ of the auxiliary fields correspond to the 
$k$-th partial fraction of the Pad\'e approximant (with the exception of $b_0$). 

Let us now discuss the structure of the matrices 
$\what{\mathcal{A}}_{k,1}$ and $\what{\mathcal{A}}_{k,2}$ introduced above.
In Sec.~\ref{sec:variational-prelim}, we discussed how to obtain the 
representation of the restriction of the field $u_{N_1,N_2}(\vv{y},t)$ to 
the boundary segments in the Lobatto basis. Following the same convention, we observe
from~\eqref{eq:restrict-seg} that $\Lambda_1\what{U}$ and $\what{U}\Lambda_2$
are representations of the restriction of the field to the boundary segments.
Therefore, the auxiliary fields $\mathcal{A}_{k,1}$ and $\mathcal{A}_{k,2}$ are
functions on the boundary segments so that their expansion in the Lobatto basis reads 
as follows:
\begin{equation}
\begin{split}
\mathcal{A}_{k,1}(\pm1,y_2,t)
&=\sum_{p_1=0,1}\sum_{p_2=0}^{N_{2}}
\what{\mathcal{A}}_{p_1,p_2}^{(k,1)}\phi_{p_1}(\pm1)\phi_{p_2}(y_2),\\
\mathcal{A}_{k,2}(y_1,\pm1,t)
&=\sum_{p_1=0}^{N_{1}}\sum_{p_2=0,1}
\what{\mathcal{A}}_{p_1,p_2}^{(k,2)}\phi_{p_2}(y_1)\phi_{p_2}(\pm1),
\end{split}
\end{equation}
where we have exploited the property of Lobatto polynomials to establish that only 
first two rows of the matrices $\what{\mathcal{A}}_{k,1}$ are non-zero while 
only first two columns of the matrices $\what{\mathcal{A}}_{k,2}$ are non-zero. 
Next, we introduce a block matrix $\what{\mathcal{A}}_1$ with the matrix entries 
$\what{\mathcal{A}}_{k,1}$ defined as follows:
\begin{equation}
\what{\mathcal{A}}_{1}=
\begin{pmatrix}
\what{\mathcal{A}}_{1,1}\\
\what{\mathcal{A}}_{2,1}\\
\vdots\\
\what{\mathcal{A}}_{K,1}
\end{pmatrix}\quad\text{where}\quad
\what{\mathcal{A}}_{k,1}=
\begin{pmatrix}
\what{\mathcal{A}}_{0,0} & \what{\mathcal{A}}_{0,1} &\ldots & \what{\mathcal{A}}_{0,N_2}\\
\what{\mathcal{A}}_{1,0} & \what{\mathcal{A}}_{1,1} &\ldots & \what{\mathcal{A}}_{1,N_2}\\
0                        &                        0 &\ldots & 0\\
\vdots                   & \vdots                   &       & \vdots\\
0                        &                        0 &\ldots & 0
\end{pmatrix}^{(k,1)}\in\field{C}^{(N_1+1)\times(N_2+1)}.
\end{equation}
Similarly, we introduce a block matrix $\what{\mathcal{A}}_2$ with the matrix entries 
$\what{\mathcal{A}}_{k,2}\in\field{C}^{(N_1+1)\times(N_2+1)}$ defined as follows:
\begin{equation}
\what{\mathcal{A}}_{2}=
\begin{pmatrix}
\what{\mathcal{A}}_{1,2}&
\ldots&
\what{\mathcal{A}}_{K,2}
\end{pmatrix}\quad\text{where}\quad
\what{\mathcal{A}}_{k,2}=
\begin{pmatrix}
\what{\mathcal{A}}_{0,0}   & \what{\mathcal{A}}_{0,1}   &     0&\ldots& 0\\
\what{\mathcal{A}}_{1,0}   & \what{\mathcal{A}}_{1,1}   &     0&\ldots& 0\\
\vdots                     & \vdots                     &\vdots& &\vdots\\
\what{\mathcal{A}}_{N_1,0} & \what{\mathcal{A}}_{N_1,1} &     0&\ldots& 0
\end{pmatrix}^{(k,2)}.
\end{equation}

%=============================================================================%
Let us recall that the index $k$ of the auxiliary fields correspond to the 
$k$-th partial fraction of the Pad\'e approximant so that one can identify the 
ODE satisfied by each of the auxiliary fields:
\begin{equation}
\begin{aligned}
&\partial_t\what{\mathcal{A}}_{0,1}=\Lambda_1\what{U},
&&\partial_t\what{\mathcal{A}}_{0,2}=\what{U}\Lambda_2,\\
&\partial_t\what{\mathcal{A}}_{k,1}=-\eta_k^2\what{\mathcal{A}}_{k,1}+\Lambda_1\what{U},
&&\partial_t\what{\mathcal{A}}_{k,2}=-\eta_k^2\what{\mathcal{A}}_{k,2}+\what{U}\Lambda_2.
\end{aligned}
\end{equation}
Let us introduce the matrices
$\mathcal{E}_K=\diag(\eta_1^2,\eta_2^2,\ldots,\eta_K^2)\in\field{R}^{(K\times K)}$
and $\vv{1}_K = (1,1,\ldots,1)^{\tp}\in\field{R}^{K}$ which can be used to 
rewrite the ODEs as follows:
\begin{equation}\label{eq:odes-auxi-pade}
\begin{aligned}
&\partial_t\what{\mathcal{A}}_{0,1}=\Lambda_1\what{U},
&&\partial_t\what{\mathcal{A}}_{0,2}=\what{U}\Lambda_2,\\
&\partial_t\what{\mathcal{A}}_{1}
=-(\mathcal{E}_K\otimes I_1)\what{\mathcal{A}}_{1}
+\vv{1}_K\otimes(\Lambda_1\what{U}),
&&\partial_t\what{\mathcal{A}}_{2}
=-\what{\mathcal{A}}_{2}(\mathcal{E}_K\otimes I_2)
+\vv{1}^{\tp}_K\otimes(\what{U}\Lambda_2).
\end{aligned}
\end{equation}
Employing the approximations developed in~\eqref{eq:pade-apprxm} 
in the (non-local) dynamical system stated in~\eqref{eq:linear-sys-hf-split}, we obtain
an effectively local system given by 
\begin{equation}\label{eq:linear-sys-hf-pade}
\begin{split}
i M_1\partial_t\what{U}M_2
&=\beta_1S_1\what{U}M_2+\beta_2M_1\what{U}S_2
+e^{-i\pi/4}b_0\left[\sqrt{\beta_1}\Lambda_1\what{U}M_2
+\sqrt{\beta_2}M_1\what{U}\Lambda_2\right]\\
&\quad-e^{-i\pi/4}\left[\sqrt{\beta_1}(\vv{b}^{\tp}\otimes I_1)\what{\mathcal{A}}_{1}M_2
+\sqrt{\beta_2}M_1\what{\mathcal{A}}_{2}(\vv{b}\otimes I_2)\right]
+d_0e^{i\pi/4}\frac{1}{2}
\left[\sqrt{\beta_1}\beta_2\what{\mathcal{A}}_{0,1}S_2
+\beta_1\sqrt{\beta_2}S_1\what{\mathcal{A}}_{0,2}\right]\\
&\quad-e^{i\pi/4}\frac{1}{2}
\left[\sqrt{\beta_1}\beta_2(\vv{d}^{\tp}\otimes I_1)\what{\mathcal{A}}_{1}S_2
+\beta_1\sqrt{\beta_2}S_1\what{\mathcal{A}}_{2}(\vv{d}\otimes I_2)\right]
+\frac{3}{4}\sqrt{\beta_1\beta_2}\,\Lambda_1\what{U}\Lambda_2.
\end{split}
\end{equation}
This dynamical system together with~\eqref{eq:odes-auxi-pade} can be 
solved further with any choice of temporal discretization methods 
(namely, one-step and multistep methods).
%=============================================================================%
\begin{rem}\label{rem:aux-fns-hf}
For the numerical implementation, only the non-zero entries of 
$\what{\mathcal{A}}_{k,1}$ and $\what{\mathcal{A}}_{k,2}$ are relevant.
Therefore, we may introduce auxiliary functions 
\begin{equation*}
\begin{aligned}
 &\Gamma_l:\;\varphi_{k,l}(y_2,t)=\mathcal{A}_{k,1}(-1,y_2,t),
&&\Gamma_r:\;\varphi_{k,r}(y_2,t)=\mathcal{A}_{k,1}(+1,y_2,t),\\
 &\Gamma_b:\;\varphi_{k,b}(y_1,t)=\mathcal{A}_{k,1}(y_1,-1,t),
&&\Gamma_t:\;\varphi_{k,t}(y_1,t)=\mathcal{A}_{k,1}(y_1,+1,t). 
\end{aligned}
\end{equation*}
Recalling $a_1\in\{l,r\}$ and $a_2\in\{b,t\}$, the expansion in terms of Lobatto 
basis reads as
\begin{equation*}
\varphi_{k,a_1}(y_2,t)
=\sum_{p_2=0}^{N_{2}}\what{\varphi}_{k,p_2,a_1}(t)\,\phi_{p_2}(y_2),\quad
\varphi_{k,a_2}(y_1,t)
=\sum_{p_1=0}^{N_{1}}\what{\varphi}_{k,p_1,a_2}(t)\,\phi_{p_1}(y_1),
\end{equation*}
so that 
\begin{equation*}
\begin{aligned}
&\what{\mathcal{A}}_{k,1}=
\begin{pmatrix}
\what{\varphi}_{k,0,l} & \what{\varphi}_{k,1,l} &\ldots & \what{\varphi}_{k,N_2,l}\\
\what{\varphi}_{k,0,r} & \what{\varphi}_{k,1,r} &\ldots & \what{\varphi}_{k,N_2,r}\\
0                        &                        0 &\ldots & 0\\
\vdots                   & \vdots                   &\dots  & \vdots\\
0                        &                        0 &\ldots & 0
\end{pmatrix}
=
\begin{pmatrix}
\what{\vs{\varphi}}_{k,l}^{\tp}\\
\what{\vs{\varphi}}_{k,r}^{\tp}\\
Z_1
\end{pmatrix},
&&
\what{\vs{\varphi}}_{k,a_1}=
\begin{pmatrix}
\what{\varphi}_{k,0,a_1}\\
\what{\varphi}_{k,1,a_1}\\
\vdots\\
\what{\varphi}_{k,N_2,a_1}
\end{pmatrix},\\
&\what{\mathcal{A}}_{k,2}=
\begin{pmatrix}
\what{\varphi}_{k,0,b} &\what{\varphi}_{k,0,t} & 0 & \ldots & 0\\ 
\what{\varphi}_{k,1,b} &\what{\varphi}_{k,1,t} & 0 & \ldots & 0\\
\vdots                 &\vdots                 & \vdots & &\vdots\\
\what{\varphi}_{k,N_1,b} &\what{\varphi}_{k,N_1,t} & 0 & \ldots & 0 
\end{pmatrix}=
\begin{pmatrix}
\what{\vs{\varphi}}_{k,b}& \what{\vs{\varphi}}_{k,t}& Z_2
\end{pmatrix},
&&
\what{\vs{\varphi}}_{k,a_2}=
\begin{pmatrix}
\what{\varphi}_{k,0,a_2}\\
\what{\varphi}_{k,1,a_2}\\
\vdots\\
\what{\varphi}_{k,N_1,a_2}
\end{pmatrix},
\end{aligned}
\end{equation*}
where $Z_1=0_{(N_1-1)\times(N_2+1)}$ and $Z_2=0_{(N_1+1)\times(N_2-1)}$.
\end{rem}
The approach presented here is referred to as the \emph{conventional Pad\'e} 
approach. Each of the methods discussed below is derived by applying a 
particular time-stepping scheme to the conventional Pad\'e approach, 
therefore, we label them as `CP--' followed by the acronym for the time-stepping method.

%=============================================================================%    
\subsubsection{CP--BDF1}
The discrete numerical scheme for the dynamical system defined
in~\eqref{eq:odes-auxi-pade} and~\eqref{eq:linear-sys-hf-pade}
is labelled as `CP--BDF1' if the underlying time-stepping method is BDF1. 

Recalling $\rho=1/(\Delta t)$, the BDF1-based discretization of the dynamical 
system in~\eqref{eq:linear-sys-hf-pade} reads as
\begin{equation}\label{eq:bdf1-main}
\begin{split}
&i \rho M_1\what{U}^{j+1}M_2 -i \rho M_1\what{U}^{j}M_2\\
&=\beta_1S_1\what{U}^{j+1}M_2+\beta_2M_1\what{U}^{j+1}S_2
+e^{-i\pi/4}b_0\left[\sqrt{\beta_1}\Lambda_1\what{U}^{j+1}M_2
+\sqrt{\beta_2}M_1\what{U}^{j+1}\Lambda_2\right]
+\frac{3}{4}(J_1J_2)^{-1}\Lambda_1\what{U}^{j+1}\Lambda_2\\
&\quad-e^{-i\pi/4}\left[\sqrt{\beta_1}(\vv{b}^{\tp}\otimes I_1)\what{\mathcal{A}}^{j+1}_{1}M_2
+\sqrt{\beta_2}M_1\what{\mathcal{A}}^{j+1}_{2}(\vv{b}\otimes I_2)\right]
+d_0e^{i\pi/4}\frac{1}{2}
\left[\sqrt{\beta_1}\beta_2\what{\mathcal{A}}^{j+1}_{0,1}S_2
+\beta_1\sqrt{\beta_2}S_1\what{\mathcal{A}}^{j+1}_{0,2}\right]\\
&\quad-e^{i\pi/4}\frac{1}{2}
\left[\sqrt{\beta_1}\beta_2(\vv{d}^{\tp}\otimes I_1)\what{\mathcal{A}}^{j+1}_{1}S_2
+\beta_1\sqrt{\beta_2}S_1\what{\mathcal{A}}^{j+1}_{2}(\vv{d}\otimes I_2)\right],
\end{split}
\end{equation}
and that of the auxiliary fields in~\eqref{eq:odes-auxi-pade} reads as
\begin{equation}\label{eq:bdf1-aux}
\begin{aligned}
\what{\mathcal{A}}^{j+1}_{0,1}&=\frac{1}{\rho}\Lambda_1\what{U}^{j+1}+\what{\mathcal{A}}^{j}_{0,1},
&&\what{\mathcal{A}}^{j+1}_{0,2}=\frac{1}{\rho}\what{U}^{j+1}\Lambda_2+\what{\mathcal{A}}^{j}_{0,2},\\
\what{\mathcal{A}}^{j+1}_{1}&= (\mathcal{G}_K\otimes I_1)
\left[\frac{1}{\rho}\vv{1}_K\otimes(\Lambda_1\what{U}^{j+1})+\what{\mathcal{A}}^{j}_{1}\right],
&&\what{\mathcal{A}}^{j+1}_{2}=\left[\frac{1}{\rho}\vv{1}^{\tp}_K\otimes(\what{U}^{j+1}\Lambda_2)
+\what{\mathcal{A}}^{j}_{2}\right]
(\mathcal{G}_K\otimes I_2),
\end{aligned}
\end{equation}
where we have introduced 
\begin{equation}
\mathcal{G}_K = \diag\left(\frac{1}{1+\ovl{\eta}_1^2},\frac{1}{1+\ovl{\eta}_2^2},\ldots,
\frac{1}{1+\ovl{\eta}_K^2}\right)\in\field{R}^{(K\times K)}, \quad \ovl{\eta}_k = \eta_k/\sqrt{\rho}.
\end{equation}
%=============================================================================%
Plugging-in the updates for the auxiliary fields from~\eqref{eq:bdf1-aux} 
in the discrete system~\eqref{eq:bdf1-main} and simplifying, we obtain
\begin{equation*}
\begin{split}
&M_1\what{U}^{j+1}M_2 - M_1\what{U}^{j}M_2\\
&=-\alpha_1^{-2}S_1\what{U}^{j+1}M_2-\alpha_2^{-2}M_1\what{U}^{j+1}S_2
-\ovl{b}_0\left[\alpha_1^{-1}\Lambda_1\what{U}^{j+1}M_2
+\alpha_2^{-1} M_1\what{U}^{j+1}\Lambda_2\right]
-\frac{3}{4\alpha_1\alpha_2}\Lambda_1\what{U}^{j+1}\Lambda_2\\
&\quad + \alpha_1^{-1}
\left[\frac{1}{\rho}\left(\ovl{\vv{b}}^{\tp}\mathcal{G}_K\vv{1}_K\right)\otimes(\Lambda_1\what{U}^{j+1})
+\left(\ovl{\vv{b}}^{\tp}\mathcal{G}_K\otimes I_1\right)\what{\mathcal{A}}^{j}_{1}\right] M_2
 +\alpha_2^{-1} M_1
\left[\frac{1}{\rho}\left(\vv{b}^{\tp}\mathcal{G}_K\vv{1}_K\right)\otimes(\what{U}^{j+1}\Lambda_2)
+\what{\mathcal{A}}^{j}_{2}\left(\mathcal{G}_K\vv{b}\otimes I_2\right)\right]\\
&\quad +\frac{1}{2} \alpha_1^{-1}\alpha_2^{-2}
\left[\left(\frac{1}{\rho}\left(\vv{d}^{\tp}\mathcal{G}_K\vv{1}_K\right)\otimes(\Lambda_1\what{U}^{j+1})
+\left(\vv{d}^{\tp}\mathcal{G}_K\otimes I_1\right)\what{\mathcal{A}}^{j}_{1}\right)
-\ovl{d}_0\left(\frac{1}{\rho}\Lambda_1\what{U}^{j+1}+\what{\mathcal{A}}^{j}_{0,1}\right)\right]S_2\\
&\quad +\frac{1}{2} \alpha_2^{-1}\alpha_1^{-2}S_1
\left[\left(\frac{1}{\rho}\left(\vv{d}^{\tp}\mathcal{G}_K\vv{1}_K\right)\otimes(\what{U}^{j+1}\Lambda_2)
+\what{\mathcal{A}}^{j}_{2}\left(\mathcal{G}_K\vv{d}\otimes I_2\right) \right)      
-\ovl{d}_0\left(\frac{1}{\rho}\what{U}^{j+1}\Lambda_2+\what{\mathcal{A}}^{j}_{0,2}\right)\right].
\end{split}
\end{equation*}
Introducing the quantities
\begin{equation}\label{eq:varpi-pm1o2}
\left\{\begin{aligned}
&\Gamma^{(1/2)}_k=\frac{\ovl{b}_k}{1+\ovl{\eta}_k^2},\quad k=1,\ldots,K,
&&\Gamma^{(-1/2)}_k=\frac{\ovl{d}_k}{1+\ovl{\eta}_k^2},
\quad k=1,\ldots,K,\\
&\varpi^{(1/2)}=\ovl{b}_0-\frac{1}{\rho}\ovl{\vv{b}}^{\tp}\left(\mathcal{G}_K\vv{1}_K\right)
=\ovl{b}_0-\frac{1}{\rho}\sum_{k=1}^K\Gamma^{(1/2)}_k,
&&\varpi^{(-1/2)}=\frac{\ovl{d}_0}{\rho}-\frac{1}{\rho}\ovl{\vv{d}}^{\tp}
\left(\mathcal{G}_K\vv{1}_K\right)
=\frac{\ovl{d}_0}{\rho}-\frac{1}{\rho}\sum_{k=1}^K\Gamma^{(-1/2)}_k.
\end{aligned}\right.
\end{equation}
and identifying the terms with the latest time-step of the field, $\what{U}^{j+1}$, we obtain
\begin{equation}
\begin{split}
& M_1\what{U}^{j+1}M_2+\alpha^{-2}_1S_1\what{U}^{j+1}M_2
+\alpha^{-2}_2M_1\what{U}^{j+1}S_2
+\varpi^{(1/2)}\left[\alpha^{-1}_1\Lambda_1\what{U}^{j+1}M_2
+\alpha_2^{-1}M_1\what{U}^{j+1}\Lambda_2\right]\\
&+\frac{1}{2}\varpi^{(-1/2)}
\left[\alpha^{-1}_1\alpha_2^{-2}\Lambda_1\what{U}^{j+1}S_2
+\alpha_1^{-2}\alpha_2^{-1}S_1\what{U}^{j+1}\Lambda_2\right]
+\frac{3}{4\alpha_1\alpha_2}\Lambda_1\what{U}^{j+1}\Lambda_2 
=M_1\what{U}^{j}M_2 +\mathcal{B}^{j+1}_+ +\mathcal{B}^{j+1}_-,
\end{split}
\end{equation}
where the history terms are given given by
\begin{equation}
\left\{\begin{aligned}
\mathcal{B}^{j+1}_+&=
\alpha^{-1}_1\left[\left(\ovl{\vv{b}}^{\tp}\mathcal{G}_K\right)\otimes I_1\right]
\what{\mathcal{A}}_1^{j}M_2+\alpha_2^{-1}M_1\what{\mathcal{A}}_2^{j}
\left[\left(\mathcal{G}_K\ovl{\vv{b}}\right)\otimes I_2\right],\\
\mathcal{B}^{j+1}_-&=
\frac{1}{2}\alpha^{-1}_1\alpha_2^{-2}\left[
\left(\ovl{\vv{d}}^{\tp}\mathcal{G}_K\otimes I_1 \right)\what{\mathcal{A}}_1^{j}
-\ovl{d}_0\what{\mathcal{A}}_{0,1}^{j}\right]S_2
+\frac{1}{2}\alpha^{-2}_1\alpha_2^{-1}S_1 \left[ 
\what{\mathcal{A}}_2^{j}\left(\mathcal{G}_K\ovl{\vv{d}}\otimes I_2\right)
-\what{\mathcal{A}}_{0,2}^{j}\ovl{d}_0 \right].
\end{aligned}\right.
\end{equation}
The latest value, $\what{U}^{j+1}$, can then be used to carry out the updates 
for the auxiliary fields as in~\eqref{eq:bdf1-aux} which will be used to 
compute the history terms in the next time-step.
%=============================================================================%
\begin{rem}\label{rem:bdf1-hist}
As noted in Remark~\ref{rem:aux-fns-hf}, only the non-zero entries of 
$\what{\mathcal{A}}_{1}$ and $\what{\mathcal{A}}_{2}$ are relevant for the 
numerical implementation so that we may rewrite the history terms 
defined above for each of the boundary segments separately as follows:
\begin{equation*}
\begin{aligned}
&\what{\vv{\mathcal{B}}}^{j+1}_{a_1,+} 
= \sum_{k=1}^K\Gamma^{(1/2)}_k\what{\vv{\varphi}}_{k,a_1},
&&\what{\vv{\mathcal{B}}}^{j+1}_{a_2,+} 
= \sum_{k=1}^K\Gamma^{(1/2)}_k\what{\vv{\varphi}}_{k,a_2},
\\
&\what{\vv{\mathcal{B}}}^{j+1}_{a_1,-} 
= \sum_{k=1}^K\Gamma^{(-1/2)}_k\what{\vv{\varphi}}_{k,a_1}
-\ovl{d}_0\what{\vv{\varphi}}_{0,a_1},
&&\what{\vv{\mathcal{B}}}^{j+1}_{a_2,-} 
=\sum_{k=1}^K\Gamma^{(-1/2)}_k\what{\vv{\varphi}}_{k,a_2}
-\ovl{d}_0\what{\vv{\varphi}}_{0,a_2}.
\end{aligned}
\end{equation*}
It is straightforward to verify that the history terms, $\mathcal{B}^{j+1}_+$ 
and $\mathcal{B}^{j+1}_-$, can now be efficiently computed as
\begin{equation}
\begin{split}
&\mathcal{B}^{j+1}_+ = \frac{1}{\alpha_1}\vv{e}^{(1)}_0\otimes
\left[\what{\vs{\mathcal{B}}}^{j+1}_{l,+}\right]^{\tp}M_2
+\frac{1}{\alpha_1}\vv{e}^{(1)}_1\otimes
\left[\what{\vs{\mathcal{B}}}^{j+1}_{r,+}\right]^{\tp}M_2
+\frac{1}{\alpha_2}
\left[M_1\what{\vs{\mathcal{B}}}^{j+1}_{b,+}\right]\otimes\left(\vv{e}_0^{(2)}\right)^{\tp}
+\frac{1}{\alpha_2}
\left[M_1\what{\vs{\mathcal{B}}}^{j+1}_{t,+}\right]\otimes\left(\vv{e}^{(2)}_1\right)^{\tp},\\
&\mathcal{B}^{j+1}_-  = \frac{1}{2\alpha_1\alpha^2_2}\left(\vv{e}^{(1)}_0\otimes
\left[\what{\vs{\mathcal{B}}}^{j+1}_{l,+}\right]^{\tp}S_2
+\vv{e}^{(1)}_1\otimes
\left[\what{\vs{\mathcal{B}}}^{j+1}_{r,+}\right]^{\tp}S_2\right)
+\frac{1}{2\alpha^2_1\alpha_2}\left(
\left[S_1\what{\vs{\mathcal{B}}}^{j+1}_{b,+}\right]\otimes\left(\vv{e}_0^{(2)}\right)^{\tp}
+\left[S_1\what{\vs{\mathcal{B}}}^{j+1}_{t,+}\right]\otimes\left(\vv{e}^{(2)}_1\right)^{\tp}
\right).
\end{split}
\end{equation}
\end{rem}
%=============================================================================%
%=============================================================================%
\subsubsection{CP--BDF2}
The discrete numerical scheme for the dynamical system defined
in~\eqref{eq:odes-auxi-pade} and~\eqref{eq:linear-sys-hf-pade} is labelled as 
`CP--BDF2' if the underlying time-stepping method is BDF2.

Recalling $\rho=3/(2\Delta t)$, the BDF2-based discretization of the dynamical 
system in~\eqref{eq:linear-sys-hf-pade} reads as
\begin{equation}\label{eq:bdf2-main}
\begin{split}
&i \rho M_1\what{U}^{j+2}M_2 
-i\rho M_1\left(\frac{4}{3}\what{U}^{j+1}-\frac{1}{3}\what{U}^{j}\right)M_2\\
&=\beta_1S_1\what{U}^{j+2}M_2+\beta_2M_1\what{U}^{j+2}S_2
+e^{-i\pi/4}b_0\left[\sqrt{\beta_1}\Lambda_1\what{U}^{j+2}M_2
+\sqrt{\beta_2}M_1\what{U}^{j+2}\Lambda_2\right]
+\frac{3}{4}(J_1J_2)^{-1}\Lambda_1\what{U}^{j+2}\Lambda_2\\
&\quad-e^{-i\pi/4}\left[\sqrt{\beta_1}(\vv{b}^{\tp}\otimes I_1)\what{\mathcal{A}}^{j+2}_{1}M_2
+\sqrt{\beta_2}M_1\what{\mathcal{A}}^{j+2}_{2}(\vv{b}\otimes I_2)\right]
+d_0e^{i\pi/4}\frac{1}{2}
\left[\sqrt{\beta_1}\beta_2\what{\mathcal{A}}^{j+2}_{0,1}S_2
+\beta_1\sqrt{\beta_2}S_1\what{\mathcal{A}}^{j+2}_{0,2}\right]\\
&\quad-e^{i\pi/4}\frac{1}{2}
\left[\sqrt{\beta_1}\beta_2(\vv{d}^{\tp}\otimes I_1)\what{\mathcal{A}}^{j+2}_{1}S_2
+\beta_1\sqrt{\beta_2}S_1\what{\mathcal{A}}^{j+2}_{2}(\vv{d}\otimes I_2)\right],
\end{split}
\end{equation}
and that of the auxiliary fields in~\eqref{eq:odes-auxi-pade} reads as
\begin{equation}\label{eq:bdf2-aux}
\begin{aligned}
&\what{\mathcal{A}}^{j+2}_{0,1}=\frac{1}{\rho}\Lambda_1\what{U}^{j+2}
+\frac{4}{3}\what{\mathcal{A}}^{j+1}_{0,1}-\frac{1}{3}\what{\mathcal{A}}^{j}_{0,1},
&&\what{\mathcal{A}}^{j+2}_{0,2}=\frac{1}{\rho}\what{U}^{j+2}\Lambda_2
+\frac{4}{3}\what{\mathcal{A}}^{j+1}_{0,2}-\frac{1}{3}\what{\mathcal{A}}^{j}_{0,2},\\
&\what{\mathcal{A}}^{j+2}_{1}=(\mathcal{G}_K\otimes I_1)
\left[\frac{1}{\rho}\vv{1}_K\otimes(\Lambda_1\what{U}^{j+2})+
\frac{4}{3}\what{\mathcal{A}}^{j+1}_{1}-\frac{1}{3}\what{\mathcal{A}}^{j}_{1}\right],
&&\what{\mathcal{A}}^{j+2}_{2} 
= \left[\frac{1}{\rho}\vv{1}^{\tp}_K\otimes(\what{U}^{j+2}\Lambda_2)
+\frac{4}{3}\what{\mathcal{A}}^{j+1}_{2}-\frac{1}{3}\what{\mathcal{A}}^{j}_{2}
\right](\mathcal{G}_K\otimes I_2).
\end{aligned}
\end{equation}
Following along similar lines as in BDF1, the discrete linear system
corresponding to~\eqref{eq:linear-sys-hf-pade} turns out to be
\begin{equation}
\begin{split}
& M_1\what{U}^{j+2}M_2+\alpha^{-2}_1S_1\what{U}^{j+2}M_2
+\alpha^{-2}_2M_1\what{U}^{j+2}S_2
+\varpi^{(1/2)}\left[\alpha^{-1}_1\Lambda_1\what{U}^{j+2}M_2
+\alpha_2^{-1}M_1\what{U}^{j+2}\Lambda_2\right]
+\frac{3}{4\alpha_1\alpha_2}\Lambda_1\what{U}^{j+2}\Lambda_2\\
&+\frac{1}{2}\varpi^{(-1/2)}
\left[\alpha^{-1}_1\alpha_2^{-2}\Lambda_1\what{U}^{j+2}S_2
+\alpha_1^{-2}\alpha_2^{-1}S_1\what{U}^{j+2}\Lambda_2\right]
=M_1\left(\frac{4}{3}\what{U}^{j+1}-\frac{1}{3}\what{U}^{j}\right)M_2
 +\mathcal{B}^{j+2}_+ +\mathcal{B}^{j+2}_-,
\end{split}
\end{equation}
where the history terms are given given by
\begin{equation}
\left\{\begin{aligned}
\mathcal{B}^{j+2}_+&=
\alpha^{-1}_1\left[\left(\ovl{\vv{b}}^{\tp}\mathcal{G}_K\right)\otimes I_1\right]
\left(\frac{4}{3}\what{\mathcal{A}}_1^{j+1}-\frac{1}{3}\what{\mathcal{A}}_1^{j}\right)M_2
+\alpha_2^{-1}M_1\left(\frac{4}{3}\what{\mathcal{A}}_2^{j+1}-\frac{1}{3}\what{\mathcal{A}}_2^{j}\right)
\left[\left(\mathcal{G}_K\ovl{\vv{b}}\right)\otimes I_2\right],\\
\mathcal{B}^{j+2}_-&=
\frac{1}{2}\alpha^{-1}_1\alpha_2^{-2}\left[
\left(\ovl{\vv{d}}^{\tp}\mathcal{G}_K\otimes I_1 \right)
\left(\frac{4}{3}\what{\mathcal{A}}_1^{j+1}-\frac{1}{3}\what{\mathcal{A}}_1^{j}\right)
-\ovl{d}_0
\left(\frac{4}{3}\what{\mathcal{A}}_{0,1}^{j+1}-\frac{1}{3}\what{\mathcal{A}}_{0,1}^{j}\right)
\right]S_2\\
& \quad +\frac{1}{2}\alpha^{-2}_1\alpha_2^{-1}S_1 \left[ 
\left(\frac{4}{3}\what{\mathcal{A}}_2^{j+1}-\frac{1}{3}\what{\mathcal{A}}_2^{j}\right)
\left(\mathcal{G}_K\ovl{\vv{d}}\otimes I_2\right)
-\left(\frac{4}{3}\what{\mathcal{A}}_{0,2}^{j+1}-\frac{1}{3}\what{\mathcal{A}}_{0,2}^{j}\right)
\ovl{d}_0\right].
\end{aligned}\right.
\end{equation}
The simplifications noted in Remark~\ref{rem:bdf1-hist} can also be applied 
here for the history terms. We omit the details because the adaption is 
straightforward. 
%=============================================================================%    
%=============================================================================%    
\subsubsection{CP--TR}
The discrete numerical scheme for the dynamical system is labelled as `CP--TR' 
if the underlying time-stepping method is TR.

Recalling $\rho=2/(\Delta t)$, the TR based discretization of the dynamical system
in~\eqref{eq:linear-sys-hf-pade} reads as
\begin{equation}\label{eq:tr-main}
\begin{split}
&i \rho M_1\what{U}^{j+1/2}M_2 -i \rho M_1\what{U}^{j}M_2\\
&=\beta_1S_1\what{U}^{j+1/2}M_2+\beta_2M_1\what{U}^{j+1/2}S_2
+e^{-i\pi/4}b_0\left[\sqrt{\beta_1}\Lambda_1\what{U}^{j+1/2}M_2
+\sqrt{\beta_2}M_1\what{U}^{j+1/2}\Lambda_2\right]
+\frac{3}{4}(J_1J_2)^{-1}\Lambda_1\what{U}^{j+1/2}\Lambda_2\\
&\quad-e^{-i\pi/4}\left[\sqrt{\beta_1}(\vv{b}^{\tp}\otimes I_1)\what{\mathcal{A}}^{j+1/2}_{1}M_2
+\sqrt{\beta_2}M_1\what{\mathcal{A}}^{j+1/2}_{2}(\vv{b}\otimes I_2)\right]
+d_0e^{i\pi/4}\frac{1}{2}
\left[\sqrt{\beta_1}\beta_2\what{\mathcal{A}}^{j+1/2}_{0,1}S_2
+\beta_1\sqrt{\beta_2}S_1\what{\mathcal{A}}^{j+1/2}_{0,2}\right]\\
&\quad-e^{i\pi/4}\frac{1}{2}
\left[\sqrt{\beta_1}\beta_2(\vv{d}^{\tp}\otimes I_1)\what{\mathcal{A}}^{j+1/2}_{1}S_2
+\beta_1\sqrt{\beta_2}S_1\what{\mathcal{A}}^{j+1/2}_{2}(\vv{d}\otimes I_2)\right],
\end{split}
\end{equation}
and that of the auxiliary fields in~\eqref{eq:odes-auxi-pade} reads as
\begin{equation}\label{eq:tr-aux}
\begin{aligned}
&\what{\mathcal{A}}^{j+1/2}_{0,1}=\frac{1}{\rho}\Lambda_1\what{U}^{j+1/2}+\what{\mathcal{A}}^{j}_{0,1},
&&\what{\mathcal{A}}^{j+1/2}_{0,2}=\frac{1}{\rho}\what{U}^{j+1/2}\Lambda_2+\what{\mathcal{A}}^{j}_{0,2},\\
&\what{\mathcal{A}}^{j+1/2}_{1} =(\mathcal{G}_K\otimes I_1)
\left[\frac{1}{\rho}\vv{1}_K\otimes(\Lambda_1\what{U}^{j+1/2})+\what{\mathcal{A}}^{j}_{1}\right],
&&\what{\mathcal{A}}^{j+1/2}_{2}=
\left[\frac{1}{\rho}\vv{1}^{\tp}_K\otimes(\what{U}^{j+1/2}\Lambda_2)+\what{\mathcal{A}}^{j}_{2}\right]
(\mathcal{G}_K\otimes I_2).
\end{aligned}
\end{equation}
Next, we plug-in these updates for auxiliary fields from~\eqref{eq:tr-aux} in the 
discrete system~\eqref{eq:tr-main} to obtain
\begin{equation}
\begin{split}
& M_1\what{U}^{j+1/2}M_2+\alpha^{-2}_1S_1\what{U}^{j+1/2}M_2
+\alpha^{-2}_2M_1\what{U}^{j+1/2}S_2
+\varpi^{(1/2)}\left[\alpha^{-1}_1\Lambda_1\what{U}^{j+1/2}M_2
+\alpha_2^{-1}M_1\what{U}^{j+1/2}\Lambda_2\right]\\
&+\frac{1}{2}\varpi^{(-1/2)}
\left[\alpha^{-1}_1\alpha_2^{-2}\Lambda_1\what{U}^{j+1/2}S_2
+\alpha_1^{-2}\alpha_2^{-1}S_1\what{U}^{j+1/2}\Lambda_2\right]
+\frac{3}{4\alpha_1\alpha_2}\Lambda_1\what{U}^{j+1/2}\Lambda_2
=M_1\what{U}^{j}M_2 +\mathcal{B}^{j+1/2}_+ +\mathcal{B}^{j+1/2}_-,
\end{split}
\end{equation}
where the history terms are given given by
\begin{equation}
\left\{\begin{aligned}
\mathcal{B}^{j+1/2}_+&=
\alpha^{-1}_1\left[\left(\ovl{\vv{b}}^{\tp}\mathcal{G}_K\right)\otimes I_1\right]
\what{\mathcal{A}}_1^{j}M_2+\alpha_2^{-1}M_1\what{\mathcal{A}}_2^{j}
\left[\left(\mathcal{G}_K\ovl{\vv{b}}\right)\otimes I_2\right],\\
\mathcal{B}^{j+1/2}_-&=
\frac{1}{2}\alpha^{-1}_1\alpha_2^{-2}\left[
\left(\ovl{\vv{d}}^{\tp}\mathcal{G}_K\otimes I_1 \right)\what{\mathcal{A}}_1^{j}
-\ovl{d}_0\what{\mathcal{A}}_{0,1}^{j}\right]S_2
+\frac{1}{2}\alpha^{-2}_1\alpha_2^{-1}S_1 \left[ 
\what{\mathcal{A}}_2^{j}\left(\mathcal{G}_K\ovl{\vv{d}}\otimes I_2\right)
-\what{\mathcal{A}}_{0,2}^{j}\ovl{d}_0 \right].
\end{aligned}\right.
\end{equation}
Once again, simplification of the history terms is possible along the lines of 
Remark~\ref{rem:bdf1-hist} but we omit the discussion for the sake of brevity of 
presentation. 

The computation of the history terms in the next time-step requires values of 
the auxiliary fields on the temporal grid which can be facilitated by the TR-based 
discretization of the ODEs in~\eqref{eq:odes-auxi-pade} as follows
\begin{equation}
\begin{aligned}
&\what{\mathcal{A}}^{j+1}_{0,1}=\frac{2}{\rho}\Lambda_1\what{U}^{j+1/2}+\what{\mathcal{A}}^{j}_{0,1},
&&\what{\mathcal{A}}^{j+1}_{0,2} =\frac{2}{\rho}\what{U}^{j+1/2}\Lambda_2+\what{\mathcal{A}}^{j}_{0,2},\\
&\what{\mathcal{A}}^{j+1}_{1} 
=\left[\frac{2}{\rho}\left(\mathcal{G}_K\vv{1}_K\right)\otimes(\Lambda_1\what{U}^{j+1/2})
+\left(\mathcal{H}_K\otimes I_1\right)\what{\mathcal{A}}^{j}_{1}\right],
&&\what{\mathcal{A}}^{j+1}_{2}=
\left[\frac{2}{\rho}\left(\mathcal{G}_K\vv{1}_K\right)^{\tp}\otimes(\what{U}^{j+1/2}\Lambda_2)
+\what{\mathcal{A}}^{j}_{2}\left(\mathcal{H}_K\otimes I_2\right)\right],
\end{aligned}
\end{equation}
where
\begin{equation}
\mathcal{H}_K = \diag\left(\frac{1-\ovl{\eta}_1^2}{1+\ovl{\eta}_1^2}
,\frac{1-\ovl{\eta}_2^2}{1+\ovl{\eta}_2^2},\ldots,
 \frac{1-\ovl{\eta}_K^2}{1+\ovl{\eta}_K^2}\right)\in\field{R}^{(K\times K)}.
\end{equation}

%=============================================================================%    
%=============================================================================%    
%=============================================================================%    
\section{Variational Formulation: TBCs}\label{sec:vform-np}
In this section, we consider the numerical solution of the IBVP 
in~\eqref{eq:2D-SE-CT} using a Legendre-Galerkin method as laid out 
in Sec.~\ref{sec:variational-prelim}. The variational formulation begins 
with~\eqref{eq:vform-2d} where the boundary terms, namely, 
$\mathcal{I}_1$ and $\mathcal{I}_2$, need to incorporate the transparent 
boundary conditions stated in~\eqref{eq:maps-cq}. In the following, we suppress 
the independent variables in the auxiliary function 
$\varphi(y_1,y_2,\tau_1,\tau_2)$ for the sake of brevity of presentation. Starting 
the with the segments denoted by $\Gamma_{a_1},\;a_1\in\{l,r\}$, we have
\begin{equation}
\begin{split}
\mathcal{I}_1&=J_1^{-2}\int_{\field{I}}\left(
[\theta_{p_1,p_2}\partial_{y_1}u]_{+1}-[\theta_{p_1,p_2}\partial_{y_1}u]_{-1}\right)dy_2\\
&=J_1^{-1}\int_{\field{I}}\left[\phi_{p_1}(+1)\left(\left.-e^{-i\pi/4}\partial^{1/2}_{\tau_1} 
\varphi\right|_{\tau_1,\tau_2=t} \right)_{y_1=+1}
+\phi_{p_1}(-1)\left(\left.-e^{-i\pi/4}\partial^{1/2}_{\tau_1}
\varphi\right|_{\tau_1,\tau_2=t}\right)_{y_1=-1}\right]\phi_{p_2}(y_2)dy_2\\
&=-e^{-i\pi/4}J_1^{-1}\partial^{1/2}_{\tau_1}\left.\left[\phi_{p_1}(+1)\int_{\field{I}}\phi_{p_2}(y_2)
\left(\varphi\right)_{y_1=+1}dy_2
+\phi_{p_1}(-1)\int_{\field{I}}\phi_{p_2}(y_2)\left(\varphi\right)_{y_1=-1}dy_2\right]
\right|_{\tau_1,\tau_2=t}.
\end{split}
\end{equation}
Similarly, for the segments denoted by $\Gamma_{a_2},\;a_2\in\{b,t\}$, we have
\begin{equation}
\begin{split}
\mathcal{I}_2&=J_2^{-2}\int_{\field{I}}\left(
[\theta_{p_1,p_2}\partial_{y_2}u]_{+1}-[\theta_{p_1,p_2}\partial_{y_2}u]_{-1}\right)dy_1\\
&=J_2^{-1}\int_{\field{I}}\left[\phi_{p_2}(+1)\left(\left.-e^{-i\pi/4}\partial^{1/2}_{\tau_2} 
\varphi\right|_{\tau_1,\tau_2=t} \right)_{y_2=+1}
+\phi_{p_2}(-1)\left(\left.-e^{-i\pi/4}\partial^{1/2}_{\tau_2}
\varphi\right|_{\tau_1,\tau_2=t}\right)_{y_2=-1}\right]\phi_{p_1}(y_1)dy_1\\
&=-e^{-i\pi/4}J_2^{-1}\partial^{1/2}_{\tau_2}\left.\left[\phi_{p_2}(+1)\int_{\field{I}}\phi_{p_1}(y_1)
\left(\varphi\right)_{y_2=+1}dy_1
+\phi_{p_2}(-1)\int_{\field{I}}\phi_{p_1}(y_1)\left(\varphi\right)_{y_2=-1}dy_1\right]
\right|_{\tau_1,\tau_2=t}.
\end{split}
\end{equation}
%=============================================================================%
Let $\what{\Phi}_{p_1,p_2}$ denote the expansion coefficient in the Lobatto 
basis for the auxiliary function $\varphi(y_1,y_2,\tau_1,\tau_2)$ so that
\begin{equation}
\varphi(y_1,y_2,\tau_1,\tau_2) =\sum_{p_1=0}^{N_1}\sum_{p_2=0}^{N_{2}}
\what{\Phi}_{p_1,p_2}(\tau_1,\tau_2)\phi_{p_1}(y_1)\phi_{p_2}(y_2).
\end{equation}
We introduce the matrix 
$\what{\Phi}=(\what{\Phi}_{p_1,p_2})\in\field{C}^{(N_1+1)\times(N_2+1)}$. 
Restrictions of the auxiliary function on the boundary segments 
$\Gamma_{a_1}$ and $\Gamma_{a_2}$ require representation of 
$\varphi(y_1=\pm1,y_2,\tau_1,\tau_2)$ and 
$\varphi(y_1,y_2=\pm1,\tau_1,\tau_2)$, respectively. Let $\what{\Phi}_k$ 
denote the restriction of $\what{\Phi}$ to $\Gamma_{a_k}$ for $k=1,2$, then
\begin{equation}
\what{\Phi}_1 = \Lambda_1\what{\Phi},\quad \what{\Phi}_2 = \what{\Phi}\Lambda_2,
\end{equation}
where $\Lambda_1$ and $\Lambda_2$ are defined in Sec.~\ref{sec:variational-prelim}.
The linear system corresponding to~\eqref{eq:vform-2d} turns out to be
\begin{equation}\label{eq:linear-sys-tbcs}
i M_1\partial_t\what{U}M_2=J_1^{-2}S_1\what{U}M_2+J_2^{-2}M_1\what{U}S_2
+e^{-i\pi/4}\left.\left[(J_1)^{-1}\partial_{\tau_1}^{1/2}\what{\Phi}_1(\tau_1,\tau_2)M_2
+(J_2)^{-1}M_1\partial_{\tau_2}^{1/2}\what{\Phi}_2(\tau_1,\tau_2)\right]\right|_{\tau_1,\tau_2=t}.
\end{equation}
Here, we observe that the implementation of the $1/2$-order temporal derivatives require the history 
of the auxiliary function from the start of the computations which is facilitated by the IVPs 
in~\eqref{eq:ivps-cq}. The variational formulation for these IVPs can be written as:
%=============================================================================%
\begin{equation}\label{eq:vform-aux}
\begin{split}
\Gamma_{a_2}:\quad
\left(i\partial_{\tau_1}\varphi+J_1^{-2}\partial^2_{y_1}\varphi, \phi_{p_1}
\right)_{\field{I}}
&=i\partial_{\tau_1}\left(\varphi,\phi_{p_1}\right)_{\field{I}}
-J_1^{-2}\left(\partial_{y_1}\varphi,\phi'_{p_1}\right)_{\field{I}}
+J_1^{-2}\left(
[\phi_{p_1}\partial_{y_1}\varphi]_{y_1=+1}
-[\phi_{p_1}\partial_{y_1}\varphi]_{y_1=-1}
\right),\\
\Gamma_{a_1}:\quad
\left(i\partial_{\tau_2}\varphi+J_2^{-2}\partial^2_{y_2}\varphi, \phi_{p_2}
\right)_{\field{I}}
&=i\partial_{\tau_2}\left(\varphi,\phi_{p_2}\right)_{\field{I}}
-J_2^{-2}\left(\partial_{y_2}\varphi,\phi'_{p_2}\right)_{\field{I}}
+J_2^{-2}\left(
[\phi_{p_2}\partial_{y_2}\varphi]_{y_2=+1}
-[\phi_{p_2}\partial_{y_2}\varphi]_{y_2=-1}\right).
\end{split}
\end{equation}
Supplying the maps present in~\eqref{eq:maps-cq-auxi} to the boundary terms for the 
segments $\Gamma_{a_2}$ above, we have
\begin{equation}
\begin{split}
J_1^{-2}\left(
 [\phi_{p_1}\partial_{y_1}\varphi]_{y_1=+1}
-[\phi_{p_1}\partial_{y_1}\varphi]_{y_1=-1}\right)
&=J_1^{-1}\left[\phi_{p_1}(+1)\left(
-e^{-i\pi/4}\partial^{1/2}_{\tau_1}\varphi\right)_{y_1=+1}
+\phi_{p_1}(-1)\left(-e^{-i\pi/4}\partial^{1/2}_{\tau_1}\varphi\right)_{y_1=-1}
\right]\\
&=-e^{-i\pi/4}J_1^{-1}\partial^{1/2}_{\tau_1}\left[\phi_{p_1}(+1)
\left(\varphi\right)_{y_1=+1}
+\phi_{p_1}(-1)\left(\varphi\right)_{y_1=-1}\right],
\end{split}
\end{equation}
Similarly, for the segments $\Gamma_{a_1}$, we have
\begin{equation}
\begin{split}
J_2^{-2}\left(
 [\phi_{p_2}\partial_{y_2}\varphi]_{+1}
-[\phi_{p_2}\partial_{y_2}\varphi]_{-1}\right)
&=J_2^{-1}\left[
 \phi_{p_2}(+1)\left(-e^{-i\pi/4}\partial^{1/2}_{\tau_2}\varphi\right)_{y_2=+1}
+\phi_{p_2}(-1)\left(-e^{-i\pi/4}\partial^{1/2}_{\tau_2}\varphi\right)_{y_2=-1}\right]\\
&=-e^{-i\pi/4}J_2^{-1}\partial^{1/2}_{\tau_2}\left[
 \phi_{p_2}(+1)\left(\varphi\right)_{y_2=+1}
+\phi_{p_2}(-1)\left(\varphi\right)_{y_2=-1}\right].
\end{split}
\end{equation}
Introducing the appropriate mass and stiffness matrices, the variational form 
in~\eqref{eq:vform-aux} can be stated as
\begin{equation}\label{eq:linear-sys-ivps}
\begin{split}
& i M_1\partial_{\tau_1}\what{\Phi}_2(\tau_1,\tau_2)
=J_1^{-2}S_1\what{\Phi}_2(\tau_1,\tau_2)
+e^{-i\pi/4}\left[\partial_{\tau_1}^{1/2}(J_1)^{-1}\Lambda_1\what{\Phi}_2(\tau_1,\tau_2)
\right],\quad\tau_1\in(\tau_2,t],\\
&i \partial_{\tau_2}\what{\Phi}_1(\tau_1,\tau_2)M_2
=J_2^{-2}\what{\Phi}_1(\tau_1,\tau_2)S_2
+e^{-i\pi/4}\left[\partial_{\tau_2}^{1/2}(J_2)^{-1}\what{\Phi}_1(\tau_1,\tau_2)\Lambda_2
\right],\quad\tau_2\in(\tau_1,t],
\end{split}
\end{equation}
where the initial conditions are 
$\what{\Phi}_2(\tau_2,\tau_2)=\what{U}(\tau_2)\Lambda_2$ and 
$\what{\Phi}_1(\tau_1,\tau_1)=\Lambda_1\what{U}(\tau_1)$.
%=============================================================================%
%=============================================================================%
\subsection{Nonlocal temporal discretization: Convolution quadrature}
In this section, we present various time-stepping method to discretize the 
dynamical system defined by~\eqref{eq:linear-sys-tbcs} 
and~\eqref{eq:linear-sys-ivps}. Here, we employ one-step methods, namely, 
BDF1 and TR, together with a compatible convolution quadrature scheme for the 
fractional operators. Such a discretization scheme is evidently nonlocal exhibiting 
the nonlocal nature of the dynamical system. 

The recipe for handling the nonlocal terms in context of CQ is presented 
in Sec~\ref{sec:NL-TD-HF}. The computation of quadrature sums require history of 
the auxiliary fields on the segments which can be facilitated by solving the 
dynamical systems established in~\eqref{eq:linear-sys-ivps}. Each of the methods 
discussed below is derived by applying a particular time-stepping scheme to the 
CQ approach, therefore, we label them as `CQ--' followed by the acronym for the 
time-stepping method.
%=============================================================================%
\subsubsection{CQ--BDF1}
The discrete numerical scheme for the dynamical system~\eqref{eq:linear-sys-tbcs} 
is labelled as `CQ--BDF1' where the underlying one-step method is BDF1.
Recalling $\rho =1/\Delta t$ and $\alpha_k=\sqrt{\rho/\beta_k}e^{-i\pi/4}$ for $k=1,2$, 
the complete discretization of~\eqref{eq:linear-sys-tbcs} reads as
\begin{equation}
M_1\what{U}^{j+1}M_2+\alpha^{-2}_1S_1\what{U}^{j+1}M_2+\alpha^{-2}_2M_1\what{U}^{j+1}S_2
+\alpha_1^{-1}\Lambda_1\what{U}^{j+1}M_2 \\
+\alpha_2^{-1}M_1\what{U}^{j+1}\Lambda_2 =M_1\what{U}^{j}M_2 - \mathcal{B}^{j+1},
\end{equation}
where we have used the fact that $\what{\Phi}_1^{j+1,j+1}=\Lambda_1\what{U}^{j+1}$ 
and $\what{\Phi}_2^{j+1,j+1}=\what{U}^{j+1}\Lambda_2$. The history term is given by
\begin{equation}
\mathcal{B}^{j+1} = \sum_{k=1}^{j}\omega^{(1/2)}_{j+1-k}\left[\alpha_1^{-1}\what{\Phi}_1^{k,j+1}M_2
+\alpha_2^{-1}M_1\what{\Phi}_2^{j+1,k}\right].
\end{equation}
The computation of the history term above requires the off-diagonal samples 
$\what{\Phi}_1^{k,j+1}$ and $\what{\Phi}_2^{j+1,k}$ which can be obtained by 
solving the dynamical system defined in~\eqref{eq:linear-sys-ivps} corresponding 
to the boundary segments. The BDF1-based temporal discretization of this 
system reads as
\begin{equation}
\begin{split}
\left(M_1 + \alpha_1^{-2}S_1+\alpha_1^{-1}\Lambda_1\right)\what{\Phi}_2^{j+1,p}
=M_1\what{\Phi}_2^{j,p}  
-\mathcal{C}_2^{j+1},\quad p=0,1,\ldots,j,\\
 \what{\Phi}_1^{q,j+1}\left(M_2 + \alpha_2^{-2}S_2
 +\alpha_2^{-1}\Lambda_2\right)= \what{\Phi}_1^{q,j}M_2  
-\mathcal{C}_1^{j+1},\quad q=0,1,\ldots,j.
\end{split}
\end{equation}
where the history terms are given by
\begin{equation}
\mathcal{C}_1^{j+1} = \sum_{k=1}^{j}\omega^{(1/2)}_{j+1-k}\left[\alpha_2^{-1}
\what{\Phi}_1^{q,k}\Lambda_2\right] ,\quad
\mathcal{C}_2^{j+1} = \sum_{k=1}^{j}\omega^{(1/2)}_{j+1-k}\left[\alpha_1^{-1}\Lambda_1
\what{\Phi}_2^{k,p}\right]. 
\end{equation}
%=============================================================================%    
\subsubsection{CQ--TR}
The discrete numerical scheme for the dynamical system~\eqref{eq:linear-sys-tbcs} 
is labelled as `CQ--TR' where the underlying one-step method is TR.
Recalling $\rho =2/\Delta t$, the complete discretization of~\eqref{eq:linear-sys-tbcs} 
reads as
\begin{equation}
M_1\what{U}^{j+1/2}M_2+\alpha^{-2}_1S_1\what{U}^{j+1/2}M_2+\alpha^{-2}_2M_1\what{U}^{j+1/2}S_2
+\alpha_1^{-1}\Lambda_1\what{U}^{j+1/2}M_2+\alpha_2^{-1}M_1\what{U}^{j+1/2}\Lambda_2\\
=M_1\what{U}^{j}M_2 -\mathcal{B}^{j+1/2},
\end{equation}
where the history term is computed in the similar manner as done in BDF1. 
Here, $\mathcal{B}^{j+1/2}=(\mathcal{B}^{j+1}+\mathcal{B}^{j})/2$. 

In order to obtain the off-diagonal samples $\what{\Phi}_1^{k,j+1}$ and 
$\what{\Phi}_2^{j+1,k}$, we solve the dynamical system defined  
in~\eqref{eq:linear-sys-ivps} corresponding to the boundary segments. The TR-based
temporal discretization of this system reads as
\begin{equation}
\begin{split}
\left(M_1+\alpha_1^{-2}S_1+\alpha_1^{-1}\Lambda_1\right)\what{\Phi}_2^{j+1/2,p}
=M_1\what{\Phi}_2^{j,p}-\mathcal{C}_2^{j+1/2},\quad p=0,1,\ldots,j,\\
\what{\Phi}_1^{q,j+1/2}
\left(M_2 + \alpha_2^{-2}S_2+\alpha_2^{-1}\Lambda_2\right)= \what{\Phi}_1^{q,j}M_2  
-\mathcal{C}_1^{j+1/2},\quad q=0,1,\ldots,j.
\end{split}
\end{equation}
where the history terms corresponding to the corners are computed in the 
similar manner as done in BDF1. Here, we use the fact that 
$\mathcal{C}_k^{j+1/2}=(\mathcal{C}_k^{j+1}+\mathcal{C}_k^{j})/2$ for $k=1,2$.  
%=============================================================================%
%=============================================================================%
\subsection{Local temporal discretization: Novel Pad\'e}
In this section, we present an effectively local form of the (nonlocal) dynamical 
system stated in~\eqref{eq:linear-sys-tbcs} and~\eqref{eq:linear-sys-ivps}. 
This is accomplished by 
employing a novel Pad\'e approximant based representation for the time-fractional 
operators in~\eqref{eq:linear-sys-tbcs} and~\eqref{eq:linear-sys-ivps}
to obtain a system of ODEs which approximate the original system. 
This approach can be contrasted with the CQ-approach where the history grows with 
each time-step making it expensive computationally as well as from a storage point of view. 

In order to write the rational approximations for the time-fractional operators
present in~\eqref{eq:linear-sys-tbcs}, we introduce the auxiliary fields
$\mathcal{A}_{k,1}$ and $\mathcal{A}_{k,2}$ where $ k=1,2,\ldots,K$ such that  
\begin{equation}\label{eq:npade-apprxm}
\begin{split}
&\left.\partial_{\tau_1}^{1/2}\what{\Phi}_1\right|_{\tau_1,\tau_2=t}
=b_0\Lambda_1\what{U}
-\sum_{k=1}^Kb_{k}\what{\mathcal{A}}_{k,1}(t,t)
=b_0\Lambda_1\what{U}-(\vv{b}^{\tp}\otimes I_1)\what{\mathcal{A}}_1(t,t),\\
&\left.\partial_{\tau_2}^{1/2}\what{\Phi}_2\right|_{\tau_1,\tau_2=t}
=b_0\what{U}\Lambda_2
-\sum_{k=1}^Kb_{k}\what{\mathcal{A}}_{k,2}(t,t)
=b_0\what{U}\Lambda_2-\what{\mathcal{A}}_2(t,t)(\vv{b}\otimes I_2).
\end{split}
\end{equation}
The auxiliary fields $\mathcal{A}_{k,1}$ and $\mathcal{A}_{k,2}$ can be 
expanded in the Lobatto basis as follows:
\begin{equation}
\begin{split}
\mathcal{A}_{k,1}(\pm1,y_2,\tau_1,\tau_2)
&=\sum_{p_1=0,1}\sum_{p_2=0}^{N_{2}}
\what{\mathcal{A}}_{p_1,p_2}^{(k,1)}(\tau_1,\tau_2)\,\phi_{p_1}(\pm1)\phi_{p_2}(y_2),\\
\mathcal{A}_{k,2}(y_1,\pm1,\tau_1,\tau_2)
&=\sum_{p_1=0}^{N_{1}}\sum_{p_2=0,1}
\what{\mathcal{A}}_{p_1,p_2}^{(k,2)}(\tau_1,\tau_2)\,\phi_{p_2}(y_1)\phi_{p_2}(\pm1).
\end{split}
\end{equation}
The structure of the matrices $\what{\mathcal{A}}_{k,1}$ and 
$\what{\mathcal{A}}_{k,2}$ are same as that defined in Sec.~\ref{sec:L-TD-HF}. 
The structure of the block matrix $\what{\mathcal{A}}_p$ with the matrix entries 
$\what{\mathcal{A}}_{k,p}$ for $p=1,2$ is also borrowed from the same section.

As noted before, the index $k$ of the auxiliary fields correspond to the 
$k$-th partial fraction of the Pad\'e approximant so that the ODE satisfied 
by each of the auxiliary fields reads as
\begin{equation}\label{eq:odes-auxi-npade}
\begin{aligned}
&\partial_{\tau_1}\what{\mathcal{A}}_{k,1}(\tau_1,\tau_2)
=-\eta_k^2\what{\mathcal{A}}_{k,1}(\tau_1,\tau_2)+\what{\Phi}_1,
&&\partial_{\tau_2}\what{\mathcal{A}}_{k,2}(\tau_1,\tau_2)
=-\eta_k^2\what{\mathcal{A}}_{k,2}(\tau_1,\tau_2) +\what{\Phi}_2,\\
&\partial_{\tau_1}\what{\mathcal{A}}_{1}(\tau_1,\tau_2)
=-(\mathcal{E}_K\otimes I_1)\what{\mathcal{A}}_{1}(\tau_1,\tau_2)
+\vv{1}_K\otimes\what{\Phi}_1,
&&\partial_{\tau_2}\what{\mathcal{A}}_{2}(\tau_1,\tau_2)
=-\what{\mathcal{A}}_{2}(\tau_1,\tau_2)(\mathcal{E}_K\otimes I_2)
+\vv{1}^{\tp}_K\otimes\what{\Phi}_2.
\end{aligned}
\end{equation}
%=============================================================================%
Let us observe that the non-zero entries of $\what{\mathcal{A}}_{k,1}$ and 
$\what{\mathcal{A}}_{k,2}$ are related to the auxiliary functions introduced in 
Sec.~\ref{sec:CT-NP} as follows:
\begin{equation*}
\begin{aligned}
 &\Gamma_l:\;\varphi_{k,l}(y_2,\tau_1,\tau_2)=\mathcal{A}_{k,1}(-1,y_2,\tau_1,\tau_2),
&&\Gamma_r:\;\varphi_{k,r}(y_2,\tau_1,\tau_2)=\mathcal{A}_{k,1}(+1,y_2,\tau_1,\tau_2),\\
 &\Gamma_b:\;\varphi_{k,b}(y_1,\tau_1,\tau_2)=\mathcal{A}_{k,1}(y_1,-1,\tau_1,\tau_2),
&&\Gamma_t:\;\varphi_{k,t}(y_1,\tau_1,\tau_2)=\mathcal{A}_{k,1}(y_1,+1,\tau_1,\tau_2). 
\end{aligned}
\end{equation*}
Following the convention laid out in Remark~\eqref{rem:aux-fns-hf}, 
expansion of the auxiliary fields in the Lobatto basis read as
\begin{equation}\label{eq:aux-np}
\begin{split}
&\varphi_{k,a_1}(y_2,\tau_1,\tau_2)
=\sum_{p_2=0}^{N_{2}}\what{\varphi}_{k,p_2,a_1}(\tau_1,\tau_2)\phi_{p_2}(y_2),\quad a_1\in\{l,r\},\\
&\varphi_{k,a_2}(y_1,\tau_1,\tau_2)
=\sum_{p_1=0}^{N_{1}}\what{\varphi}_{k,p_1,a_2}(\tau_1,\tau_2)\phi_{p_1}(y_1),\quad a_2\in\{b,t\},
\end{split}
\end{equation}
and the matrices $\what{\mathcal{A}}_{k,1}$ and 
$\what{\mathcal{A}}_{k,2}$ can be written as
\begin{equation}
\what{\mathcal{A}}_{k,1}=
\begin{pmatrix}
\what{\vs{\varphi}}_{k,l}^{\tp}\\
\what{\vs{\varphi}}_{k,r}^{\tp}\\
Z_1
\end{pmatrix},\quad 
\what{\mathcal{A}}_{k,2}=
\begin{pmatrix}
\what{\vs{\varphi}}_{k,b}& \what{\vs{\varphi}}_{k,t}& Z_2
\end{pmatrix}.
\end{equation}
Employing the approximations developed in~\eqref{eq:npade-apprxm} 
in the (non-local) dynamical system stated in~\eqref{eq:linear-sys-tbcs}, we obtain
an effectively local system given by 
\begin{equation}\label{eq:linear-sys-npade}
\begin{split}
i M_1\partial_t\what{U}M_2
&=J_1^{-2}S_1\what{U}M_2+J_2^{-2}M_1\what{U}S_2
+e^{-i\pi/4}b_0\left[(J_1)^{-1}\Lambda_1\what{U}M_2
+(J_2)^{-1}M_1\what{U}\Lambda_2\right]\\
&\quad -e^{-i\pi/4}\left[(J_1)^{-1}(\vv{b}^{\tp}\otimes I_1)\what{\mathcal{A}}_1(t,t)M_2
+(J_2)^{-1}M_1\what{\mathcal{A}}_2(t,t)(\vv{b}\otimes I_2)\right].
\end{split}
\end{equation}
%=============================================================================%
As noted in Sec.~\ref{sec:CT-NP}, the auxiliary fields $\mathcal{A}_{k,1}$ and 
$\mathcal{A}_{k,2}$ satisfy the same IVPs as that of ${\Phi}_{1}$ and $\Phi_{2}$ 
on the boundary segments. Therefore, the dynamical system for 
$\mathcal{A}_{k,1}$ and $\mathcal{A}_{k,2}$ work out to be (similar to that 
in~\eqref{eq:linear-sys-ivps}):
\begin{equation}\label{eq:ivp-pade-auxi-mat}
\begin{split}
\Gamma_{a_1}:\quad
& i \partial_{\tau_2}\what{\mathcal{A}}_{k,1}(\tau_1,\tau_2)M_2
=J_2^{-2}\what{\mathcal{A}}_{k,1}(\tau_1,\tau_2)S_2
+e^{-i\pi/4}\left[\partial_{\tau_2}^{1/2}(J_2)^{-1}
\what{\mathcal{A}}_{k,1}(\tau_1,\tau_2)\Lambda_2\right],\\
\Gamma_{a_2}:\quad
&i \partial_{\tau_1}M_1\what{\mathcal{A}}_{k,2}(\tau_1,\tau_2)
=J_1^{-2}S_1\what{\mathcal{A}}_{k,2}(\tau_1,\tau_2)
+e^{-i\pi/4}\left[\partial_{\tau_1}^{1/2}(J_1)^{-1}
\Lambda_1\what{\mathcal{A}}_{k,2}(\tau_1,\tau_2)\right],
\end{split}
\end{equation}
which can be restated in terms of block matrices corresponding to these
auxiliary fields as
\begin{equation}\label{eq:ivp-pade-auxi-mat2}
\begin{split}
\Gamma_{a_1}:\quad
& i \partial_{\tau_2}\what{\mathcal{A}}_{1}(\tau_1,\tau_2)M_2
=J_2^{-2}\what{\mathcal{A}}_{1}(\tau_1,\tau_2)S_2
+e^{-i\pi/4}\left[\partial_{\tau_2}^{1/2}(J_2)^{-1}
\what{\mathcal{A}}_{1}(\tau_1,\tau_2)\Lambda_2\right],\\
\Gamma_{a_2}:\quad
&i \partial_{\tau_1}M_1\what{\mathcal{A}}_{2}(\tau_1,\tau_2)
=J_1^{-2}S_1\what{\mathcal{A}}_{2}(\tau_1,\tau_2)
+e^{-i\pi/4}\left[\partial_{\tau_1}^{1/2}(J_1)^{-1}
\Lambda_1\what{\mathcal{A}}_{2}(\tau_1,\tau_2)\right].
\end{split}
\end{equation} 
In order to write the rational approximations for the (non-local) time-fractional operators
present in~\eqref{eq:ivp-pade-auxi-mat}, we introduce the auxiliary fields
$\mathcal{C}_{k,k'}$ where $ k=k'=1,2,\ldots,K$ such that  
\begin{equation}\label{eq:npade-apprxm-corners}
\begin{split}
&\partial_{\tau_1}^{1/2}(\Lambda_1\what{\mathcal{A}}_{k,2})
=b_0\Lambda_1\what{\mathcal{A}}_{k,2}
-\sum_{k'=1}^Kb_{k'}\what{\mathcal{C}}_{k',k}(\tau_1,\tau_2)
=b_0\Lambda_1\what{\mathcal{A}}_{2} 
-(\vv{b}^{\tp}\otimes I_1)\what{\mathcal{C}}(\tau_1,\tau_2),\\
&\partial_{\tau_2}^{1/2}(\what{\mathcal{A}}_{k,1}\Lambda_2)
=b_0\what{\mathcal{A}}_{k,1}\Lambda_2
-\sum_{k'=1}^Kb_{k'}\what{\mathcal{C}}_{k,k'}(\tau_1,\tau_2)
=b_0\what{\mathcal{A}}_{1}\Lambda_2 
-\what{\mathcal{C}}(\tau_1,\tau_2)(\vv{b}\otimes I_2).
\end{split}
\end{equation}
The auxiliary field $\mathcal{C}_{k,k'}$ can be expanded in the Lobatto basis as follows: 
\begin{equation}
\mathcal{C}_{k,k'}(\pm1,\pm1,\tau_1,\tau_2)
=\sum_{p_1=0,1}\sum_{p_2=0,1}
\what{\mathcal{C}}_{p_1,p_2}^{(k,k')}\phi_{p_1}(\pm1)\phi_{p_2}(\pm1),\quad
\end{equation}
Exploiting the property of Lobatto polynomials, we can establish that only 
four entries in the matrix $\what{\mathcal{C}}_{k,k'}$ are non-zero. Next
we introduce a block matrix $\what{\mathcal{C}}$ with the matrix entries 
$\what{\mathcal{C}}_{k,k'}\in\field{C}^{(N_1+1)\times (N_2+1)}$ defined as follows:
\begin{equation}
\what{\mathcal{C}}=
\begin{pmatrix}
\what{\mathcal{C}}_{1,1} & \what{\mathcal{C}}_{1,2} &\ldots & \what{\mathcal{C}}_{1,K} \\
\what{\mathcal{C}}_{2,1} & \what{\mathcal{C}}_{2,2} &\ldots & \what{\mathcal{C}}_{2,K}\\
\vdots                   & \vdots                   &       & \vdots\\
\what{\mathcal{C}}_{K,1} & \what{\mathcal{C}}_{K,2} &\ldots &\what{\mathcal{C}}_{K,K}
\end{pmatrix},\quad
\what{\mathcal{C}}_{k,k'}=
\begin{pmatrix}
\what{\mathcal{C}}_{0,0} & \what{\mathcal{C}}_{0,1} &\ldots & 0\\
\what{\mathcal{C}}_{1,0} & \what{\mathcal{C}}_{1,1} &\ldots & 0\\
0                        &                        0 &\ldots & 0\\
\vdots                   & \vdots                   &       & \vdots\\
0                        &                        0 &\ldots & 0
\end{pmatrix}^{(k,k')}
=
\begin{pmatrix}
{\psi}_{k,k',lb} & {\psi}_{k,k',lt} &\ldots & 0\\
{\psi}_{k,k',rb} & {\psi}_{k,k',rt} &\ldots & 0\\
0                        &                        0 &\ldots & 0\\
\vdots                   & \vdots                   &\dots  & \vdots\\
0                        &                        0 &\ldots & 0
\end{pmatrix}.
\end{equation}
where the auxiliary functions ${\psi}_{k,k',a_1a_2}$ and 
${\psi}_{k,k',a_1a_2}$ are defined in Sec.~\ref{sec:CT-NP}.
Finally, identifying the correspondence with the partial fractions 
of the Pad\'e approximant, the ODEs satisfied by the auxiliary field 
$\mathcal{C}_{k,k'}$ work out to be
\begin{equation}
\begin{split}
\partial_{\tau_2}\what{\mathcal{C}}_{k,k'}(\tau_1,\tau_2)
=-\eta_{k'}^2\what{\mathcal{C}}_{k,k'}(\tau_1,\tau_2)
+\what{\mathcal{A}}_{k,1} \Lambda_2,\quad
\partial_{\tau_1}\what{\mathcal{C}}_{k,k'}(\tau_1,\tau_2)
=-\eta_{k}^2\what{\mathcal{C}}_{k,k'}(\tau_1,\tau_2)
+\Lambda_1\what{\mathcal{A}}_{k',2}.
\end{split}
\end{equation}
Introducing the block matrix $\what{\mathcal{C}}(\tau_1,\tau_2)$, the above 
equations can be stated as follows:
\begin{equation}\label{eq:odes-aux-corners}
\begin{split}
\partial_{\tau_2}\what{\mathcal{C}}(\tau_1,\tau_2)
=-\what{\mathcal{C}}(\tau_1,\tau_2)(\mathcal{E}_K\otimes I_2)
+\vv{1}^{\tp}_K\otimes (\what{\mathcal{A}}_{1} \Lambda_2),\quad
\partial_{\tau_1}\what{\mathcal{C}}(\tau_1,\tau_2)
=-(\mathcal{E}_K\otimes I_1) \what{\mathcal{C}}(\tau_1,\tau_2)
+\vv{1}_K\otimes(\Lambda_1\what{\mathcal{A}}_{2}).
\end{split}
\end{equation}

%=============================================================================%
Employing the approximations developed in~\eqref{eq:npade-apprxm-corners} 
in the (non-local) dynamical system stated in~\eqref{eq:ivp-pade-auxi-mat2}, 
we obtain an effectively local system given by 
\begin{equation}\label{eq:ls-auxi-npade}
\begin{split}
&i \partial_{\tau_2}\what{\mathcal{A}}_{1}M_2
=J_2^{-2}\what{\mathcal{A}}_{1}S_2
+e^{-i\pi/4}(J_2)^{-1} \left[
b_0\what{\mathcal{A}}_{1}\Lambda_2
-\what{\mathcal{C}}(\tau_1,\tau_2)(\vv{b}\otimes I_2)
\right],\\
&i \partial_{\tau_1}M_1\what{\mathcal{A}}_{2}
=J_1^{-2}S_1\what{\mathcal{A}}_{2}
+e^{-i\pi/4}(J_1)^{-1} \left[
b_0\Lambda_1\what{\mathcal{A}}_{2} -(\vv{b}^{\tp}\otimes I_1)\what{\mathcal{C}}(\tau_1,\tau_2)
\right].
\end{split}
\end{equation}
The approach presented here is referred to as the \emph{novel Pad\'e} 
approach. Each of the methods discussed below is derived by applying a 
particular time-stepping scheme to the novel Pad\'e approach, 
therefore, we label them as `NP--' followed by the acronym for 
the time-stepping method.
%=============================================================================%
%=============================================================================%    
\subsubsection{NP--BDF1}
The discrete numerical scheme for the dynamical system defined 
in~\eqref{eq:linear-sys-npade} and~\eqref{eq:ls-auxi-npade} is labelled as 
`NP--BDF1' if the underlying time-stepping method is BDF1. 

Recalling $\rho=1/(\Delta t)$, the BDF1-based discretization of the dynamical 
system in~\eqref{eq:linear-sys-npade} reads as
\begin{equation}\label{eq:bdf1-main-np}
\begin{split}
M_1\what{U}^{j+1}M_2 - M_1\what{U}^{j}M_2
&=-\alpha_1^{-2}S_1\what{U}^{j+1}M_2 - \alpha_2^{-2} M_1\what{U}^{j+1}S_2
-\ovl{b}_0\left[\alpha_1^{-1}\Lambda_1\what{U}^{j+1}M_2 +\alpha_2^{-1}M_1\what{U}^{j+1}\Lambda_2\right]\\
&\quad +\left[\alpha_1^{-1}(\ovl{\vv{b}}^{\tp}\otimes I_1)\what{\mathcal{A}}^{j+1,j+1}_{1}M_2
+\alpha^{-1}_2 M_1\what{\mathcal{A}}^{j+1,j+1}_{2}(\ovl{\vv{b}}\otimes I_2)\right],
\end{split}
\end{equation}
and that of the auxiliary fields in~\eqref{eq:odes-auxi-npade} reads as
\begin{equation}\label{eq:bdf1-aux-np}
\begin{split}
\what{\mathcal{A}}^{j+1,j+1}_{1}
&=(\mathcal{G}_K\otimes I_1)
\left[\frac{1}{\rho}\vv{1}_K\otimes(\Lambda_1\what{U}^{j+1})+\what{\mathcal{A}}^{j,j+1}_{1}\right]
=\left[\frac{1}{\rho}\left(\mathcal{G}_K\vv{1}_K\right)\otimes(\Lambda_1\what{U}^{j+1})
+\left(\mathcal{G}_K\otimes I_1\right)\what{\mathcal{A}}^{j,j+1}_{1}\right],\\
\what{\mathcal{A}}^{j+1,j+1}_{2}
&=\left[\frac{1}{\rho}\vv{1}^{\tp}_K\otimes(\what{U}^{j+1}\Lambda_2)
+\what{\mathcal{A}}^{j+1,j}_{2}\right](\mathcal{G}_K\otimes I_2)
=\left[\frac{1}{\rho}\left(\mathcal{G}_K\vv{1}_K\right)^{\tp}\otimes(\what{U}^{j+1}\Lambda_2)
+\what{\mathcal{A}}^{j+1,j}_{2}\left(\mathcal{G}_K\otimes I_2\right)\right].
\end{split}
\end{equation}
Plugging-in the updates for the auxiliary fields from~\eqref{eq:bdf1-aux-np} 
in the discrete system~\eqref{eq:bdf1-main-np} to obtain 
\begin{equation}
\begin{split}
M_1\what{U}^{j+1}M_2 
+\alpha_1^{-2}S_1\what{U}^{j+1}M_2 &+ \alpha_2^{-2} M_1\what{U}^{j+1}S_2
+\ovl{b}_0\left[\alpha_1^{-1}\Lambda_1\what{U}^{j+1}M_2 +\alpha_2^{-1}M_1\what{U}^{j+1}\Lambda_2\right]\\
&= M_1\what{U}^{j}M_2 
+\alpha_1^{-1} \left[
\left(\frac{1}{\rho}\vv{b}^{\tp}\mathcal{G}_K\vv{1}_K\otimes\Lambda_1\what{U}^{j+1}\right)M_2
+\left(\vv{b}^{\tp}\mathcal{G}_K\otimes I_1\right)\what{\mathcal{A}}^{j,j+1}_{1}M_2\right]\\
&\quad+\alpha_2^{-1} \left[ M_1
\left(\frac{1}{\rho}\left(\mathcal{G}_K\vv{1}_K\right)^{\tp}\vv{b}\otimes(\what{U}^{j+1}\Lambda_2)\right)
+M_1\what{\mathcal{A}}^{j+1,j}_{2}\left(\mathcal{G}_K\vv{b}\otimes I_2\right)
\right].\\
\end{split}
\end{equation}
%=============================================================================%
Introducing $\varpi^{(1/2)}$ as defined in~\eqref{eq:varpi-pm1o2} and 
collecting the terms containing $\what{U}^{j+1}$, the resulting discrete system 
becomes
\begin{equation}\label{eq:sys-main-1}
M_1\what{U}^{j+1}M_2 
+\alpha_1^{-2}S_1\what{U}^{j+1}M_2 + \alpha_2^{-2} M_1\what{U}^{j+1}S_2
+\varpi^{(1/2)} \left[\alpha_1^{-1}\Lambda_1\what{U}^{j+1}M_2 
+\alpha_2^{-1}M_1\what{U}^{j+1}\Lambda_2\right]
=M_1\what{U}^{j}M_2 +\mathcal{B}^{j+1},
\end{equation}
where the history term is given by
\begin{equation}
\mathcal{B}^{j+1} =
\alpha_1^{-1} \left(\vv{b}^{\tp}\mathcal{G}_K\otimes I_1\right)\what{\mathcal{A}}^{j,j+1}_{1}M_2
+ \alpha_2^{-1}M_1\what{\mathcal{A}}^{j+1,j}_{2}
\left(\mathcal{G}_K\vv{b}\otimes I_2\right).\\
\end{equation}
%=============================================================================%
\begin{rem}\label{rem:np-bdf1-hist}
As noted in Remark~\ref{rem:aux-fns-hf}, only the non-zero entries of 
$\what{\mathcal{A}}_{1}$ and $\what{\mathcal{A}}_{2}$ are relevant for the 
numerical implementation so that we may rewrite the history terms 
defined above for each of the boundary segments separately as follows:
\begin{equation*}
\what{\vv{\mathcal{B}}}^{j+1}_{a_1} 
= \sum_{k=1}^K\Gamma^{(1/2)}_k\what{\vv{\varphi}}^{j,j+1}_{k,a_1},\quad
\what{\vv{\mathcal{B}}}^{j+1}_{a_2} 
= \sum_{k=1}^K\Gamma^{(1/2)}_k\what{\vv{\varphi}}^{j+1,j}_{k,a_2},
\end{equation*}
where the auxiliary fields are as defined in~\eqref{eq:aux-np}.
It is straightforward to verify that the history term, $\mathcal{B}^{j+1}$ 
, can now be efficiently computed as
\begin{equation*}
\mathcal{B}^{j+1} = \frac{1}{\alpha_1}\vv{e}^{(1)}_0\otimes
\left[\what{\vs{\mathcal{B}}}^{j+1}_{l}\right]^{\tp}M_2
+\frac{1}{\alpha_1}\vv{e}^{(1)}_1\otimes
\left[\what{\vs{\mathcal{B}}}^{j+1}_{r}\right]^{\tp}M_2
+\frac{1}{\alpha_2}
\left[M_1\what{\vs{\mathcal{B}}}^{j+1}_{b}\right]\otimes\left(\vv{e}_0^{(2)}\right)^{\tp}
+\frac{1}{\alpha_2}
\left[M_1\what{\vs{\mathcal{B}}}^{j+1}_{t}\right]\otimes\left(\vv{e}^{(2)}_1\right)^{\tp}.
\end{equation*}
where $\vv{e}^{(p)}_0=(\delta_{k,0})_{k\in\field{J}_p}$ and 
$\vv{e}^{(p)}_1=(\delta_{k,1})_{k\in\field{J}_p}$ for $p=1,2$ are the orthogonal 
unit vectors in $\field{R}^{N_p+1}$.
\end{rem}
%=============================================================================%
The computation of the history term above requires the off-diagonal samples 
$\what{\mathcal{A}}^{j,j+1}_{1}$ and $\what{\mathcal{A}}^{j+1,j}_{2}$ which can 
be obtained by advancing the corresponding IVPs satisfied by these auxiliary fields.
The BDF1-based discretization of the dynamical system in~\eqref{eq:ls-auxi-npade} reads as
\begin{equation}\label{eq:ls-auxi-np-disc}
\begin{split}
&\left(M_1+\alpha_1^{-2}S_1\right)\what{\mathcal{A}}^{j+1,j}_{2}
=M_1\what{\mathcal{A}}^{j,j}_{2}
-\alpha_1^{-1} \left[ b_0\Lambda_1\what{\mathcal{A}}^{j+1,j}_{2} 
-(\vv{b}^{\tp}\otimes I_1)\what{\mathcal{C}}^{j+1,j} \right],\\
& \what{\mathcal{A}}^{j,j+1}_{1}\left(M_2  +\alpha_2^{-2}S_2\right)
= \what{\mathcal{A}}^{j,j+1}_{1}M_2 -\alpha_2^{-1} \left[
b_0\what{\mathcal{A}}^{j,j+1}_{1}\Lambda_2
-\what{\mathcal{C}}^{j,j+1}(\vv{b}\otimes I_2)
\right],
\end{split}
\end{equation}
and that of the auxiliary fields in~\eqref{eq:odes-aux-corners} work out to be 
\begin{equation}
\begin{split}
\what{\mathcal{C}}^{j+1,j}
&=(\mathcal{G}_K\otimes I_1)
\left[\frac{1}{\rho}\vv{1}_K\otimes(\Lambda_1\what{\mathcal{A}_2}^{j+1,j})+\what{\mathcal{C}}^{j,j}\right]
=\left[\frac{1}{\rho}\left(\mathcal{G}_K\vv{1}_K\right)\otimes(\Lambda_1\what{\mathcal{A}_2}^{j+1,j})
+\left(\mathcal{G}_K\otimes I_1\right)\what{\mathcal{C}}^{j,j}\right],\\
\what{\mathcal{C}}^{j,j+1}
&=\left[\frac{1}{\rho}\vv{1}^{\tp}_K\otimes(\what{\mathcal{A}_1}^{j,j+1}\Lambda_2)
+\what{\mathcal{C}}^{j,j}\right](\mathcal{G}_K\otimes I_2) 
=\left[\frac{1}{\rho}\left(\mathcal{G}_K\vv{1}_K\right)^{\tp}\otimes(\what{\mathcal{A}_1}^{j,j+1}\Lambda_2)
+\what{\mathcal{C}}^{j,j}\left(\mathcal{G}_K\otimes I_2\right)\right].
\end{split}
\end{equation}
Plugging-in these updates in the discrete system~\eqref{eq:ls-auxi-np-disc} and 
collecting the latest time-step of the field yields
\begin{equation}\label{eq:sys-aux-1}
\begin{split}
\Gamma_{a_2}:\quad
&\left(M_1 +\alpha_1^{-2}S_1+\alpha_1^{-1}\varpi^{(1/2)}\Lambda_1\right)\what{\mathcal{A}}^{j+1,j}_{2}
=  M_1\what{\mathcal{A}}^{j,j}_{2}
-\alpha_1^{-1}\left(\vv{b}^{\tp}\mathcal{G}_K\otimes I_1\right)\what{\mathcal{C}}^{j,j},\\
\Gamma_{a_1}:\quad
& \what{\mathcal{A}}^{j,j+1}_{1}
\left(M_2  +\alpha_2^{-2}S_2+\alpha_2^{-1}\varpi^{(1/2)}\Lambda_2\right)
= \what{\mathcal{A}}^{j,j}_{1}M_2 
-\alpha_2^{-1} \what{\mathcal{C}}^{j,j}\left(\mathcal{G}_K\vv{b}\otimes I_2\right).
\end{split}
\end{equation}
The computation of the history sums require diagonal-to-diagonal update of the 
field $\what{\mathcal{C}}^{j,j}$ which can be achieved via 
$\what{\mathcal{C}}^{j,j}\rightarrow\what{\mathcal{C}}^{j,j+1}\rightarrow\what{\mathcal{C}}^{j+1,j+1}$.
The BDF1-based discretization of the auxiliary field in~\eqref{eq:odes-aux-corners} provides: 
\begin{equation}\label{eq:corner-update-bdf1}
\what{\mathcal{C}}^{j+1,j+1}
=(\mathcal{G}_K\otimes I_1)
\left[\frac{1}{\rho}\vv{1}_K\otimes(\Lambda_1\what{\mathcal{A}_2}^{j+1,j+1})+\what{\mathcal{C}}^{j,j+1}\right]
=\left[\frac{1}{\rho}\left(\mathcal{G}_K\vv{1}_K\right)\otimes(\Lambda_1\what{\mathcal{A}_2}^{j+1,j+1})
+\left(\mathcal{G}_K\otimes I_1\right)\what{\mathcal{C}}^{j,j+1}\right].
\end{equation}
%=============================================================================%    
The algorithmic steps for the solving the IBVP in~\eqref{eq:2D-SE-CT} using
NP--BDF1 method are enumerated below.
\begin{enumerate}[%
,leftmargin=*
%,topsep=0mm
%,itemsep=-2mm
%,partopsep=1ex
%,parsep=1ex
,label={\bfseries Step \arabic*:}
]
\item Solve the dynamical system in~\eqref{eq:sys-aux-1} corresponding to auxiliary 
fields for one step, say, $j\rightarrow (j+1)$ using the previously computed values 
of the auxiliary field on the segments (as described in Fig.~\ref{fig:IVP-auxi-pade}).

\item Compute the history functions on each of the boundary segments of the 
domain $\Omega_i$ as in Remark~\ref{rem:np-bdf1-hist} to obtain
\begin{equation*}
\mathcal{B}^{j+1} =
\alpha_1^{-1} \left(\vv{b}^{\tp}\mathcal{G}_K\otimes I_1\right)\what{\mathcal{A}}^{j,j+1}_{1}M_2
+ \alpha_2^{-1}M_1\what{\mathcal{A}}^{j+1,j}_{2}
\left(\mathcal{G}_K\vv{b}\otimes I_2\right).\\
\end{equation*}
\item Solve the discrete linear system in~\eqref{eq:sys-main-1} to compute the solution
at $(j+1)$-th time step, i.e, $\what{U}^{j+1}$. 
\item Update the arrays storing the values of the auxiliary fields 
$\what{\mathcal{A}}_{1}$ and $\what{\mathcal{A}}_{2}$ with the $(j+1)$-th time 
step values using~\eqref{eq:bdf1-aux-np}. The auxiliary field 
$\what{\mathcal{C}}$ at the corners is updated using~\eqref{eq:corner-update-bdf1}. 
\end{enumerate}

%=============================================================================%    
%=============================================================================%    
\subsubsection{NP--TR}
The discrete numerical scheme for the dynamical system defined 
in~\eqref{eq:linear-sys-npade} and~\eqref{eq:ls-auxi-npade} is labelled as `NP--TR' 
if the underlying time-stepping method is TR. 

Recalling $\rho=2/(\Delta t)$, the TR-based discretization of the dynamical 
system in~\eqref{eq:linear-sys-npade} reads as
\begin{equation}\label{eq:tr-main-np}
\begin{split}
M_1\what{U}^{j+1/2}M_2 - M_1\what{U}^{j}M_2
&=-\alpha_1^{-2}S_1\what{U}^{j+1/2}M_2 - \alpha_2^{-2} M_1\what{U}^{j+1/2}S_2
-\ovl{b}_0\left[\alpha_1^{-1}\Lambda_1\what{U}^{j+1/2}M_2 +\alpha_2^{-1}M_1\what{U}^{j+1/2}\Lambda_2\right]\\
&\quad +\left[\alpha_1^{-1}(\ovl{\vv{b}}^{\tp}\otimes I_1)\what{\mathcal{A}}^{j+1/2,j+1/2}_{1}M_2
+\alpha^{-1}_2 M_1\what{\mathcal{A}}^{j+1/2,j+1/2}_{2}(\ovl{\vv{b}}\otimes I_2)\right],
\end{split}
\end{equation}
and that of the auxiliary fields in~\eqref{eq:odes-auxi-npade} reads as
\begin{equation}\label{eq:tr-aux-np}
\begin{split}
\what{\mathcal{A}}^{j+1/2,j+1/2}_{1}
&=\frac{1}{2}\left[ \left(\mathcal{H}_K\otimes I_1\right) \what{\mathcal{A}}^{j,j+1}_{1}
+\what{\mathcal{A}}^{j,j}_{1}\right]+ \left(\mathcal{G}_K\otimes I_1\right) 
\left[\frac{1}{\rho}\vv{1}_K\otimes \left(\Lambda_1\what{U}^{j+1/2}
+\frac{\what{\Phi}^{j,j+1}_{1}-\what{\Phi}^{j,j}_{1}}{2}\right)\right],\\
\what{\mathcal{A}}^{j+1/2,j+1/2}_{2}
&=\frac{1}{2}\left[ \what{\mathcal{A}}^{j+1,j}_{2}\left(\mathcal{H}_K\otimes I_2\right) 
+\what{\mathcal{A}}^{j,j}_{2}\right]+  
\left[\frac{1}{\rho}\vv{1}^{\tp}_K\otimes \left(\what{U}^{j+1/2}\Lambda_2
+\frac{\what{\Phi}^{j+1,j}_{2}-\what{\Phi}^{j,j}_{2}}{2}\right)\right]
\left(\mathcal{G}_K\otimes I_2\right).
\end{split}
\end{equation}
Plugging-in the updates for the auxiliary fields from~\eqref{eq:tr-aux-np} 
in the discrete system~\eqref{eq:tr-main-np} and simplifying, yields
\begin{equation}\label{eq:sys-main-2}
\begin{split}
M_1\what{U}^{j+1/2}M_2 
+\alpha_1^{-2}S_1\what{U}^{j+1/2}M_2  + \alpha_2^{-2} M_1\what{U}^{j+1/2}S_2
& +\varpi^{(1/2)} \left[\alpha_1^{-1}\Lambda_1\what{U}^{j+1/2}M_2 
+\alpha_2^{-1}M_1\what{U}^{j+1/2}\Lambda_2\right]\\
& = M_1\what{U}^{j}M_2 + \mathcal{B}_{a_1}^{j+1/2} +\mathcal{B}_{a_2}^{j+1/2},
\end{split}
\end{equation}
where the history terms are given by
\begin{equation}
\begin{split}
& \mathcal{B}_{a_1}^{j+1/2} =\alpha_1^{-1}
\frac{1}{\rho}\left(\vv{b}^{\tp}\mathcal{G}_K\vv{1}_K\right)
\left(\frac{\what{\Phi}^{j+1,j}_{1}-\what{\Phi}^{j,j}_{1}}{2}\right)M_2
+ \frac{1}{2}\alpha_1^{-1}\left[
\left(\vv{b}^{\tp}\mathcal{H}_K\otimes I_1\right)\what{\mathcal{A}}^{j,j+1}_{1}
+\left(\vv{b}^{\tp}\otimes I_1\right)\what{\mathcal{A}}^{j,j}_{1}
\right]M_2,\\
&\mathcal{B}_{a_2}^{j+1/2} = \alpha_2^{-1}
\frac{1}{\rho}\left(\vv{b}^{\tp}\mathcal{G}_K\vv{1}_K\right)
 M_1\left(\frac{\what{\Phi}^{j+1,j}_{2}-\what{\Phi}^{j,j}_{2}}{2}\right)
+\frac{1}{2}\alpha_2^{-1}M_1\left[
\what{\mathcal{A}}^{j+1,j}_{2}\left(\mathcal{H}_K\vv{b}\otimes I_2\right)
+\what{\mathcal{A}}^{j+1,j}_{2}\left(\vv{b}\otimes I_2\right)\right].
\end{split}
\end{equation}
The simplifications noted in Remark~\ref{rem:np-bdf1-hist} can also be applied 
here for the history terms. We omit the details because the adaption is 
straightforward. 
%=============================================================================%
The computation of the history terms require diagonal-to-diagonal update of 
fields $\what{\mathcal{A}}^{j,j}_{1}$ and $ \what{\mathcal{A}}^{j,j}_{2}$ which 
can be obtained by TR-based discretization of~\eqref{eq:odes-auxi-npade} as follows:
\begin{equation}\label{eq:tr-aux2-np}
\begin{split}
\what{\mathcal{A}}^{j+1,j+1}_{1}
&=\left(\mathcal{H}_K\otimes I_1\right) \what{\mathcal{A}}^{j,j+1}_{1}
+ \left[\frac{2}{\rho}\left(\mathcal{G}_K\vv{1}_K\right)\otimes
\left(\Lambda_1\what{U}^{j+1/2} +\frac{\what{\Phi}^{j,j+1}_{1}-\what{\Phi}^{j,j}_{1}}{2}\right) \right],\\
\what{\mathcal{A}}^{j+1,j+1}_{2}
&=\what{\mathcal{A}}^{j+1,j}_{2}\left(\mathcal{H}_K\otimes I_2\right) 
+ \left[\frac{2}{\rho}\left(\mathcal{G}_K\vv{1}_K\right)^{\tp}\otimes
\left(\what{U}^{j+1/2}\Lambda_2 +\frac{\what{\Phi}^{j+1,j}_{2}-\what{\Phi}^{j,j}_{2}}{2}\right) \right].
\end{split}
\end{equation}
Moreover, the updates above require off-diagonal samples of the auxiliary fields 
which can be obtained by solving the dynamical system in~\eqref{eq:ls-auxi-npade}. 
The TR-based discretization of this system reads as
\begin{equation}\label{eq:ls-auxi-np-disc-tr}
\begin{split}
&\left(M_1+\alpha_1^{-2}S_1\right)\what{\mathcal{A}}^{j+1/2,j}_{2}
 =  M_1\what{\mathcal{A}}^{j,j}_{2}
-\alpha_1^{-1} \left[ b_0\Lambda_1\what{\mathcal{A}}^{j+1/2,j}_{2} 
-(\vv{b}^{\tp}\otimes I_1)\what{\mathcal{C}}^{j+1/2,j} \right],\\
& \what{\mathcal{A}}^{j,j+1/2}_{1}\left(M_2  +\alpha_2^{-2}S_2\right)
= \what{\mathcal{A}}^{j,j}_{1}M_2 -\alpha_2^{-1} \left[
b_0\what{\mathcal{A}}^{j,j+1/2}_{1}\Lambda_2
-\what{\mathcal{C}}^{j,j+1/2}(\vv{b}\otimes I_2)
\right],
\end{split}
\end{equation}
and that of the auxiliary field in~\eqref{eq:odes-aux-corners} work out to be 
\begin{equation}
\begin{split}
\what{\mathcal{C}}^{j+1/2,j}
&=
\left[\frac{1}{\rho}\left(\mathcal{G}_K\vv{1}_K\right)\otimes(\Lambda_1\what{\mathcal{A}_2}^{j+1/2,j})
+\left(\mathcal{H}_K\otimes I_1\right)\what{\mathcal{C}}^{j,j}\right],\\
\what{\mathcal{C}}^{j,j+1/2}
&=
\left[\frac{1}{\rho}\left(\mathcal{G}_K\vv{1}_K\right)^{\tp}\otimes(\what{\mathcal{A}_1}^{j,j+1/2}\Lambda_2)
+\what{\mathcal{C}}^{j,j}\left(\mathcal{H}_K\otimes I_2\right)\right].
\end{split}
\end{equation}
%=============================================================================%
Plugging-in these updates in the discrete system~\eqref{eq:ls-auxi-np-disc-tr} 
and collecting the latest time-step of the field yields
\begin{equation}\label{eq:sys-aux-2}
\begin{split}
\Gamma_{a_2}:\quad
&\left(M_1 +\alpha_1^{-2}S_1+\alpha_1^{-1} \varpi^{(1/2)}\Lambda_1\right)\what{\mathcal{A}}^{j+1/2,j}_{2}
=  M_1\what{\mathcal{A}}^{j,j}_{2}
+\alpha_1^{-1}\left(\vv{b}^{\tp}\mathcal{H}_K\otimes I_1\right)\what{\mathcal{C}}^{j,j}, \\
\Gamma_{a_1}:\quad
& \what{\mathcal{A}}^{j,j+1/2}_{1}\left(M_2+\alpha_2^{-2}S_2
+\alpha_2^{-1}\varpi^{(1/2)}\Lambda_2\right)
= \what{\mathcal{A}}^{j,j}_{1}M_2 
+\alpha_2^{-1} \what{\mathcal{C}}^{j,j}\left(\mathcal{H}_K\vv{b}\otimes I_2\right).
\end{split}
\end{equation}
The computation of the history terms above will require diagonal-to-diagonal update of 
auxiliary field $\what{\mathcal{C}}^{j,j}$ which can be achieved via 
$\what{\mathcal{C}}^{j,j}\rightarrow\what{\mathcal{C}}^{j,j+1}\rightarrow\what{\mathcal{C}}^{j+1,j+1}$.
The TR-based discretization of the systems~\eqref{eq:odes-aux-corners} leads to
\begin{equation}\label{eq:corner-update-tr}
\begin{split}
&\what{\mathcal{C}}^{j,j+1}
=\left[\frac{2}{\rho}\left(\mathcal{G}_K\vv{1}_K\right)^{\tp}\otimes(\what{\mathcal{A}_1}^{j,j+1/2}\Lambda_2)
+\what{\mathcal{C}}^{j,j}\left(\mathcal{H}_K\otimes I_2\right)\right],\\
&\what{\mathcal{C}}^{j+1,j+1}
=\left[\frac{2}{\rho}\left(\mathcal{G}_K\vv{1}_K\right)\otimes
\frac{1}{2}\left(\Lambda_1\what{\mathcal{A}_2}^{j+1,j+1}+\Lambda_1\what{\mathcal{A}_2}^{j,j+1} \right)
+\left(\mathcal{H}_K\otimes I_1\right)\what{\mathcal{C}}^{j,j+1}\right].
\end{split}
\end{equation}
The computation of history terms $\mathcal{B}_{a_1}^{j+1/2}$ and $\mathcal{B}_{a_2}^{j+1/2}$
also require the off-diagonal samples $\what{\Phi}^{j,j+1}_{1}$ and $\what{\Phi}^{j+1,j}_{2}$ 
which can be obtained by solving the dynamical system in~\eqref{eq:linear-sys-ivps}. The 
TR-based discretization of this system reads as
\begin{equation}\label{eq:ls-disc-phi}
\begin{split}
&\left(M_1+\alpha_1^{-2}S_1\right)\what{\Phi}^{j+1/2,j}_{2}
 =  M_1\what{\Phi}^{j,j}_{2}
-\alpha_1^{-1} \left[ b_0\Lambda_1\what{\Phi}^{j+1/2,j}_{2} 
-(\vv{b}^{\tp}\otimes I_1)\what{\mathcal{A}}^{j+1/2,j}_1\Lambda_2 \right],\\
&  \what{\Phi}^{j,j+1/2}_{1}\left(M_2  +\alpha_2^{-2}S_2\right)
=  \what{\Phi}^{j,j}_{1}M_2 -\alpha_2^{-1} \left[
b_0\what{\Phi}^{j,j+1/2}_{1}\Lambda_2
-\Lambda_1\what{\mathcal{A}}^{j,j+1/2}_2(\vv{b}\otimes I_2)
\right],
\end{split}
\end{equation}
and that of the auxiliary fields in~\eqref{eq:odes-auxi-npade} reads as
\begin{equation}
\begin{split}
\what{\mathcal{A}}_1^{j+1/2,j}
&=
\left[\frac{1}{\rho}\left(\mathcal{G}_K\vv{1}_K\right)\otimes(\Lambda_1\what{\Phi}^{j+1/2,j})
+\left(\mathcal{G}_K\otimes I_1\right)\what{\mathcal{A}}_1^{j,j}\right],\\
\what{\mathcal{A}}_2^{j,j+1/2}
&=
\left[\frac{1}{\rho}\left(\mathcal{G}_K\vv{1}_K\right)^{\tp}\otimes(\what{\Phi}^{j,j+1/2}\Lambda_2)
+\what{\mathcal{A}}_2^{j,j}\left(\mathcal{G}_K\otimes I_2\right)\right].
\end{split}
\end{equation}
%=============================================================================%
Plugging-in these updates in the discrete system~\eqref{eq:ls-disc-phi} and 
collecting the latest time-step of the field yields
\begin{equation}\label{eq:sys-aux-3}
\begin{split}
\Gamma_{a_2}:\quad
&\left(M_1 +\alpha_1^{-2}S_1+\alpha_1^{-1}\varpi^{(1/2)}\Lambda_1\right)\what{\Phi}^{j+1/2,j}_{2}
=M_1\what{\Phi}^{j,j}_{2}
+\alpha_1^{-1}\left(\vv{b}^{\tp}\mathcal{G}_K\otimes I_1\right)\what{\mathcal{A}}_1^{j,j}\Lambda_2, \\
\Gamma_{a_1}:\quad
& \what{\Phi}^{j,j+1/2}_{1}\left(M_2+\alpha_2^{-2}S_2
+\alpha_2^{-1}\varpi^{(1/2)}\Lambda_2\right)
= \what{\Phi}^{j,j}_{1}M_2 
+\alpha_2^{-1}\Lambda_1\what{\mathcal{A}}_2^{j,j}
\left(\mathcal{G}_K\vv{b}\otimes I_2\right).
\end{split}
\end{equation}
%=============================================================================%    
The algorithmic steps for the solving the IBVP in~\eqref{eq:2D-SE-CT} using
NP--TR method are enumerated below.
\begin{enumerate}[%
,leftmargin=*
%,topsep=0mm
%,itemsep=-2mm
%,partopsep=1ex
%,parsep=1ex
,label={\bfseries Step \arabic*:}
]
\item Solve the dynamical system in~\eqref{eq:sys-aux-2} and~\eqref{eq:sys-aux-3} 
corresponding to auxiliary fields for one step, say, $j\rightarrow (j+1)$ using 
the previously computed values of the auxiliary field on the segments 
(as described in Fig.~\ref{fig:IVP-auxi-pade}).
\item Compute the history functions on each of the boundary segments of the 
domain $\Omega_i$ to obtain
\begin{equation*}
\begin{split}
& \mathcal{B}_{a_1}^{j+1/2} =\alpha_1^{-1}
\frac{1}{\rho}\left(\vv{b}^{\tp}\mathcal{G}_K\vv{1}_K\right)
\left(\frac{\what{\Phi}^{j+1,j}_{1}-\what{\Phi}^{j,j}_{1}}{2}\right)M_2
+ \frac{1}{2}\alpha_1^{-1}\left[
\left(\vv{b}^{\tp}\mathcal{H}_K\otimes I_1\right)\what{\mathcal{A}}^{j,j+1}_{1}
+\left(\vv{b}^{\tp}\otimes I_1\right)\what{\mathcal{A}}^{j,j}_{1}
\right]M_2,\\
&\mathcal{B}_{a_2}^{j+1/2} = \alpha_2^{-1}
\frac{1}{\rho}\left(\vv{b}^{\tp}\mathcal{G}_K\vv{1}_K\right)
 M_1\left(\frac{\what{\Phi}^{j+1,j}_{2}-\what{\Phi}^{j,j}_{2}}{2}\right)
+\frac{1}{2}\alpha_2^{-1}M_1\left[
\what{\mathcal{A}}^{j+1,j}_{2}\left(\mathcal{H}_K\vv{b}\otimes I_2\right)
+\what{\mathcal{A}}^{j+1,j}_{2}\left(\vv{b}\otimes I_2\right)\right].
\end{split}
\end{equation*}
\item Solve the discrete linear system in~\eqref{eq:sys-main-2} to compute the solution
at $(j+1)$-th time step, i.e, $\what{U}^{j+1/2}$. 
\item Update the arrays storing the values of the auxiliary fields 
$\what{\mathcal{A}}_{1}$ and $\what{\mathcal{A}}_{2}$ with the $(j+1)$-th time 
step values using~\eqref{eq:tr-aux2-np}. The auxiliary field 
$\what{\mathcal{C}}$ at the corners is updated using~\eqref{eq:corner-update-tr}. 
\end{enumerate}

%=============================================================================%    
%=============================================================================%    
\section{Numerical Experiments}\label{sec:numerical-experiments}
In this section, we will present several numerical tests to analyze the 
accuracy of the numerical schemes developed in this work to solve the 
IBVP~\eqref{eq:2D-SE-CT}. These numerical tests consist of analyzing the error 
evolution behaviour and empirical convergence of the considered IBVP. In order to
visualize the reflections in the numerical solution, we also present some contour 
plots at different instants of time.

%=============================================================================%
%=============================================================================%
\subsection{Exact solutions}\label{sec:exact-solution}
The exact solutions admissible for our numerical experiments are such that the 
initial profile must be effectively supported within the computational domain.
We primarily consider the wavepackets which are a modulation of a Gaussian
envelop so that the requirement of effective initial support can be easily met.
In this class of solutions, we consider the Chirped-Gaussian and the
Hermite-Gaussian profiles.

%=============================================================================%
\def\arraystretch{2}
\setlength{\tabcolsep}{1mm}
\begin{table}[!htbp]
\centering
\caption{\label{tab:cg2d} Chirped-Gaussian profile with $A_0=2$ and 
$c_0\in\{4,8,12,16\}$.}
\begin{tabular}{m{10mm}m{5mm}m{30mm}m{20mm}m{40mm}}
\hline
\multicolumn{5}{c}{
$G\left(\vv{x},t\right) 
= A_0\sum_{j=1}^n G(\vv{x},t;\vv{a}_j,\vv{b}_j,\vv{c}_j),
\quad\vv{c}_j=c_0(\cos\theta_j,\sin\theta_j)$}\\
\hline
Type &$n$ & $(\vv{a}_j\in\field{R}^2)_{j=1}^n$
&$(\vv{b}_j\in\field{R}^2)_{j=1}^n$ & $\vs{\theta}\in\field{R}^n$\\
\hline
IIA,\newline{IIB} & $4$ & 
$\vv{a}_1 = (1/2.5,1/2.4)$,\newline 
$\vv{a}_2 = (1/2.3,1/2.2)$,\newline
$\vv{a}_3 = (1/2.7,1/2.6)$,\newline
$\vv{a}_4 = (1/2.2,1/2.5)$
& $\vv{b}_j=(1/2)\vv{1}$
&IIA:~$\vs{\theta}_A=(0,\pi/2,\pi,3\pi/2)$,\newline{IIB}:~$\vs{\theta}_B=\vv{\theta}_A+(\pi/4)\vv{1}$\\
\hline
\end{tabular}
\end{table}
%=============================================================================%
%=============================================================================%
%=============================================================================%
\subsubsection{Chirped-Gaussian profile}
\begin{figure}[!htbp]
\begin{center}
\def\myscale{1}
\includegraphics[width=\textwidth]{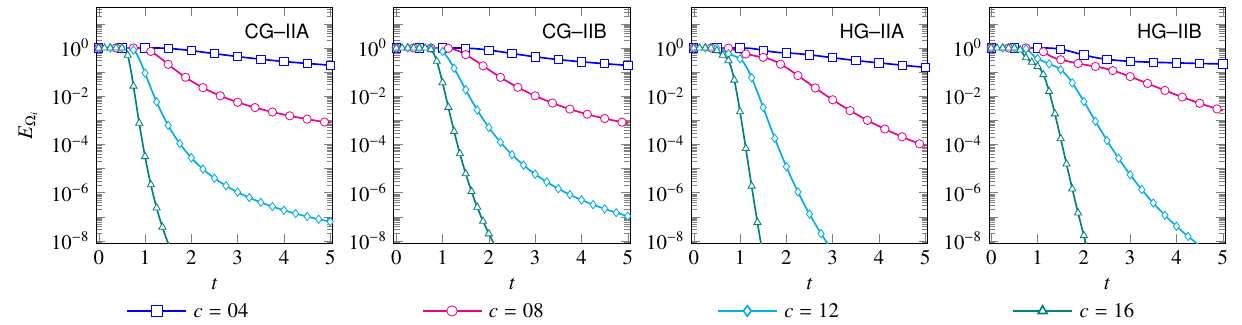}
\end{center}
\caption{\label{fig:wp2d-energy-content}The figure shows the evolution of the 
relative energy content as defined in~\eqref{eq:cg2d-energy-content} of the 
chirped-Gaussian and Hermite-Gaussian profiles considered in Table~\ref{tab:cg2d}
and Table~\ref{tab:hg2d}. Here the
computational domain is $\Omega_i=(-10,10)^2$.}
\end{figure}

Using the function defined by
\begin{equation}
\mathcal{G}(x,t;a,b)=\frac{1}{\sqrt{1+4i(a+ib)t}}
\exp\left[-\frac{(a+ib)}{1+4i(a+ib)t}x^2\right],\quad a>0,
\end{equation}
one can define the family solutions referred to as chirped-Gaussian profile by
\begin{equation}
G(\vv{x},t;\vv{a},\vv{b},\vv{c})
=\mathcal{G}(x_1-c_1t,t;a_1,b_1)\mathcal{G}(x_2-c_2t,t;a_2,b_2)
\exp\left(+i\frac{1}{2}\vv{c}\cdot\vv{x}-i\frac{1}{4}\vv{c}\cdot\vv{c}\,t\right),
\end{equation}
where $\vv{a}\in\field{R}^2_+$ determines the effective support of the 
profile at $t=0$, $\vv{b}\in\field{R}^2$ is the \emph{chirp} parameter and 
$\vv{c}\in\field{R}^2$ is the velocity of the profile. Using linear 
combination, one can further define a more general family of solutions with 
parameters $A_0,c_0\in\field{R}$, $\vs{\theta}\in\field{R}^n$, 
$(\vv{a}_j\in\field{R}^2_+)_{j=1}^n$ and 
$(\vv{b}_j\in\field{R}^2)_{j=1}^n$ given by
\begin{equation}
G\left(\vv{x},t;\;c_0,A_0;\;(\vv{a}_j)_{j=1}^n,(\vv{b}_j)_{j=1}^n,\vv{\theta}\right) 
= A_0\sum_{j=1}^n G(\vv{x},t;\vv{a}_j,\vv{b}_j,\vv{c}_j),
\quad\vv{c}_j=c_0(\cos\theta_j,\sin\theta_j).
\end{equation}
Let the constant vector $(1,1,\ldots)\in\field{R}^n$ be denoted by $\vv{1}$,
then the specific values of the parameters of the solutions used in the
numerical experiments can be summarized as in Table~\ref{tab:cg2d}. The energy 
content of the profile within the computational domain $\Omega_i$ 
over time is
\begin{equation}\label{eq:cg2d-energy-content}
E_{\Omega_i}(t)
=\left.{\int_{\Omega_i}|G(\vv{x},t)|^2d^2\vv{x}}\right/
{\int_{\Omega_i}|G(\vv{x},0)|^2d^2\vv{x}},
\quad t\geq0.
\end{equation}
The profiles are chosen with non-zero speed $c_0$ so that the field hits the
boundary of $\Omega_i$, In case of a square computational domain, the type `A' class of 
solutions are directed to the segments of the boundary normally while type `B' 
class of solutions are directed to the corners. The behavior of the 
$E_{\Omega_i}(t)$ for $t\in[0,5]$ is shown in Fig.~\ref{fig:wp2d-energy-content}.

%=============================================================================%
\def\arraystretch{2}
\setlength{\tabcolsep}{1mm}
\begin{table}[!htbp]
\centering
\caption{\label{tab:hg2d} Hermite-Gaussian profile with $A_0=2$ and 
$c_0\in\{4,8,12,16\}$.}
\begin{tabular}{m{10mm}m{5mm}m{30mm}m{20mm}m{40mm}}
\hline
\multicolumn{5}{c}{
$G\left(\vv{x},t\right) 
= A_0\sum_{j=1}^n G(\vv{x},t;\vv{m}_j,\vv{a}_j,\vv{c}_j),\;
\vv{c}_j=c_0(\cos\theta_j,\sin\theta_j)$
}\\
\hline
Type &$n$ & $(\vv{a}_j\in\field{R}^2)_{j=1}^n$
&$(\vv{m}_j\in\field{N}_0^2)_{j=1}^n$ & $\vs{\theta}\in\field{R}^n$\\
\hline
IIA,\newline{IIB} & $4$ & 
$\vv{a}_1 = (1/2.5,1/2.4)$,\newline 
$\vv{a}_2 = (1/2.3,1/2.2)$,\newline
$\vv{a}_3 = (1/2.7,1/2.6)$,\newline
$\vv{a}_4 = (1/2.2,1/2.5)$ &
$\vv{m}_1 = (1,2)$,\newline
$\vv{m}_2 = (2,1)$,\newline
$\vv{m}_3 = (2,1)$,\newline
$\vv{m}_4 = (1,2)$
&$\vs{\theta}_A=(0,\pi/2,\pi,3\pi/2)$,\newline%
$\vs{\theta}_B=\vv{\theta}_A+(\pi/4)\vv{1}$\\
\hline
\end{tabular}
\end{table}
%=============================================================================%

%=============================================================================%
%=============================================================================%

%=============================================================================%
%=============================================================================%
\subsubsection{Hermite-Gaussian profile}
Consider the class of normalized Hermite-Gaussian functions defined by 
\begin{equation}
\mathcal{G}_{m}(x,t;a)=\gamma_m^{-1}H_{m}\left(\frac{\sqrt{2a}\,x }{w(t)}\right)
\sqrt{\frac{\mu(t)}{a}}
\exp\left[-\mu(t)x^2-i\,m\theta(t)\right],
\quad m\in\field{N}_0=\{0,1,\ldots\},
\end{equation}
where $a>0$,
\begin{equation}
w(t)=\sqrt{1+(4at)^2},\quad
\frac{1}{\mu(t)}=\frac{1}{a}+i{4t}=\frac{1}{a}w(t)\exp[i\theta(t)],
\end{equation}
and the normalization factor is given by 
$\gamma^2_{m}=2^m(m!)\sqrt{\pi}(2a)^{-1/2}$. 
The Hermite polynomials are evaluated using the following relations
\begin{equation}
H_{n+1}(x)=2xH_n(x)-2nH_{n-1}(x),\quad H_{n+1}'(x)=2(n+1)H_{n}(x),
\end{equation}
with $H_0(x) = 1$ and $H_1(x)=2x$. Using these functions, we
can define a family of solutions referred to as Hermite-Gaussian profile by
\begin{equation}
G(\vv{x},t;\vv{m},\vv{a},\vv{c})
=\mathcal{G}_{m_1}(x_1-c_1t,t;a_1)\mathcal{G}_{m_2}(x_2-c_2t,t;a_2)
\exp\left(+i\frac{1}{2}\vv{c}\cdot\vv{x}-i\frac{1}{4}\vv{c}\cdot\vv{c}\,t\right),
\end{equation}
where  $\vv{m}\in\field{N}_0^2$ is the order parameter, 
$\vv{a}\in\field{R}^2_+$ determines the effective support of the profile 
at $t=0$ and $\vv{c}\in\field{R}^2$ is the velocity of the profile. Once again,
using linear combination, one can further define a more general family of solutions with 
parameters $A_0,c_0\in\field{R}$, $\vs{\theta}\in\field{R}^n$, 
$(\vv{m}_j\in\field{N}^2_0)_{j=1}^n$ and
$(\vv{a}_j\in\field{R}^2_+)_{j=1}^n$ given by
\begin{equation}
G\left(\vv{x},t;\;c_0,A_0;\;(\vv{m}_j)_{j=1}^n,(\vv{a}_j)_{j=1}^n,\vv{\theta}\right) 
= A_0\sum_{j=1}^n G(\vv{x},t;\vv{m}_j,\vv{a}_j,\vv{c}_j),
\quad\vv{c}_j=c_0(\cos\theta_j,\sin\theta_j).
\end{equation}
As in the last case, we set $\vv{1}=(1,1,\ldots)\in\field{R}^n$, then the specific 
values of the parameters of the solutions used in the
numerical experiments can be summarized as in Table~\ref{tab:hg2d}. The energy 
content of the profile within the computational domain $\Omega_i$ 
over time is $E_{\Omega_i}(t)$ as defined by~\eqref{eq:cg2d-energy-content} with
the appropriate profile in the integrand. The profiles are chosen with non-zero 
speed $c_0$ so that the field hits the
boundary of $\Omega_i$. In case of a square computational domain, the type `A' 
class of solutions are directed to the segments of the boundary normally while 
type `B' class of solutions are directed to the corners. The behavior of the 
$E_{\Omega_i}(t)$ for $t\in[0,5]$ is shown in Fig.~\ref{fig:wp2d-energy-content}.
%=============================================================================%
\def\arraystretch{2}
\setlength{\tabcolsep}{1mm}
\begin{table}[!htbp]
\centering
\caption{\label{tab:ee-params}Numerical parameters for studying the evolution 
error}
\begin{tabular}{m{80mm}m{50mm}}\hline
Computational domain ($\Omega_i$) & $(-10,10)\times(-10,10)$\\\hline
Maximum time ($T_{max}$)          & $5$\\\hline
No. of time-steps ($N_t$)         & $5000+1$\\\hline
Time-step ($\Delta t$)            & $10^{-3}=T_{max}/(N_t-1)$\\\hline
\>Number of LGL-points ($(N+1)\times (N+1)$) & $200\times 200$\\\hline
\end{tabular}
\end{table}
%=============================================================================%

%=============================================================================%
\def\arraystretch{2}
\setlength{\tabcolsep}{1mm}
\begin{table}[!htbp]
\centering
\caption{\label{tab:ca-params}Numerical parameters for studying the convergence 
error for the HF methods}
\begin{tabular}{m{80mm}m{50mm}}\hline
Computational domain ($\Omega_i$) & $(-10,10)\times(-10,10)$\\\hline
Maximum time ($T_{max}$)          & $5$\\\hline
Set of no. of time-steps ($\field{N}_t$)& $\{2^8,2^9,\ldots,2^{16}\}$\\\hline
Time-step               & $\{\Delta t=T_{max}/(N_t-1),\;N_t\in\field{N}_t\}$\\\hline
\>Number of LGL-points ($(N+1)\times (N+1)$) & $200\times 200$\\\hline
\end{tabular}
\end{table}
%=============================================================================%

%=============================================================================%
%=============================================================================%
\subsection{Tests for evolution error}\label{sec:tests-ee}
In this section, we consider the IBVP~\eqref{eq:2D-SE-CT} where the initial condition 
corresponds to the exact solutions described in Sec.~\ref{sec:exact-solution}. 
In this work, we have developed two discrete versions of the DtN maps, namely, 
HF and NP. Each of these methods have a variant determined by the choice of 
time-stepping method used in the temporal discretization so that the complete list of
methods to be tested can be labelled as NP--BDF1, HF--BDF1 (corresponding to the
one-step method BDF1), HF--BDF2 (corresponding to the
multi-step method BDF2), and, NP--TR, HF--TR (corresponding to the
one-step method TR). Moreover, HF methods have two variants namely, CQ and CP,
depending on the discretization scheme used for time-fractional operators. For all 
Pad\'e based methods, we
use diagonal Pad\'e approximants of order 30.
The numerical parameters used in this section are summarized in Table~\ref{tab:ee-params}. 
The numerical solution is labelled according to that of the boundary maps. The error 
in the evolution of the profile computed numerically is quantified by the 
relative $\fs{L}^2(\Omega_i)$-norm: 
\begin{equation}\label{eq:error-ibvp}
e(t_j)=
\left.\left(\int_{\Omega_i}
\left|u(\vv{x},t_j)-[u(\vv{x},t_j)]_{\text{num.}}
\right|^2d^2\vv{x}\right)^{1/2}\,\right/
\left(\int_{\Omega_i}\left|u_0(\vv{x})\right|^2d^2\vv{x}\right)^{1/2},
\end{equation}
for $t_j\in[0,T_{max}]$. The integral above will be approximated by Gauss 
quadrature over LGL-points. 
%=============================================================================%
\begin{figure}[!h]
\begin{center}
\includegraphics[width=\textwidth]{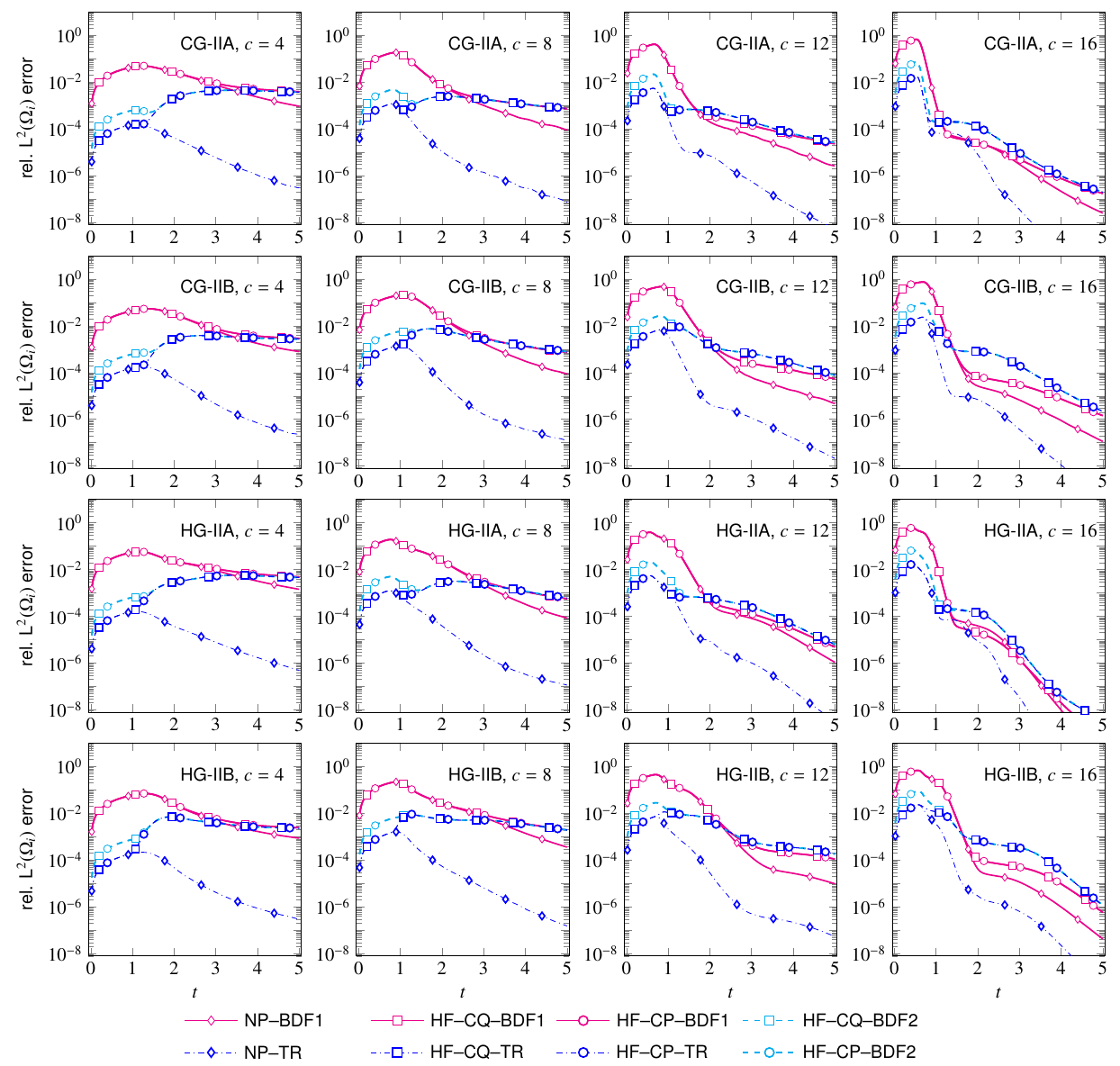}
\end{center}
\caption{\label{fig:wpcghg2d-ee}The figure shows a comparison of evolution of error 
in the numerical solution of the IBVP~\eqref{eq:2D-SE-CT} with various 
approximations of the TBCs for the chirped-Gaussian and Hermite-Gaussian profiles 
with different values of the speed `c' (see Table~\ref{tab:cg2d} and Table~\ref{tab:hg2d}). 
The numerical parameters and the labels are described in 
Sec.~\ref{sec:tests-ee} where the error is quantified by~\eqref{eq:error-ibvp}.}
\end{figure}
%=============================================================================%

The numerical results for the evolution error on $\Omega_i$ corresponding to the 
chirped-Gaussian and Hermite-Gaussian profiles are shown in Fig.~\ref{fig:wpcghg2d-ee}.
Note that the TR based methods perform better than the BDF based methods which is 
clear from the error peaks in Fig.~\ref{fig:wpcghg2d-ee}. Let us observe that BDF1 
method is incapable of distinguishing between NP and HF methods which points towards 
the fact that higher order method is needed to exploit the superior accuracy of NP 
method to that of HF. Also, for some initial time-steps, second order HF methods 
follow NP--TR but soon they become less accurate due to accumulation of error.
Particularly for HF methods, the error behaviour improves for the case of faster 
moving profiles. The performance of NP--TR turns out to be best among all the 
methods considered in this test. There is not much significant difference in the results for 
evolution error on $\Omega_i$ corresponding to the chirped-Gaussian 
profile and Hermite-Gaussian profile.  
%=============================================================================%
%=============================================================================%
\begin{figure}[!h]
\begin{center}
\includegraphics[width=\textwidth]{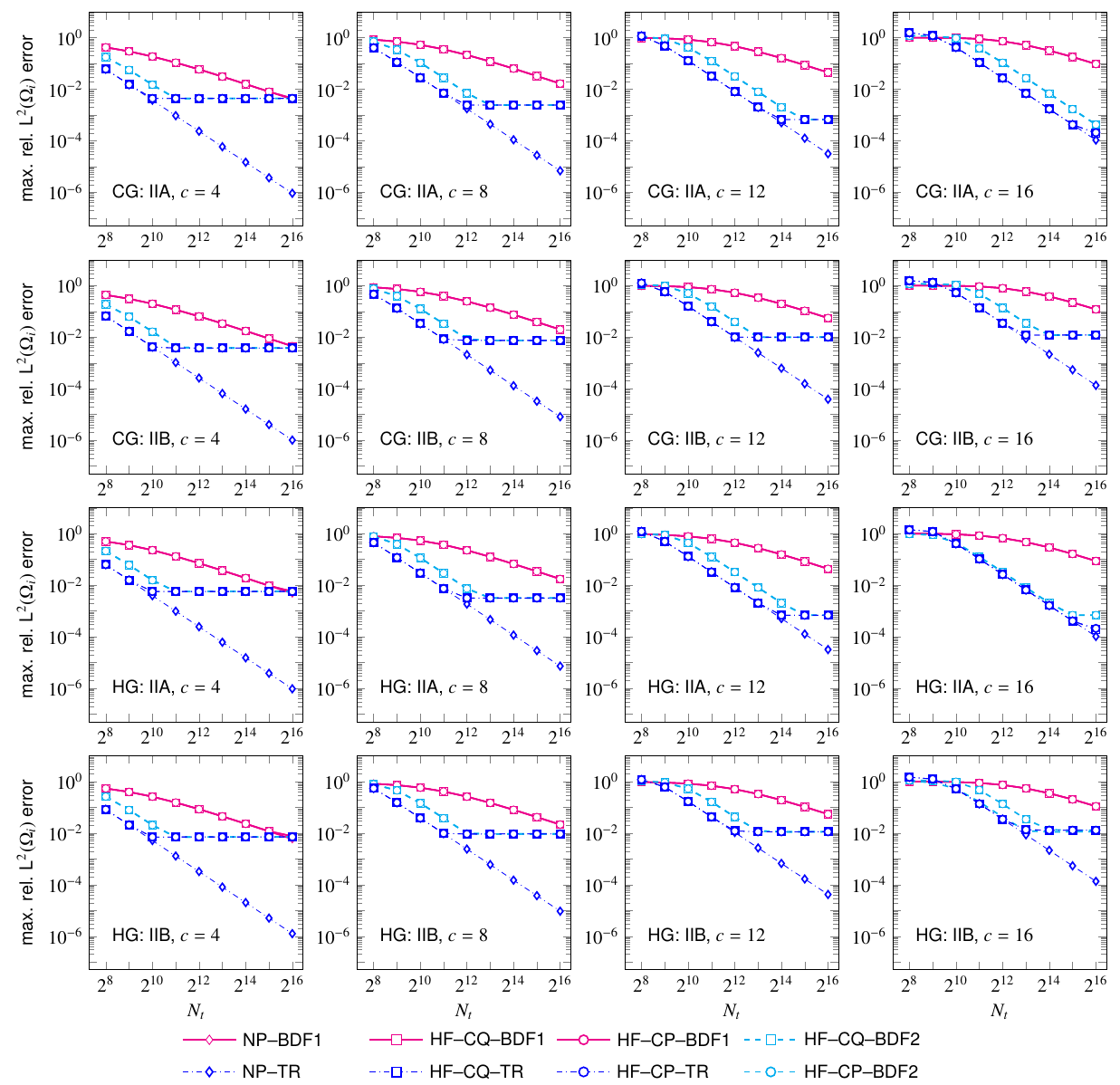}
\end{center}
\caption{\label{fig:wpcghg2d-ca}The figure shows a comparison of convergence
behaviour in the numerical solution of the IBVP~\eqref{eq:2D-SE-CT} with various 
approximations of the TBCs for the chirped-Gaussian profile
and Hermite-Gaussian profile with different values of the speed `c' 
(see Table~\ref{tab:cg2d} and Table~\ref{tab:hg2d}). The numerical parameters 
and the labels are described in Sec.~\ref{sec:tests-ca} where the error is 
quantified by~\eqref{eq:max-error-ibvp}.}
\end{figure}
%=============================================================================%    

\begin{figure}[!htbp]
\begin{center}
\includegraphics[scale=0.95]{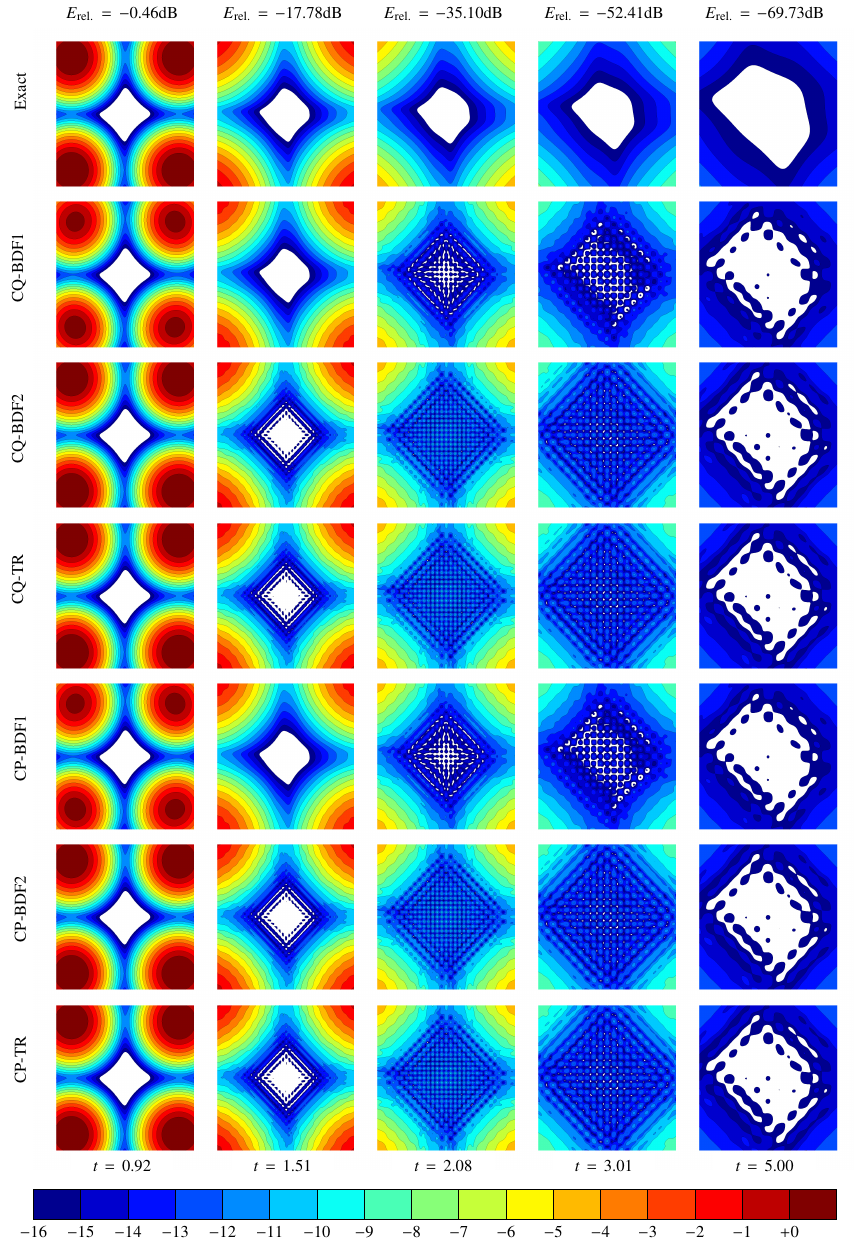}
\end{center}
\caption{\label{fig:wpcg2d_ee_typeIIB_c_12_contour}The figure shows contour 
plots of $4\log_{10}|u(\vv{x},t)|$ with various TBCs for the IIB-type 
chirped-Gaussian profile (see Table~\ref{tab:cg2d}) for $c=12$ at different 
instants of time. The numerical parameters and the labels are described in 
Sec.~\ref{sec:tests-ee}.}
\end{figure}

\begin{figure}[!htbp]
\begin{center}
\includegraphics[scale=0.95]{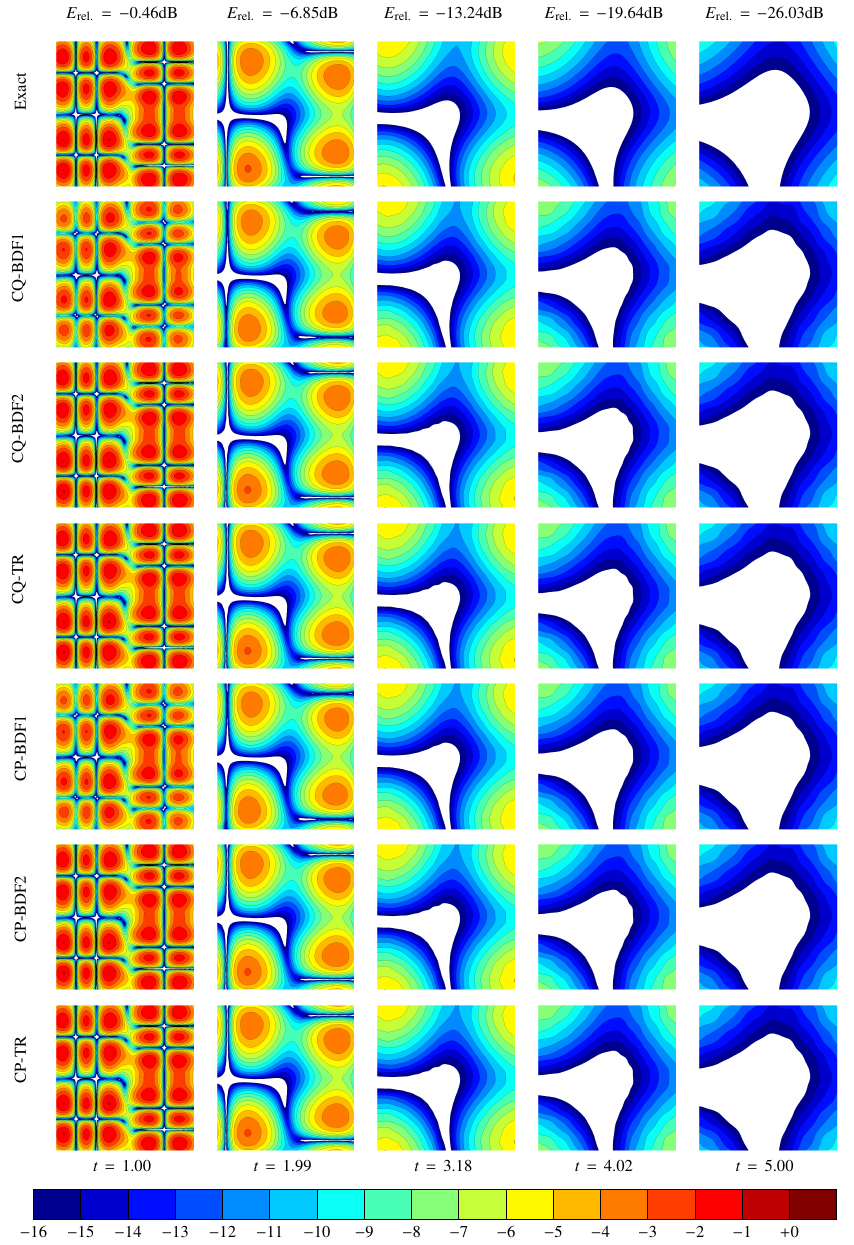}
\end{center}
\caption{\label{fig:wphg2d_ee_typeIIB_c_08_contour}The figure shows contour 
plots of $8\log_{10}|u(\vv{x},t)|$ with various TBCs for the IIB-type
Hermite-Gaussian profile (see Table~\ref{tab:hg2d}) for $c=8$ at different 
instants of time. The numerical parameters and the labels are described in 
Sec.~\ref{sec:tests-ee}.}
\end{figure}

The contour plots depicting $4\log_{10}|u(\vv{x},t)|$ and $8\log_{10}|u(\vv{x},t)|$
corresponding to the chirped-Gaussian and Hermite-Gaussian profiles are shown in 
Figs.~\ref{fig:wpcg2d_ee_typeIIB_c_12_contour} 
and~\ref{fig:wphg2d_ee_typeIIB_c_08_contour}, respectively. 
There are some visible reflections in the contour plots as expected for the maps 
obtained under high-frequency approximation. Moreover, TR based methods
better resolve the contour levels in comparison to BDF based methods.
%=============================================================================%
%=============================================================================%
\subsection{Tests for convergence}\label{sec:tests-ca}
In this section, we analyze the convergence behaviour of the numerical schemes
for the IBVP in~\eqref{eq:2D-SE-CT} where the initial condition corresponds to the 
exact solutions described in Sec.~\ref{sec:exact-solution}. The error used to 
study the convergence behaviour is quantified by 
the maximum relative $\fs{L}^2(\Omega_i)$-norm: 
\begin{equation}\label{eq:max-error-ibvp}
e=\max\left\{e(t_j)|\;j=0,1,\ldots,N_t-1\right\},
\end{equation}
for $t_j\in[0,T_{max}]$.
To analyze the convergence behaviour of NP and CP--HF
methods, we use diagonal Pad\'e approximants of order 30. The 
Table~\ref{tab:ca-params} lists the numerical parameters for studying 
the convergence behaviour of NP and HF methods. 
The complete list of the methods tested for convergence is as follows:
NP--BDF1, NP--TR, CQ--BDF1, CQ--BDF2, CQ--TR, CP30-BDF1, CP30--BDF2, CP--TR.

The numerical results for the convergence behaviour on $\Omega_i$ corresponding to the 
chirped-Gaussian and the Hermite-Gaussian profiles are shown in Fig.~\ref{fig:wpcghg2d-ca}. 
It can be seen that the diagonal Pad\'e approximant based methods, namely, NP30, CP30, 
and CQ based methods show stable behaviour for the
parameters specified in Table~\ref{tab:ca-params}. In the log-log scale, we can
clearly identify the error curve to be a straight line before it plateaus. The
orders can be recovered from the slope which we found consistent with the order
of the underlying time-stepping methods. The TR methods perform better than the 
BDF1 and BDF2 methods which is obvious from the slope of the error curves in 
Fig.~\ref{fig:wpcghg2d-ca}. 
Particularly for HF-based methods, we note that the error 
background decreases with increase in the value of $c$ for IIA-type profiles which 
are moving towards the segments of the domain while the error 
background remains unchanged with the increase in the value of $c$ for IIB-type 
profiles which are directed towards the corners of the domain. In addition, faster
moving profiles show better convergence in comparison to slow moving profiles
as evident from the error curves in Fig.~\ref{fig:wpcghg2d-ca}, particularly 
for IIA-type profiles.
Our novel Pad\'e approach based methods perform equally well for both IIA and 
IIB-type profiles.
Finally, let us observe that there is no significant difference in the 
results for convergence on $\Omega_i$ corresponding to the 
chirped-Gaussian profile and the Hermite-Gaussian profile.  

%=============================================================================%
%=============================================================================%

%=============================================================================%
%=============================================================================%

%=============================================================================%    
%=============================================================================%    
%=============================================================================%    
%=============================================================================%    
\section{Conclusion}\label{sec:conclusion}
% problem statement
In this work, we presented an efficient numerical implementation of the 
transparent boundary operator and its various approximations for the case of 
free Schr\"{o}dinger equation on a rectangular computational domain which 
presents unique challenges on account of the presence of corners.
More specifically, we addressed the spatial and temporal discretization of the 
IBVP obtained as a result of imposing BCs obtained under high-frequency 
approximation of the TBCs along with suitable 
corner conditions. The recipe was then extended to the IBVP obtained as a 
result of imposing the TBCs.

In quite a standard fashion, the use of convolution quadrature scheme for the 
fractional operators yields a numerical algorithm which scales in complexity 
(in terms of number of operations as well 
as memory storage) with the time-steps. In contrast, the use of Pad\'e approximant 
based rational approximation (in a manner which treated the corners of the boundary 
adequately) yields algorithms that do not scale in complexity with the number of 
time-steps. The performance of these algorithms were demonstrated using 
time-stepping methods upto second order of convergence. We were able to 
empirically confirm the stability and convergence behavior of the proposed
algorithms. Further, the approach developed in this work makes the system amenable to 
higher-order time-stepping methods in the class of 
linear one-step and multi-step methods\footnote{It remains to be explored if one 
can also employ the exponential integrators~\cite{BCR2000,CM2002,HL2003} which 
allow for larger step-sizes for the time-stepping scheme.}. Let us emphasize 
that for the first time, we are able to carry out a comparison of the 
high-frequency approximation of the transparent operator with the exact form 
confirming the nature of underlying approximations. We believe that these insights 
would be useful in understanding the behavior of the high-frequency approximations 
for more general systems for which exact TBC is not available. We note that, 
within the pseudo-differential approach, corners in the boundary of the 
computational domain can be adequately treated only under the high-frequency 
approximation~\cite{V2019}; therefore, we hope that the approach developed in 
this work would be relevant for the case of general Schr\"odinger equation.

Finally, let us observe that there are several open problems that have not been 
addressed in this work. For instance, the well-posedness of the IBVP under the 
nonreflecting BCs and its various approximations on the rectangular domain is left 
as an open problem. A rigorous analysis of the proposed numerical schemes for 
stability and convergence are deferred to a future publication. We note that the 
high-frequency BCs do not yield a convergent numerical scheme for the IBVP. 
However, the stability of the Pad\'e approximant based approach for the spatially 
discrete system follows from the standard results for autonomous dynamical systems. 
Extension of this approach to the TBC is left as an open problem.

%=============================================================================%    
%=============================================================================%    
%=============================================================================%    
%=============================================================================%    

\section*{Acknowledgements}
I, senior research fellow (grant no. 09/086(1431)/2019-EMR-I), am thankful to CSIR, 
India, for providing me financial assistance.  I would like to thank IIT Delhi
for providing the computational resources.

%=============================================================================%    
\bibliographystyle{elsarticle-num} 
%\bibliography{mybib}

\begin{thebibliography}{10}
\expandafter\ifx\csname url\endcsname\relax
  \def\url#1{\texttt{#1}}\fi
\expandafter\ifx\csname urlprefix\endcsname\relax\def\urlprefix{URL }\fi
\expandafter\ifx\csname href\endcsname\relax
  \def\href#1#2{#2} \def\path#1{#1}\fi

\bibitem{CiCP2008}
X.~Antoine, A.~Arnold, C.~Besse, M.~Ehrhardt, A.~Sch\"{a}dle,
  \href{https://hal.archives-ouvertes.fr/hal-00347884}{A review of transparent
  and artificial boundary conditions techniques for linear and nonlinear
  {S}chr\"{o}dinger equations}, Comm. Comput. Phys. 4~(4) (2008) 729--796.
\newline\urlprefix\url{https://hal.archives-ouvertes.fr/hal-00347884}

\bibitem{LM1998}
D.~Lee, S.~T. McDaniel, Ocean Acoustic Propagation by Finite Difference
  Methods, Modern Applied Mathematics and Computer Science, Pergmon Press, New
  York, 1988.

\bibitem{Levy2000}
M.~Levy, Parabolic equation methods for electromagnetic wave propagation,
  Vol.~45 of IEE Electromagnetic Waves Series, Institution of Electrical
  Engineers (IEE), London, 2000.

\bibitem{AK2003}
Y.~S. Kivshar, G.~P. Agrwal, Optical Solitons: From Fibers to Photonic
  Crystals, 1st Edition, Academic Press, San Diego, California, 2003.

\bibitem{S2002}
A.~Sch\"adle, Non-reflecting boundary conditions for the two-dimensional
  {S}chr\"odinger equation, Wave Motion 35~(2) (2002) 181--188.
\newblock \href {https://doi.org/10.1016/S0165-2125(01)00098-1}
  {\path{doi:10.1016/S0165-2125(01)00098-1}}.

\bibitem{HH2004}
H.~Han, Z.~Huang, \href{https://projecteuclid.org/euclid.cms/1250880210}{Exact
  artificial boundary conditions for {S}chr\"odinger equation in
  $\mathbb{R}^2$}, Comm. Math. Sci. 2~(1) (2004) 79--94.
\newline\urlprefix\url{https://projecteuclid.org/euclid.cms/1250880210}

\bibitem{S2005}
J.~Szeftel, Design of absorbing boundary conditions for {S}chr\"odinger
  equations in $\mathbb{R}^d$, SIAM Journal on Numerical Analysis 42~(4) (2005)
  1527--1551.

\bibitem{V2019}
V.~Vaibhav, On the nonreflecting boundary operators for the general two
  dimensional {S}chr\"odinger equation, J.~Math.~Phys. 60~(1) (2019) 011509.
\newblock \href {https://doi.org/10.1063/1.5030875}
  {\path{doi:10.1063/1.5030875}}.

\bibitem{AB2001}
X.~Antoine, C.~Besse, Construction, structure and asymptotic approximations of
  a microdifferential transparent boundary condition for the linear
  {S}chr\"{o}dinger equation, J. Math. Pures Appl. 80~(7) (2001) 701--738.
\newblock \href {https://doi.org/10.1016/S0021-7824(01)01213-2}
  {\path{doi:10.1016/S0021-7824(01)01213-2}}.

\bibitem{ABM2004}
X.~Antoine, C.~Besse, V.~Mouysset, Numerical schemes for the simulation of the
  two-dimensional {S}chr\"{o}dinger equation using non-reflecting boundary
  conditions, Math. Comput. 73~(248) (2004) 1779--1799.
\newblock \href {https://doi.org/10.1090/S0025-5718-04-01631-X}
  {\path{doi:10.1090/S0025-5718-04-01631-X}}.

\bibitem{ABK2012}
X.~Antoine, C.~Besse, P.~Klein, Absorbing boundary conditions for the
  two-dimensional {S}chr\"odinger equation with an exterior potential. {P}art
  {I}: {C}onstruction and a priori estimates, Math. Models Methods Appl. Sci.
  22~(10) (2012) 1250026.
\newblock \href {https://doi.org/10.1142/S0218202512500261}
  {\path{doi:10.1142/S0218202512500261}}.

\bibitem{ABK2013}
X.~Antoine, C.~Besse, P.~Klein, Absorbing boundary conditions for the
  two-dimensional {S}chr\"odinger equation with an exterior potential. {P}art
  {II}: {D}iscretization and numerical results, Numer. Math. 125~(2) (2013)
  191--223.
\newblock \href {https://doi.org/10.1007/S00211-013-0542-8}
  {\path{doi:10.1007/S00211-013-0542-8}}.

\bibitem{Menza1997}
L.~D. Menza, Transparent and absorbing boundary conditions for the
  {S}chr\"{o}dinger equation in a bounded domain, Numer. Funct. Anal. Optim.
  18~(7) (1997) 759--775.
\newblock \href {https://doi.org/10.1080/01630569708816790}
  {\path{doi:10.1080/01630569708816790}}.

\bibitem{FP2011}
R.~M. Feshchenko, A.~V. Popov, Exact transparent boundary condition for the
  parabolic equation in a rectangular computational domain, J. Opt. Soc. Am. A
  28~(3) (2011) 373--380.
\newblock \href {https://doi.org/10.1364/JOSAA.28.000373}
  {\path{doi:10.1364/JOSAA.28.000373}}.

\bibitem{YV2024}
S.~Yadav, V.~Vaibhav, Transparent boundary condition and its effectively local
  approximation for the {S}chr\"{o}dinger equation on a rectangular
  computational domain (2024).
\newblock \href {http://arxiv.org/abs/2403.07787} {\path{arXiv:2403.07787}}.

\bibitem{SV2023}
S.~Yadav, V.~Vaibhav, Nonreflecting boundary condition for the free
  {S}chr{\"o}dinger equation in 2d, in: 2023 Photonics and Electromagnetics
  Research Symposium (PIERS), 2023, pp. 328--337.
\newblock \href {https://doi.org/10.1109/PIERS59004.2023.10221299}
  {\path{doi:10.1109/PIERS59004.2023.10221299}}.

\bibitem{BCR2000}
S.~Blanes, F.~Casas, J.~Ros, Improved high order integrators based on the
  {M}agnus expansion, BIT Numer. Math. 40~(3) (2000) 434--450.
\newblock \href {https://doi.org/10.1023/A:1022311628317}
  {\path{doi:10.1023/A:1022311628317}}.

\bibitem{CM2002}
S.~M. Cox, P.~C. Matthews, Exponential time differencing for stiff systems, J.
  Comput. Phys. 176~(2) (2002) 430--455.
\newblock \href {https://doi.org/10.1006/jcph.2002.6995}
  {\path{doi:10.1006/jcph.2002.6995}}.

\bibitem{HL2003}
M.~Hochbruck, C.~Lubich, On {M}agnus integrators for time-dependent
  {S}chr\"odinger equations, SIAM J. Numer. Anal. 41~(3) (2003) 945--963.
\newblock \href {https://doi.org/10.1137/S0036142902403875}
  {\path{doi:10.1137/S0036142902403875}}.

\bibitem{CHQZ2007}
C.~Canuto, M.~Y. Hussaini, A.~Quarteroni, T.~A. Zang, {S}pectral {M}ethods:
  {F}undamentals in {S}ingle {D}omains, Springer-Verlag, Heidelberg, 2007.
\newblock \href {https://doi.org/10.1007/978-3-540-30726-6}
  {\path{doi:10.1007/978-3-540-30726-6}}.

\bibitem{Gautschi2012}
W.~Gautschi, Numerical Analysis, Birkh{\"a}user, Boston, 2012.
\newblock \href {https://doi.org/10.1007/978-0-8176-8259-0}
  {\path{doi:10.1007/978-0-8176-8259-0}}.

\bibitem{HLW2006}
E.~Hairer, C.~Lubich, G.~Wanner, Geometric Numerical Integration:
  Structure-Preserving Algorithms for Ordinary Differential Equations, 2nd
  Edition, Springer Series in Computational Mathematics, Springer-Verlag,
  Berlin, 2006.
\newblock \href {https://doi.org/10.1007/3-540-30666-8}
  {\path{doi:10.1007/3-540-30666-8}}.

\end{thebibliography}
\providecommand{\noopsort}[1]{}\providecommand{\singleletter}[1]{#1}%

\pagebreak

%=============================================================================%    
%=============================================================================%
\appendix

\section{Autonomous System}\label{app:auto-sys}
The dynamical system defined by~\eqref{eq:odes-auxi-pade} 
and~\eqref{eq:linear-sys-hf-pade} can be stated as an autonomous ODE given by
\begin{equation}\label{eq:ODE-auto}
\partial_t\vs{\mathcal{V}}=\mathcal{M}\vs{\mathcal{V}},
\end{equation}
where the vector-valued dependent variable 
$\vs{\mathcal{V}}$ comprises of blocks 
\begin{equation}
\begin{split}
&\vs{\mathcal{V}}_0=\tvec\left(\what{U}\right)=\what{\vv{U}},\\
&\vs{\mathcal{V}}_{k+1}=\tvec\left(\what{\mathcal{A}}_{k,1}\right)
=\what{\vs{\mathcal{A}}}_{k,1},\quad
\vs{\mathcal{V}}_{K+k+2}=\tvec\left(\what{\mathcal{A}}_{k,2}\right)
=\what{\vs{\mathcal{A}}}_{k,2},\quad k=0,1,\ldots,K,
\end{split}
\end{equation}
and the constant matrix $\mathcal{M}$ is defined in the following paragraphs. 
Let us label the coefficients
appearing in~\eqref{eq:linear-sys-hf-pade} as follows:
\begin{equation}
\left\{\begin{aligned}
&\gamma_{0,1}=e^{i\pi/4}\frac{1}{2}d_0\sqrt{\beta_1}\beta_2,
&&\gamma_{0,2}=e^{i\pi/4}\frac{1}{2}d_0\sqrt{\beta_2}\beta_1,&&\\
&\gamma_{k,1}=-e^{-i\pi/4}\sqrt{\beta_1}\left({b}_k+i\frac{1}{2}\beta_2d_k\right),
&&\gamma_{k,2}=-e^{-i\pi/4}\sqrt{\beta_2}\left({b}_k+i\frac{1}{2}\beta_1d_k\right),
&&k=1,\ldots,K.
\end{aligned}\right.
\end{equation}
Setting $\eta_0=0$, the dynamical system corresponding to the auxiliary 
field~\eqref{eq:odes-auxi-pade} can be stated as
\begin{equation}\label{eq:vec-aux-sys-hf}
\partial_t\what{\vs{\mathcal{A}}}_{k,1}=-\eta_k^2\what{\vs{\mathcal{A}}}_{k,1}
+\left(I_2\otimes\Lambda_1\right)\what{\vv{U}},\quad
\partial_t\what{\vs{\mathcal{A}}}_{k,2}=-\eta_k^2\what{\vs{\mathcal{A}}}_{k,2}
+\left(\Lambda_2\otimes I_1\right)\what{\vv{U}},\quad k=0,1,\ldots,K.
\end{equation}
Introducing the modified stiffness matrices 
$S'_p=S_p+e^{-i\pi/4}(b_0/\sqrt{\beta_p})\,\Lambda_p,\;p=1,2$, 
the dynamical system in~\eqref{eq:linear-sys-hf-pade} can be stated as
\begin{equation}\label{eq:vec-main-sys-hf}
\begin{split}
i \left(M_2\otimes M_1\right)\partial_t\what{\vv{U}}
&=\left[\beta_1\left(M_2\otimes S'_1\right)
+\beta_2\left(S'_2\otimes M_1\right)
+\frac{3}{4}\sqrt{\beta_1\beta_2}\left(\Lambda_2\otimes\Lambda_1\right)
\right]\what{\vv{U}}\\
&\quad+\gamma_{0,1}\left(S_2\otimes I_1\right)\what{\vs{\mathcal{A}}}_{0,1}
+\gamma_{0,2}\left(I_2\otimes S_1\right)\what{\vs{\mathcal{A}}}_{0,2}
+\sum_{k=1}^K\left[\gamma_{k,1}\left(M_2\otimes I_1\right)\what{\vs{\mathcal{A}}}_{k,1}
+\gamma_{k,2}\left(I_2\otimes M_1\right)\what{\vs{\mathcal{A}}}_{k,2}\right].
\end{split}
\end{equation}
%=============================================================================%
The form of the dynamical system~\eqref{eq:vec-aux-sys-hf} 
and~\eqref{eq:vec-main-sys-hf} dictates the matrix entries of the block matrix 
$\mathcal{M}$ which comprises non-null matrix entries on the diagonal, the first row and
the first column. It is straightforward to see that the diagonal entries of 
$\mathcal{M}$ are given by
\begin{equation}
\begin{split}
&\mathcal{M}_{0,0}
=-i\left(M_2\otimes M_1\right)^{-1}\left[\beta_1\left(M_2\otimes S'_1\right)+\beta_2\left(S'_2\otimes M_1\right)
+\frac{3}{4}\sqrt{\beta_1\beta_2}\left(\Lambda_2\otimes\Lambda_1\right)\right],\\
&\mathcal{M}_{k+1,k+1}
=-\eta_k^2\left(I_2\otimes I_1\right),\quad
\mathcal{M}_{K+k+1,K+k+1}=-\eta_k^2\left(I_2\otimes I_1\right),
\quad k=0,1,\ldots,K.
\end{split}
\end{equation}
The entries in the first row of $\mathcal{M}$ are given by
\begin{equation}
\begin{aligned}
&\mathcal{M}_{0,1}=-i\gamma_{0,1}\left(M_2\otimes M_1\right)^{-1}\left(S_2\otimes I_1\right),
&&\mathcal{M}_{0,K+2}
=-i\gamma_{0,2}\left(M_2\otimes M_1\right)^{-1}\left(I_2\otimes S_1\right), &&\\
&\mathcal{M}_{0,k+1}
=-i\gamma_{k,1}\left(M_2\otimes M_1\right)^{-1}\left(M_2\otimes I_1\right),
&&\mathcal{M}_{0,K+k+2}
=-i\gamma_{k,2}\left(M_2\otimes M_1\right)^{-1}\left(I_2\otimes M_1\right), 
&& k=1,\ldots,K.
\end{aligned}
\end{equation}
The entries in the first column of $\mathcal{M}$ are given by
\begin{equation}
\mathcal{M}_{k+1,0}=\left(I_2\otimes\Lambda_1\right),\quad
\mathcal{M}_{K+k+1,0}=\left(\Lambda_2\otimes I_1\right),\quad 
k=0,1,\ldots,K.
\end{equation}
The boundedness of the system matrix $\mathcal{M}$ under any suitable matrix 
norm follows from the standard results of Legendre-Galerkin method~\cite{CHQZ2007}.
Based on the formalism, it suffices to employ any 
zero-stable method~\cite{Gautschi2012,HLW2006} in order 
to arrive at a stable time-stepping scheme for the ODE in~\eqref{eq:ODE-auto}. A detailed 
discussion of these aspects is deffered to a future publication.
Before we conclude, let us acknowledge that a similar analysis for the novel-Pad\'e 
method can be carried out yielding a system of three coupled non-autonomous systems 
on account of the presence of three different 
temporal variables. The stability analysis of such system is not so straightforward and, therefore, 
deferred to a future publication which would also explore the possibility of 
employing exponential integrators~\cite{BCR2000,CM2002,HL2003} for such systems.

\end{document}